\documentclass[a4paper,reqno,12pt]{scrartcl}

\KOMAoption{parskip}{half-}

\usepackage [latin1]{inputenc} 

\usepackage{mathptmx} 
\usepackage[scaled=0.9]{helvet}
\usepackage{courier}

\usepackage{graphicx} 
\usepackage{float} 
\graphicspath{{files/gmpics/}} 

\usepackage {amsmath} 
\usepackage{amssymb} 
\usepackage{amscd} 
\usepackage{amsthm}

\usepackage{hyperref} 
\usepackage{color} 
\hypersetup{colorlinks=true, breaklinks=true} 

\usepackage{hyphenat}

\usepackage{csquotes} 
\usepackage[style=ieee,
						sorting=nyt,
						maxnames=4,
						backend=bibtex,
						url=false,
						isbn=false,
						doi=false,
						clearlang=true]
						{biblatex}
\bibliography{biblio}
\DeclareRedundantLanguages{english,English}{english,german,ngerman,french}

\usepackage{booktabs}
\usepackage{multirow}
\usepackage{longtable}

\DeclareMathAlphabet{\mathds}{U}{dsrom}{m}{n}
\newcommand{\C}{{\mathbb C}}

\renewcommand{\P}{{\mathbb P}}

\newcommand{\R}{{\mathbb R}}
\newcommand{\T}{{\mathbb T}}

\newcommand{\Z}{{\mathbb Z}}

\DeclareMathOperator {\ind}{ind}

\DeclareMathOperator {\val}{val}

\newcommand {\todo}[1]{\textbf{TODO #1}}

\newtheorem {theorem}{Theorem}[section]
\newtheorem {proposition}[theorem]{Proposition}
\newtheorem {lemma}[theorem]{Lemma}

\newtheorem {corollary}[theorem]{Corollary}

\newtheorem {definition}[theorem]{Definition}

\newtheorem{thm}[theorem]{Theorem}
\newtheorem{prop}[theorem]{Proposition}
\newtheorem{defn}[theorem]{Definition}
\newtheorem{problem}[theorem]{Problem}

\theoremstyle {definition}
\newtheorem {example}[theorem]{Example}

\newtheorem {remark}[theorem]{Remark}

\def \cp {\mathbb{CP}}

\def \rp {\mathbb{RP}}

\usepackage[matrix,arrow,curve]{xy}

\usepackage{longtable}

\newcommand{\ignore}[1]{\relax}
\newcommand{\conj}{\operatorname{conj}}
\newcommand{\rot}{\operatorname{rot}}
\newcommand{\St}{\operatorname{St}}
\newcommand{\Or}{\operatorname{Or}}
\newcommand{\dd}{\partial} 

\begin {document}

\title{Rational quintics in the real plane}
\author {Ilia Itenberg, Grigory Mikhalkin, Johannes Rau
\footnote{Part of the research was conducted during the stay of all three authors 
at the Max-Planck-Institut f\"ur Mathematik in Bonn.
Research is supported in part by
the FRG Collaborative Research grant 
DMS-1265228 of the U.S. National Science Foundation (I.I.), 
the grants 141329, 159240
and the NCCR SwissMAP project of the Swiss National Science Foundation (G.M.).}} 

\maketitle

\begin{abstract}
From a topological viewpoint, a rational curve in the real projective plane is
generically a smoothly immersed circle and a finite collection of isolated points. 
We give an isotopy classification of generic rational quintics in $\mathbb{RP}^2$
in the spirit of Hilbert's 16th problem.
\end{abstract}

\section{Introduction}
\subsection{Smooth real curves in the plane and Hilbert's 16th problem}
Topological classification of smooth real algebraic curves
is one of the most classical problems in real algebraic geometry.
It was included by D. Hilbert \cite{Hi} into his famous list of
problems made at the dawn of the XXth century. Let us recall
this problem as well as the basic conventions related to it.

\begin{problem}[Modern interpretation of
Hilbert's 16th problem, Part I, cf. \cite{Wilson}]
Given an integer number $d>0$,
describe possible topological types of the pair $(\rp^2,\R C)$,
where $\R C\subset \rp^2$ is the real point set of a smooth algebraic curve of degree $d$
in the real projective plane $\rp^2$. 
\end{problem}

By a real algebraic curve $C$ in $\rp^2$ we mean a real homogeneous polynomial
in $3$ variables which is considered up to multiplication by a non-zero real constant.
Such a polynomial has a zero locus $\R C \subset \rp^2$ (called the {\it real point set} of $C$) 
and a zero locus $\C C \subset \cp^2$ (called the {\it complex point set} of $C$, or {\it complexification} of $\R C$).
By abuse of language, speaking about a real algebraic curve $C$ in $\rp^2$,
we often mention only the real point set $\R C$.

The degree of a real algebraic curve $C$ in $\rp^2$ is the degree of a polynomial defining $C$.
A real algebraic curve in $\rp^2$ is called {\it non-singular} or {\it smooth},
if a polynomial defining this curve does not have critical points in $\C^3 \setminus \{0\}$.
The real point set $\R C \subset \rp^2$ of a non-singular curve $C$
is either empty or a smooth $1$-dimensional submanifold of $\rp^2$, 
that is, a disjoint
union of $l$ copies of the circle $S^1$. Furthermore, 
by Harnack's inequality \cite{Ha} we have $l\le \frac{(d-1)(d-2)}2+1$,
where $d$ is the degree of $C$.
Notice that $\frac{(d-1)(d-2)}2$ is the genus of the complexification
$\C C\subset\cp^2$. 

If $d$ is even, then every connected component $Z\subset\R C$
is homologically trivial in $\rp^2$. Such a component 
is called an {\em oval}. 
The complement $\rp^2\setminus Z$
consists of two open domains: the one homeomorphic to a disk
is called the {\em interior} of $Z$, the one homeomorphic to
a M\"obius band is called the {\em exterior} of $Z$. 

If $d$ is odd, then all but one connected components of $\R C$ 
are ovals. The remaining component $Y\subset\R C$ is isotopic to
a line $\rp^1\subset\rp^2$ and is called a {\em pseudoline}.
The complement $\rp^2\setminus Y$ is
connected, so we cannot talk of interior or exterior of a pseudoline.
Thus, for even $d$ the real point set $\R C$ of $C$ consists of $l$ ovals, 
while for odd $d$ we have one pseudoline and $l-1$ ovals. 

There is somewhat more than just the number of ovals in the topology of $(\rp^2,\R C)$. 
Two ovals $Z,Z'$ are called {\em disjoint} if their interiors
are disjoint in the set-theoretical sense.
Otherwise, they are called {\em nested}. 
More generally, we say that a collection of ovals is a {\em nest}
if any pair of ovals in this collection is nested.
The {\em depth} of a nest is the total number of ovals in the collection. 
We say that an oval $Z'$ is inside
of an oval $Z$ if $Z'$ is contained in the interior of $Z$. 
The oval is called {\em empty} if there are no ovals inside of it.

The topology of $(\rp^2,\R C)$ is completely determined by
the number of ovals together with information on each pair of ovals
whether they are disjoint or one is inside of the other.
From a combinatorial viewpoint this information is encoded
with a rooted tree whose vertices are connected components of
$\rp^2\setminus\R C$ and whose edges are the ovals of $\R C$.
The root of this tree is placed at the only non-orientable
component of  $\rp^2\setminus\R C$ if $d$ is even and
the only component adjacent to the pseudoline if $d$ is odd. 

These rooted trees, and thus the topology of $(\rp^2,\R C)$,
are traditionally encoded with the following 
system of notations introduced in \cite{Viro}.
The symbol $J$ stands for a pseudoline.
An oval is denoted with $1$.
Each rooted tree is bordered with brackets $<>$.
E.g. the topological type of a line $\rp^1\subset\rp^2$ is denoted with $<J>$,
and the topological type of an ellipse is denoted with $<1>$.

The notations are built inductively. Let $Z_1,\dots,Z_m$ be
the non-empty ovals
adjacent to the rooted component of the complement of $\rp^2\setminus\R C$.
The intersection of $\R C$ with the interior of $Z_k$ can itself
be considered as
an arrangement of ovals which is already encoded by $<x_k>$
with some symbolic notations $<x_k>$ by induction.

The topological type of $(\rp^2,\R C)$ is denoted
with
$$<J\sqcup a\sqcup 1<x_1>\sqcup\dots\sqcup 1<x_m>> \quad
\text{or}\quad
<a\sqcup 1<x_1>\sqcup\dots\sqcup 1<x_m>>,$$
depending whether $\R C$ contains a pseudoline or not
({\it i.e.}, whether $d$ is odd or even). 
Here, $a$ is the number of empty ovals 
adjacent to 
the rooted component
of $\rp^2\setminus\R C$.
The symbol $\sqcup$ is interpreted as a commutative operation,
{\it i.e.}, we do not distinguish $<x_1>\sqcup <x_2>$ 
from $<x_2>\sqcup<x_1>$. The empty curve is denoted with $<0>$. 

The following two examples were the starting
point for the classification quest known already
in the XIXth century.
\begin{example}\label{classd4}
There are $6$ topological types of non-singular curves of degree 4 in $\rp^2$. 
These are $<\alpha>$, where $0\le\alpha\le 4$, and $<1<1>>$. 
\end{example}
\begin{example}\label{classd5}
There are $8$ topological types of non-singular curves of degree 5 in $\rp^2$. 
These are $<J \ \sqcup \alpha>$, where $0\le \alpha\le 6$,
and $<J \ \sqcup 1<1>>$. 
\end{example} 

Real algebraic curves of degree $d$ in $\rp^2$ with maximal number of ovals 
allowed by Harnack's inequality, {\it i.e.},  with $\frac{(d-1)(d-2)}2+1$
connected components of the real point set, 
are called {\em $M$-curves}, cf. \cite{Petrovsky}.
Note that in degrees 4 and 5 there
are unique topological arrangements of $M$-curves:
$<4>$ and $<J\sqcup 6>$, respectively. 

\begin{defn}[F. Klein]
A real algebraic curve $C$
in $\rp^2$
is said to be {\em of type I}
if $\R C$ is null-homologous in $H_1(\C C;\Z_2)$.
\end{defn}
Note that the involution $\conj:\C C\to\C C$
of complex conjugation
has $\R C$ as its fixed point set.
Therefore, a non-singular real algebraic curve $C$ in $\rp^2$ is of type I
if and only if $\C C\setminus\R C$ is disconnected. 

It is easy to show that any $M$-curve must have type I.
Furthermore, for each degree $d$ there is a
curve of type I that has only $[\frac{d}2]$ (the integer part of $\frac{d}{2}$)
ovals.
\begin{defn}
A real curve $C$ of degree $d$ in $\rp^2$ is called {\em hyperbolic} if $\R C$ has
a nest of depth $[\frac{d}2]$. 
\end{defn}
The B\'ezout theorem implies that the topological arrangement 
of the hyperbolic curve is unique: there are no other
ovals except for those from the nest of depth $[\frac{d}2]$.
Indeed, if there is any other oval, then we may draw a straight line 
through that oval and the innermost oval in the nest. 
Such a line would intersect all ovals from the nest
and the additional oval at least in two points each.
This gives $2+2[\frac{d}2]>d$ points of intersection
between a curve of degree $d$ and a line.
In particular, $[\frac{d}2]$ is the maximal possible depth
of a nest for a curve of degree $d$ in $\rp^2$. 

It can be shown that all hyperbolic curves are of type I.
The complex orientation formula \cite{Rokhlin}
(which is reviewed in Section \ref{complex_orientations}) 
implies the following classical statement
which was known already to Klein. 

\begin{thm}
There are two possible topological
types of 
non-singular curves of degree 4 and type I in $\rp^2$:
$<4>$ {\rm (}the $M$-quartic{\rm )} and $<1<1>>$ {\rm (}the hyperbolic quartic{\rm )}. 

There are three possible topological
types 
of non-singular curves of degree 5 and type I in $\rp^2$:
$<J\sqcup 6>$ {\rm (}the M-quintic{\rm )}, $<J\sqcup 1<1>>$
{\rm (}the hyperbolic quintic{\rm )}, and $<J\sqcup 4>$
{\rm (}the four-oval quintic of type I{\rm )}. 
\end{thm} 

Currently the classification of topological arrangements
of non-singular curves of degree $d$ in $\rp^2$ is
known up to degree $7$ (Viro \cite{Viro-deg7}). As this classification
is rather large, below we list only the possible topological
types of $M$-curves. Note that
for $d\ge 6$ the topological type of an $M$-curve
is no longer unique.
\begin{thm}[Gudkov]
There are three possible topological types of
non-singular $M$-curves of degree $6$ in $\rp^2$:
$<9\sqcup 1<1>>$ {\rm (}the so-called {\em Harnack sextic}{\rm )},
$<1\sqcup 1<9>>$ {\rm (}the so-called {\em Hilbert sextic}{\rm )}, and
$<5\sqcup 1<5>>$ {\rm (}the so-called {\em Gudkov sextic}{\rm )}. 
\end{thm} 

\begin{thm}[Viro]
There are fourteen possible topological types of
non-singular $M$-curves of degree $7$ in $\rp^2$:
$<J\sqcup 15>$ and $<J\sqcup \alpha\sqcup 1<14-\alpha>>$, where $1\le\alpha\le 13$. 
\end{thm} 

\subsection{Generic rational curves in the plane}\label{sec-rational}
This paper is mainly devoted to {\it rational curves} in $\rp^2$. 
The complex point set of such a real rational curve of degree $d$ 
can be described as the image of the map $\varphi: \cp^1 \to \cp^2$ defined by 
$$
(z_0 : z_1) \mapsto (P(z_0, z_1) : Q(z_0, z_1) : R(z_0, z_1)),
$$
where $P$, $Q$, $R$ are real homogeneous polynomials of degree $d$
which do not have common zeros in $\cp^1$.

A generic rational curve in $\rp^2$ is {\it nodal}, which means
that the only possible singular points of the curve are {\it non-degenerate double points} (also called {\it nodes}). 


We may distinguish three types of nodes of a nodal curve $C$ in $\rp^2$: 
{\em hyperbolic}, {\em elliptic} and
{\em imaginary nodes}. 
Hyperbolic nodes are formed by intersections of pairs of real branches
of $\R C$.
These are points of $\rp^2$ such that $C$ is given by $x^2-y^2=0$
in some local coordinates $(x,y)$ near these points.
Elliptic nodes are formed by real intersections of pairs of
complex conjugated branches of $\C C$.
These are points of $\rp^2$ such that $C$ is given by $x^2+y^2=0$
in some local coordinates $(x,y)$ near these points.
Finally, imaginary nodes are nodes of $C$ 
in $\cp^2\setminus\rp^2$.
Such points come in pairs of complex conjugate points.
We denote the number of hyperbolic (respectively, elliptic, imaginary) nodes
by $h$ (respectively, $e$, $c$). 

The following proposition is straightforward.
\begin{prop}
The real point set of any nodal rational curve 
in $\rp^2$ is the disjoint union of
a circle 
generically immersed in $\rp^2$ and a finite set of 
elliptic nodes.
\end{prop} 

Similarly to Hilbert's 16th problem, one can ask for a topological classification
of pairs $(\rp^2, \R C)$, where $C$ is a nodal rational curve of a given degree $d$ in $\rp^2$. 
Since any self-homeomorphism of $\rp^2$ is isotopic to the identity,
such a topological classification of pairs $(\rp^2, \R C)$ provides an {\it isotopy classification}
of the real point sets $\R C$ of nodal rational curves of a given degree in $\rp^2$.
The isotopy classification in question is known up to degree $4$
(see \cite{DMello} and Section \ref{sec:degree4} for details concerning the classification for degree $4$). 
Among other results related to the isotopy classification of real rational curves,
one can mention, for example, 
the study of maximally inflected real rational curves in the context of the Shapiro-Shapiro conjecture 
and the real 
Schubert calculus; see \cite{KhSo}. 
In this paper we study nodal rational curves of degree $5$ in $\rp^2$. 

\subsection{Classification of generic rational curves of degree \texorpdfstring{$5$}{5}}\label{main-classification}

In this section we state our main result, namely the isotopy classification of nodal rational curves of degree $5$ 
in $\rp^2$. 
The classification is presented below in the form of the list of {\it smoothing diagrams} (see Section \ref{backtosmoothovals} 
for the precise definitions) 
of the curves under consideration. Each smoothing diagram describes an isotopy type
of a 
nodal rational curve of degree $5$ in $\rp^2$. The isotopy type is obtained by contracting the {\it vanishing cycles} 
(that is, edges) 
of a smoothing diagram creating a hyperbolic node for each vanishing cycle. 

\begin{theorem} \label{maintheorem}
The isotopy types of nodal rational curves of degree $5$ 
in $\rp^2$ are exactly those listed in the tables 
\ref{list1} -- \ref{list6}. 
\end{theorem}

\begin{remark}\label{mainremark}
In tables \ref{list2} and \ref{list4} we sometimes incorporated several smoothing diagrams in
the same picture. 
The merged smoothing diagrams only differ by the attachment of a single oval by a
single edge {\rm (}drawn with dashed lines{\rm )}.  
The factor next to these pictures {\rm (}e.g.\ $2${\rm x)} indicates how many smoothing 
diagrams are merged.

There are exactly $121$ isotopy types of nodal rational curves of degree $5$ in $\rp^2$. 
	With regard to $l,h,e,c$, 
	we get the following numbers of isotopy types: 
	
	\begin{center}
	$c=0$ \hspace{7ex}
	\begin{tabular}{@{}lllllllll@{}}
		\toprule
		$l$    & $e$ & 0  & 1  & 2  & 3 & 4 & 5 & 6 \\ \cmidrule(l){2-9} 
		3      &     & 13 & 3  &    &   &   &   &   \\
		5      &     & 24 & 12 & 4  & 1 &   &   &   \\
		7      &     & 9  & 9  & 6  & 4 & 2 & 1 & 1 \\ \cmidrule(l){2-9} 
		$\sum$ &     & 46 & 24 & 10 & 5 & 2 & 1 & 1 \\ \bottomrule
	\end{tabular}		
	
	$c=2$ \hspace{7ex}
	\begin{tabular}{@{}lllllll@{}}
		\toprule
		$l$    & $e$ & 0  & 1 & 2 & 3 & 4 \\ \cmidrule(l){2-7} 
		3      &     & 7  & 2 &   &   &   \\
		5      &     & 8  & 4 & 2 & 1 & 1 \\ \cmidrule(l){2-7} 
		$\sum$ &     & 15 & 6 & 2 & 1 & 1 \\ \bottomrule
	\end{tabular}

	$c=4$ \hspace{7ex}
	\begin{tabular}{@{}lllll@{}}
		\toprule
		$l$ & $e$ & 0 & 1 & 2 \\ \cmidrule(l){2-5} 
		3   &     & 3 & 2 & 1 \\ \bottomrule
	\end{tabular}
	\end{center} 

   Additionally, we have exactly one isotopy type for $c=6$ {\rm (}a non-contractible loop in $\rp^2${\rm )}. 
\end{remark}

\begin{table}%
\centering
\input{files/pics/nested1a.TpX}
\input{files/pics/nested01a.TpX}
\input{files/pics/nested001a.TpX}
\input{files/pics/nested001b.TpX}

\input{files/pics/nested111a.TpX}
\input{files/pics/nested0001a.TpX}
\input{files/pics/nested0111a.TpX}
\input{files/pics/nested0111b.TpX}

\input{files/pics/nested00001a.TpX}
\input{files/pics/nested00001b.TpX}
\input{files/pics/nested01011a.TpX}
\input{files/pics/nested01011b.TpX}

\input{files/pics/nested11111a.TpX}

\input{files/pics/nested00001ell.TpX}
\input{files/pics/nested01011ell.TpX}
\input{files/pics/nested11111ell.TpX}

	\caption{Smoothing diagrams of all isotopy types with $c=0, l=3$}%
	\label{list1}%
\end{table}

\begin{table}%
\centering
\input{files/pics/4oval300+.TpX}
\input{files/pics/4oval310a+.TpX}
\input{files/pics/4oval310a++.TpX}
\input{files/pics/4oval310b+.TpX}
\input{files/pics/4oval310b++.TpX}

\input{files/pics/4oval311c+.TpX}
\input{files/pics/4oval311e+.TpX}
\input{files/pics/4oval311f+.TpX}
\input{files/pics/4oval21A.TpX}
\input{files/pics/4oval21B.TpX}

\input{files/pics/4oval21Ca.TpX}
\input{files/pics/4oval21Cb+.TpX}
\input{files/pics/4oval21Cc+.TpX}
\input{files/pics/4oval111.TpX}

\input{files/pics/4oval300ell.TpX}
\input{files/pics/4oval310aell.TpX}
\input{files/pics/4oval310bell.TpX}
\input{files/pics/4oval311cell.TpX}
\input{files/pics/4oval311eell.TpX}

\input{files/pics/4oval311fell.TpX}
\input{files/pics/4oval21Aell.TpX}
\input{files/pics/4oval22Cell.TpX}

\input{files/pics/4oval300ell2.TpX}
\input{files/pics/4oval310ell2.TpX}
\input{files/pics/4oval22ell2.TpX}
\input{files/pics/4oval300ell3.TpX}

	\caption{Smoothing diagrams of all isotopy types with $c=0, l=5$}%
  \label{list2}%
\end{table}

\begin{table}%
\centering
  \input{files/pics/harnacka.TpX}
  \input{files/pics/harnackb.TpX}
  \input{files/pics/harnackc.TpX}
  \input{files/pics/harnackd.TpX}
  \input{files/pics/harnacke.TpX}
	
  \input{files/pics/harnackf.TpX}
  \input{files/pics/harnackg.TpX}
  \input{files/pics/harnackh.TpX}
  \input{files/pics/harnacki.TpX}
	
\input{files/pics/harnackella.TpX}
\input{files/pics/harnackellb.TpX}
\input{files/pics/harnackellc.TpX}
\input{files/pics/harnackelld.TpX}
\input{files/pics/harnackelle.TpX}

\input{files/pics/harnackellf.TpX}
\input{files/pics/harnackellg.TpX}
\input{files/pics/harnackellh.TpX}
\input{files/pics/harnackelli.TpX}

\input{files/pics/harnackellj.TpX}
\input{files/pics/harnackellk.TpX}
\input{files/pics/harnackelll.TpX}
\input{files/pics/harnackellm.TpX}
\input{files/pics/harnackelln.TpX}

\input{files/pics/harnackello.TpX}

\input{files/pics/harnackellp.TpX}
\input{files/pics/harnackellq.TpX}
\input{files/pics/harnackellr.TpX}
\input{files/pics/harnackells.TpX}

\input{files/pics/harnackellt.TpX}
\input{files/pics/harnackellu.TpX}
\input{files/pics/harnackellv.TpX}
\input{files/pics/harnackellw.TpX}

	\caption{Smoothing diagrams of all isotopy types with $c=0, l=7$}%
  \label{list3}%
\end{table}

\begin{table}%
\centering
\input{files/pics/cc1o.TpX}
\input{files/pics/cc1p.TpX}
\input{files/pics/cc1q.TpX}
\input{files/pics/cc1r.TpX}

\input{files/pics/cc1lnew.TpX}
\input{files/pics/cc1n.TpX}

\input{files/pics/cc1v.TpX}
\input{files/pics/cc1vv.TpX}
	
	\caption{Smoothing diagrams of all isotopy types with $c=2, l=3$}%
  \label{list4}%
\end{table}

\begin{table}%
\centering
\input{files/pics/cc1i.TpX}
\input{files/pics/cc1c.TpX}
\input{files/pics/cc1f.TpX}
\input{files/pics/cc1g.TpX}

\input{files/pics/cc1a.TpX}
\input{files/pics/cc1b.TpX}
\input{files/pics/cc1d.TpX}
\input{files/pics/cc1e.TpX}

\input{files/pics/cc1s.TpX}
\input{files/pics/cc1t.TpX}
\input{files/pics/cc1u.TpX}
\input{files/pics/cc1uu.TpX}

\input{files/pics/cc1w.TpX}
\input{files/pics/cc1x.TpX}
\input{files/pics/cc1y.TpX}
\input{files/pics/cc1z.TpX}
	
	\caption{Smoothing diagrams of all isotopy types with $c=2, l=5$}%
  \label{list5}%
\end{table}

\begin{table}%
\centering
\input{files/pics/cc2a.TpX}
\input{files/pics/cc2b.TpX}
\input{files/pics/cc2ccc.TpX}

\input{files/pics/cc2c.TpX}
\input{files/pics/cc2cc.TpX}
\input{files/pics/cc2d.TpX}
\input{files/pics/cc3.TpX}

	\caption{Smoothing diagrams of all isotopy types with $c=4$ (and hence $l=3$)
	and of the unique isotopy type with $c=6$ (and hence $h=e=0, l=1$)}%
  \label{list6}%
\end{table}

The parameter $l$ in the above list is the number of connected components of the real part $\R C_\circ$
of an appropriate small perturbation $C_\circ$ of a nodal rational curve $C$ of degree $5$ in $\rp^2$ 
(see Section \ref{possible_schemes}). 
The proof of Theorem \ref{maintheorem} is presented in Section \ref{sec:restrictions} 
(restrictions on the topology of nodal rational curves of degree $5$ in $\rp^2$) 
and Section \ref{sec:constructions} (constructions).

\section{Generically immersed curves and their smoothings}\label{new-section}

\subsection{Generically immersed curves instead of
smoothly embedded curves}
In contrast with smooth curves in the plane,
a connected component of an immersed curve may
be quite complicated topologically.
In the space of all immersions of a circle
to the plane ({\it i.e.}, differentiable maps from $S^1$ to the plane such that
the differential never vanishes)
we may distinguish {\em generic immersions} (see below)
that only have transverse double points as its self-intersections.
This is the only type of singularity of an immersed curve
that survives under all small perturbations in the class of smooth maps.

It was noted by V. Arnold \cite{A} that generic immersions of
a circle into a plane have some common behaviour with knots
in a 3-space, particularly from the viewpoint of finite-type invariants.
In conventional knot theory a knot is an embedding of a circle $K\approx S^1$
to the 3-space $\R^3$. Such a knot is commonly depicted with the help
of a linear projection $\pi:\R^3\to\R^2$ (normally referred to as
a {\em vertical projection}).
The image $\pi(K)$ is immersed to $\R^2$ and all of its self-crossing
points are non-degenerate double points (often called {\em crossings} in this context)
if the projection $\pi$ is chosen generically.
Once we specify at every crossing which of the two branches
is above and which is below we get a presentation of a knot by
the so-called {\em knot diagram}.
Different knot diagrams may give the same knot if they are connected with
a sequence of the 
Reidemeister moves for knot diagrams.

\begin{defn}
Let $X$ be a smooth surface.
An immersion $i: S^1\to X$ is called {\em generic} if
all its self-crossing points
are
non-degenerate double points {\rm (}called \emph{nodes} here{\rm )}.
Here $S^1$ is an oriented circle
{\rm (}e.g. the unit circle in $\C$ oriented counterclockwise{\rm )}.
\end{defn}

A knot diagram, after forgetting which branch is above and which is below at the nodes,
is an example of a generic immersion of a circle in $\R^2$.

We consider two generic immersions to be equivalent if they are homotopic
in the class of generic immersions.
Obviously, such equivalence classes are uniquely determined by the isotopy
type of $i(S^1) \subset X$.
By abuse of language,
we often identify a generic immersion $i: S^1 \to X$ with its image $K=i(S^1)$. 

Let $T_1(X)$ denote the unit tangent bundle of $X$.
Using the standard orientation of $S^1$
we may lift the immersion $K\subset\ X$ to
$\tilde{K}\subset T_1(X)$ by associating to each point of
the circle $S^1$ parameterizing $K$
the unit tangent vector to $K$ according to the orientation.
Thus, $\tilde K$ is the image of the Gauss-type map
$\tilde i:S^1\to T_1(X)$.
\begin{defn}
The homology class $[\tilde K]\in H_1(T_1(X);\Z)$ is called
the {\em rotation number} of $K$ and is denoted with $\rot(K)$.
\end{defn}
\begin{example}
Let us fix an orientation for $\R^2$ and let $K$ be a positively oriented embedded circle in $\R^2$. 
We have $H_1(T_1(\R^2);\Z)\cong\Z$, and we fix the isomorphism by setting $\rot(K) = +1$.
\end{example}

\begin{example}\label{rot-rp2}
We have $H_1(T_1(\rp^2);\Z)\cong\Z_4$.
We fix the isomorphism by
setting $\rot(\rp^1)=1\in \Z_4$, where $\rp^1 \subset \rp^2$ is
a line.
Note that there is no need to specify the orientation
of this line as the two choices of orientation are isotopic.
Thus, if $X=\rp^2$, then $\rot(K)\in \Z_4$. Furthermore, it is easy to see that
if $[K]=0\in H_1(\rp^2;\Z) = \Z_2$, then $\rot(K)$ is even,
while if $[K]\neq 0\in H_1(\rp^2;\Z) = \Z_2$, then $\rot(K)$ is odd.
\end{example}

The following
statement is classical.
\begin{thm}[Whitney \cite{Whitney}]
\label{whitney}
Two immersions are homotopic in the class of
{\rm (}not necessarily generic{\rm )} immersions
if and only if their rotation numbers coincide.
\end{thm} 

It is easy to see that if two
generic immersions are homotopic in the class of
(not necessarily generic) immersions, then they are obtained
from each other
by a series of the following planar immersion
counterparts of two of the knot theory
Reidemeister moves: namely, the second and the third Reidemeister moves.

\begin{figure}[ht]
\centerline{
\includegraphics[height=25mm]{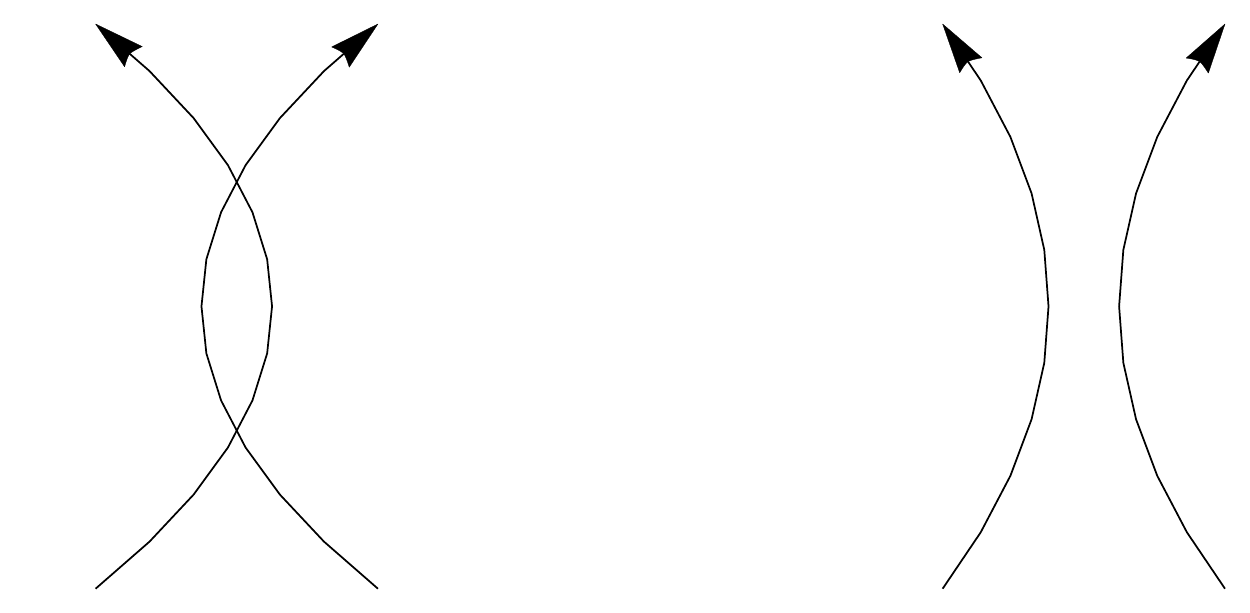}
}
\caption{Direct self-tangency perestroika}
\label{j+}
\end{figure}

\begin{figure}[ht]
\centerline{
\includegraphics[height=25mm]{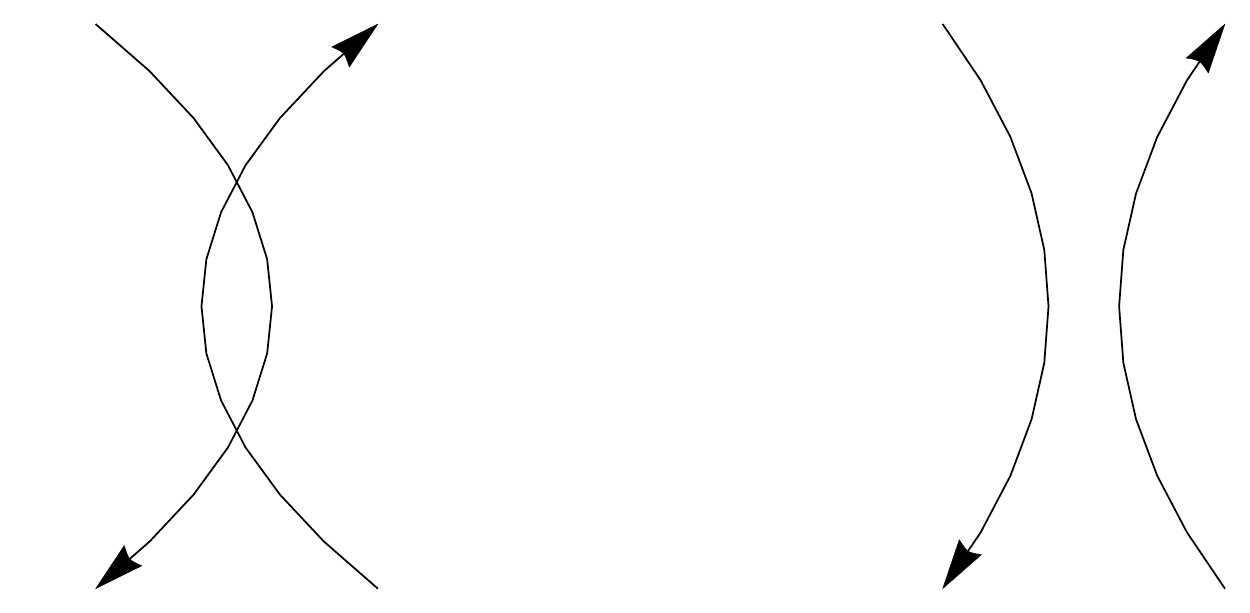}
}
\caption{Inverse self-tangency perestroika}
\label{j-}
\end{figure}

The second Reidemeister move
corresponds to passing through
a generic double tangency point. Here we distinguish
two cases: when the orientations of the two tangent
branches agree and when they disagree.
The first case is called the {\em direct self-tangency perestroika},
see Figure \ref{j+},
while the second one is called the {\em inverse self-tangency perestroika},
see Figure \ref{j-},
cf. \cite{A}, \cite{Viro-immersions}.

\begin{figure}[ht]
\centerline{
\includegraphics[height=25mm]{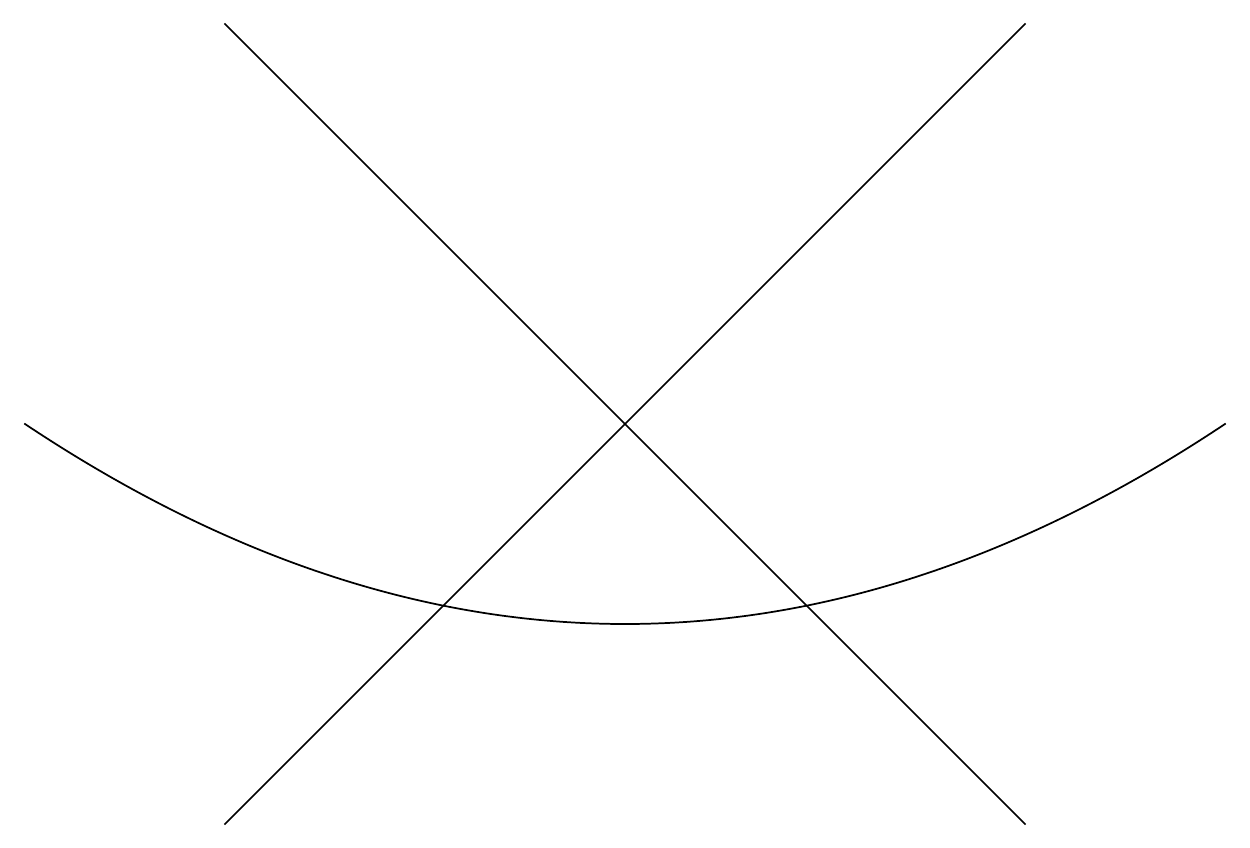}\hspace{45pt}
\includegraphics[height=25mm]{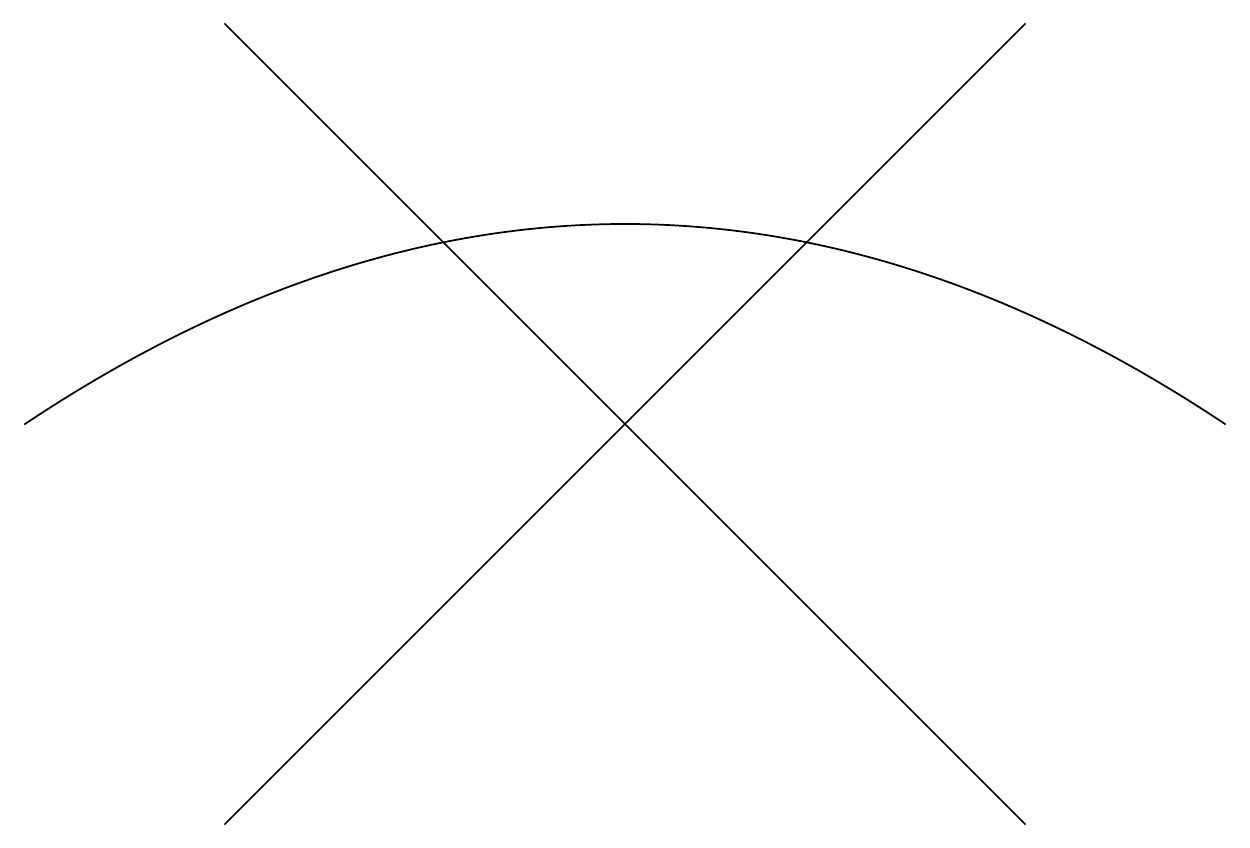}
}
\caption{Triple point perestroika}
\label{reidemeister3}
\end{figure}

The move from Figure \ref{reidemeister3} corresponds to
passing through a triple point.
Such a move is called the {\em triple point perestroika}.

Arnold \cite{A} has shown that in the case $X=\R^2$
there are three {\em degree 1} (in Vassiliev's sense)
invariants corresponding to these moves.
We do not specify the orientations at Figure \ref{reidemeister3}
as all such choices combine to a single degree 1 invariant.

The invariant corresponding to the direct self-tangency perestroika
is called $J_+$. It increases by 2 if such a move is performed in
the direction when the number of nodes is increased by 2
and remains invariant if we perform either an inverse self-tangency
or a triple point perestroika. (If we perform the direct self-tangency
move in the opposite direction, then $J_+$ decreases by 2 accordingly.)

Similarly, $J_-$ is increased by 2 if an inverse self-tangency
perestroika is performed in the direction when the number of nodes
is decreased by 2
and remains invariant if we perform a direct self-tangency
or a triple point perestroika. 

Furthermore, there is a consistent choice of direction
for performing the triple point perestroika. This choice
allows us to declare that the third invariant (called
the {\em strangeness} $St$) increases by 1 when the triple point
perestroika is performed in this direction and does not
change when any of the self-tangency perestroikas is performed.
The rule for specifying this direction
is indeed quite strange, although it was clarified by Shumakovich
in \cite{Shu} with the help of an explicit formula.

It was shown in \cite{A} that if we start from a generic
immersion, perform a number of the moves discussed above
and return to the same generic immersion, then the total
increment for each of the numbers $J_+$, $J_-$ and $St$ is zero.
Thus, to turn $J_+$, $J_-$ and $St$ to conventional integer-valued
invariants it is sufficient to choose their normalization
on one generic immersion for each possible $\rot(K)$.
This was done in \cite{A} (in such a way that the resulting invariants
are additive with respect to the connected sum).

From the definition, $J_+-J_-$ equals to the number of nodes
of the generic immersion. It was noted by Viro \cite{Viro-immersions}
that $J_-$ can be easily computed from the complex
orientation formula.
In the same paper, Viro gave
a well-defined adaptation of $J_-$ for immersed curves in $\rp^2$,
the situation we consider in this paper.


\subsection{Back to smooth ovals: smoothing of an immersion} \label{backtosmoothovals}
Let $K\subset\rp^2$ be an immersion to $\rp^2$ of a disjoint
union of a collection of oriented circles. 
By Example \ref{rot-rp2} we have $\rot(K)\in\Z_4$
(if $K$ is multicomponent, then its rotation number
is the sum of the rotation numbers of the components of $K$),
while the parity of $\rot(K)$ is determined by
the parity of the homological class $[K]\in H_1(\rp^2;\Z) = \Z_2$
which we call {\em the degree} $d=d(K)$ of $K$.

\begin{figure}[h]
\centerline{
\includegraphics[height=25mm]{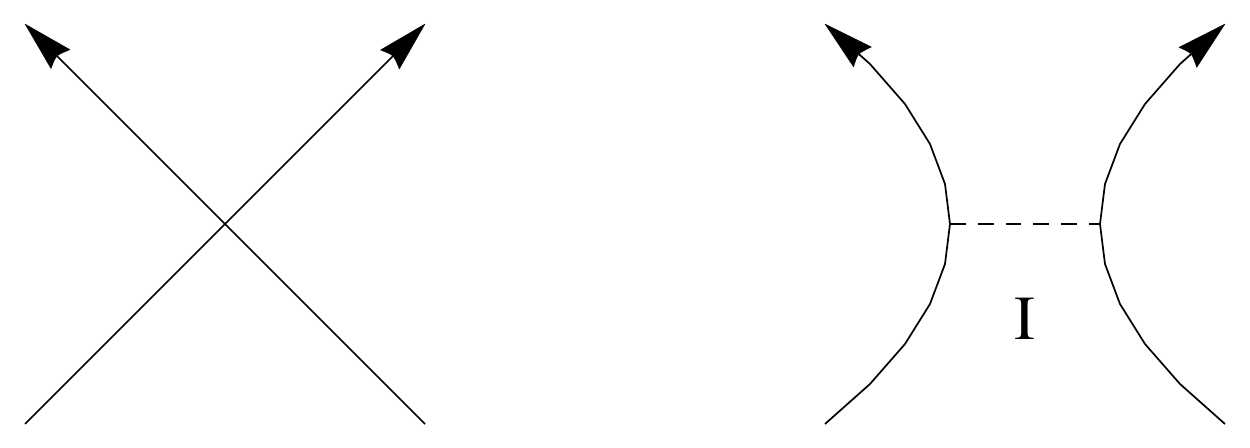}
}
\caption{Smoothing according to the orientation and
appearance of the vanishing cycle $I$}
\label{or-smoothing}
\end{figure}

Assume that the immersed curve $K$ is generic
(as before, it means that all self-crossing points are nodes).
Let $n=n(K)$ be the number of nodes of $K$, 
and let $K_\circ\subset\rp^2$ be the collection of 
smoothly embedded oriented circles
(defined up to isotopy) obtained from $K$ by smoothing
every node of $K$ as shown on
Figure \ref{or-smoothing}.
Note that to each node $p\in K$ we associate a 
closed embedded path $I_p \subset \rp^2$ with two endpoints
on the corresponding two arcs of the smoothing
and not intersecting $K_\circ$ at
inner points of the path.
\begin{defn} \label{smoothingdiagram}
The oriented curve $K_\circ$ is called {\em smoothing} of $K$.
The intervals $I_p$ are called {\em vanishing cycles}
corresponding to nodes $p$.
The diagram
\begin{equation}\label{DeltaK}
\Delta_K=(K_\circ;\bigcup\limits_p I_p),
\end{equation}
where $p$
runs over all nodes of $K$, is called
the {\em smoothing diagram} of $K$.
\end{defn}
Note that all vanishing cycles $I_p$ in the diagram $\Delta_K$
are disjoint.

\begin{defn}\label{def-coherent-membrane}
Let $L\subset\rp^2$ be
a collection of disjoint oriented embedded circles.
An {\em $L$-membrane} $I_p\subset\rp^2$ is
a smoothly embedded interval such that
$I_p\cap L$ coincides with the two endpoints of $I_p$
and at these endpoints the interval $I_p$ is transverse to $L$.

We say that an $L$-membrane $I_p$ is
{\em coherent} {\rm (}with respect to the orientation on $K_\circ${\rm )}
if there exists a deformation of $I_p$ in the class of
$L$-membranes disjoint from $I_p$ so that the endpoints
move according to the orientations.
In other words, $I_p$ is coherent if the orientations
of $L$ are as shown on Figure \ref{or-smoothing}.

An {\rm (}abstract{\rm )} smoothing diagram $\Delta=(L;I)$
consists of a collection $L \subset \rp^2$ of disjoint oriented embedded circles
and a collection $I=\bigcup_p I_p$ of disjoint embedded
closed intervals $I_p\subset\rp^2$ so that
the number of intervals $I_p$ is finite and each interval $I_p$ is a coherent
$L$-membrane. We consider such diagrams $(L;I)$ and $(L';I')$ equivalent if
they are isotopic
{\rm (}meaning the existence of an isotopy identifying
$L$ and $L'$ as well as $I$ and $I'${\rm )}.
\end{defn}

\begin{prop}\label{prop-DeltaK}
For every smoothing diagram $\Delta$
there exists a generic immersion $K$ of a disjoint
union of oriented circles to $\rp^2$ such that
$\Delta=\Delta_K$. Moreover, $K$ is unique up to isotopy.
\end{prop}
The proof of this proposition is straightforward.
We collapse every coherent membrane by performing
the move opposite to the smoothing depicted at Figure \ref{or-smoothing}.

Let $K \subset \rp^2$ be a generic immersion of a disjoint union of oriented circles.
We have a well defined lift $\tilde{K_\circ}\in T_1(\rp^2)$
as well as its homology class $\rot(K_\circ)=
[\tilde{K}_\circ]\in H_1(T_1(\rp^2))$.
Denote with $l(K)$ the number of connected components of $K_\circ$
and with $k$ the number of circles in the immersion $K$
({\it i.e.}, the number of connected components of $\tilde K$).

\begin{lemma}\label{parity-Kcirc}
We have $$\rot(K_\circ)=\rot(K).$$
Furthermore,
$$l(K)\le n(K)+k,$$ $$l(K) = n(K)+k \mod 2.$$
\end{lemma}
\begin{proof}
The homology class in $H_1(T_1(\rp^2)$ is preserved under each smoothing.
Let us consecutively smooth nodes as in Figure \ref{or-smoothing},
one by one.
At each 
step we increase
or decrease the number of components of $\tilde K$ by 1, depending on whether
the two branches of the node belong to the same or distinct
components.
\end{proof} 

The following statement can be viewed as
the $\rp^2$-version of the Whitney
formula \cite{Whitney}. It determines $\rot(K)$ once we know
the degree of $K$ and $n(K)$.
\begin{lemma}
If the degree $d(K)$ of $K \subset \rp^2$ is even, then
$$\rot(K)=2n(K)+2k\in\Z_4.$$
If
$d(K)$ is odd, then $$\rot(K)=2n(K)+2k-1\in\Z_4.$$
\end{lemma}

\begin{proof}
When smoothing, at each step we change the parity of the number of
components of $\tilde{K}$ by 1, but preserve its homology
class $[\tilde K]=\rot(K)$.
Once all nodes are smoothed we have $a$ even degree
components and $b$ odd degree components.
Each such even degree component is isotopic to a small circle
and thus its rotation number is 2.
Each odd degree component is isotopic to a line and
thus its rotation number is 1.

We have $a + b = n(K) + k \mod 2$,
$b = d(K) \mod 2$, and  $b \le 1$,
while $\rot(K) = 2a + b \mod 4$ which implies the statement of the lemma.
\end{proof}

Theorem \ref{whitney} implies the following statement.
\begin{corollary}
Two generic immersions of a circle to $\rp^2$
with the same parity of the number $n$ of nodes and the same degree $d$
are homotopic in the class of immersions.
\end{corollary} 


Let $K\subset\rp^2$ be a generic immersion of a
disjoint union of oriented circles,
and let $p\in\rp^2\setminus K$ be a point.
We choose an isomorphism between the group $H_1(\rp^2\setminus\{p\};\Z)$
and $\frac{1}{2} \Z$.
The only ambiguity in the choice of this isomorphism
is the sign.
The index $\pm\ind_K(p)$ is the half-integer well-defined up to sign
given by $[K] \in H_1(\rp^2\setminus\{p\};\Z)$.
Using half-integers guarantees that the index jumps by one
each time we pass a branch of $K$.
The absolute value $|\ind_K(p)|$ and the square $\ind^2_K(p)$
are well-defined half- resp.\ quarter-integer numbers.
Hence $|\ind_K|$ and $\ind^2_K$
are locally constant functions on $\rp^2\setminus K$.

\begin{figure}[!ht]
\centerline{
\input{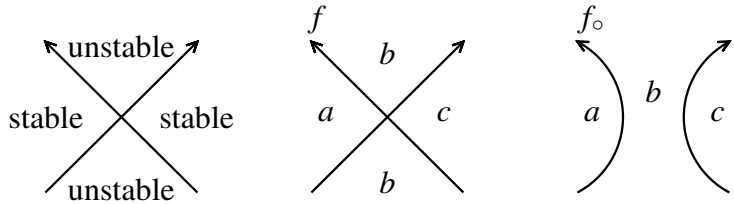}
}
\caption{(Un-)stable quadrants and a smoothable function $f$}
\label{smoothablefunction}
\end{figure}

Given a small neighbourhood $U$ of a node $p$ of $K$,
the set $U \setminus K$ consists of four connected components,
called \emph{quadrants} here. When smoothing $K$ to $K_\circ$,
the two opposite quadrants which stay disconnected are called \emph{stable},
the other two quadrants which get connected are called \emph{unstable}.
A locally constant function $f:\rp^2\setminus K\to\R$
which at each node $p$ of $K$ takes the same value
on the two unstable quadrants is called \emph{smoothable}.
The functions $|\ind_K|$ and $\ind^2_K$ are of this form.
Obviously, in this case $f$ descends to a unique function
$f_\circ:\rp^2\setminus K_\circ\to\R$
(see Figure \ref{smoothablefunction}).

\begin{figure}[!ht]
\centerline{
\input{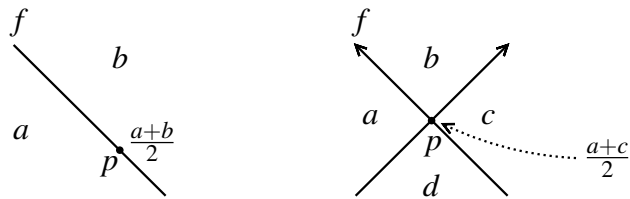}
}
\caption{Extension of a function $f$ to $\rp^2$}
\label{extendfunction}
\end{figure}

Let $f:\rp^2\setminus K\to\R$ be a locally constant function.
We extend $f$ to the whole plane $\rp^2$ as follows.
Suppose that $p\in K$. If $p$ is not one of the
nodes of $K$, we define $f(p)$ as the average of the values
of $f$ at the two regions of $\rp^2\setminus K$ adjacent to $p$.
If $p\in K$ is a node, we define $f(p)$ to be the average of the values of
$f$ on the two stable quadrants (see Figure \ref{extendfunction}).
Note that when we apply this construction to $|\ind_K|$ and $\ind^2_K$,
in general we get $(|\ind_K|(p))^2 \neq \ind^2_K(p)$ for $p \in K$.

Following \cite{Viro-chi} we define the integral $\int\limits_{\rp^2}fd\chi$ with respect to Euler characteristic
for any function
 $f:\rp^2\to\R$ with the finite image $f(\rp^2)\subset\R$
and such that there exists a cellular decomposition
of $\rp^2$, where the inverse image
$f^{-1}(j)$ is a finite union of open cells for any $j \in \R$.
Namely, if a set $A\subset\rp^2$ is a finite union of open cells, we set
$\chi(A)\in\Z$ to be the difference of the number of even-dimensional
cells in $A$ with the number of odd-dimensional cells in $A$. This is nothing
else but the Euler characteristic of $A$, but if we are to
compute it homologically, we have to take homology with closed
support. E.g. for a single
$n$-dimensional cell $A\approx\R^n$ we have $\chi(A)=(-1)^n$.
\begin{defn}[Viro, \cite{Viro-chi}]
We define
$$\int_{\rp^2}fd\chi=\sum\limits_{j\in f(\rp^2)}j \cdot
\chi(f^{-1}(j)) = \sum\limits_{\sigma \text{ cell}} (-1)^{\dim(\sigma)} f(\sigma).$$
\end{defn}
The latter sum is taken over the open cells $\sigma$ such
that $f|_\sigma$ is constant.

\begin{prop}\label{f0}
For any smoothable locally constant function $f:\rp^2\setminus K \to \R$
we have
$$\int\limits_{\rp^2}fd\chi=\int\limits_{\rp^2}f_\circ d\chi,$$
where the integrals are considered for the extensions of $f$ and $f_\circ$ to $\rp^2$.
\end{prop}
\begin{proof}
We may choose the cellular decompositions of $\rp^2$ for $f$ and $f_\circ$
so that one decomposition can be obtained from the other by replacing each node
$p$ of $K$ with the vanishing cycle $I_p$ subdivided into three cells:
the two endpoints and the relative interior.
The contribution of $I_p$ to $\int\limits_{\rp^2}f_\circ d\chi$
coincides with the contribution of $p$ to $\int\limits_{\rp^2}f d\chi$.
\end{proof}

The invariants $J_\pm$ and $\St$ from \cite{A}
also have
counterparts for immersions to $\rp^2$.

\begin{defn}[Viro, \cite{Viro-immersions}] \label{Jminus}
Let $K \subset \rp^2$ be a generic immersion of an oriented circle.
Put
$$J_-(K)=1-\int\limits_{\rp^2}\ind^2_Kd\chi.$$
\end{defn}

Furthermore, in \cite{Viro-immersions} it was shown that
the number $J_-(K)$ defined in this way
does not change under the direct
double tangency perestroika or under
the triple point perestroika. It decreases
by 2 under the inverse self-tangency perestroika
when the number of nodes is increased.

As in the case of immersions in $\R^2$ we can use
the equality $J_+(K) - J_-(K)=n(K)$ as the definition of $J_+(K)$.
In our context it is more convenient to use the
integral $\int\limits_{\rp^2}\ind^2_Kd\chi$ itself
as the invariant of an immersion instead of either $J_+(K)$ or $J_-(K)$.
We set up the following definition accordingly.

\begin{defn} \label{Orinv}
The {\em complex orientation invariant} of $K$
is the number
$$\Or(K)=\int\limits_{\rp^2}\ind^2_Kd\chi.$$
\end{defn}

\begin{figure}[!ht]
\centerline{
\input{files/pics/exampleK2.TpX}
}
\caption{An example for $\Or(K) = 2 - \frac12 - \frac52 + 1 + 4 = 4$ and $\St(K) = -1 + \frac52 - \frac12 = 1$.}
\label{exampleK2}
\end{figure}

An example of this invariant is given in Figure \ref{exampleK2}.
Note that this invariant makes sense not only for
generic immersions of single circle, but also for generic immersions
of disjoint collections of circles.
In particular, it makes sense for $K_\circ$
(an embedded collection of circles obtained from $K$).

\begin{corollary}
We have $\Or(K_\circ)=\Or(K)$.
\end{corollary}
\begin{proof}
The corollary follows from Proposition \ref{f0} since
$\ind^2_{K_\circ}=(\ind^2_K)_\circ$.
%
\end{proof}

For the rest of this subsection we assume that
$K\subset\rp^2$ is an immersion of a single circle.
Following \cite{Viro-immersions} we introduce
the strangeness invariant for generic
immersions in $\rp^2$ with the help
of the Shumakovich' formula \cite{Shu}.
Let us choose a base point $q\in K$ that
is not a node of $K$.
Each node $p\in K$ now gets a canonical local orientation
since we have the order on
the oriented branches of $K$
as they can be traced from $q$ following
the orientation of $K$. This
local orientation can be used to fix the sign
of $\ind_{K,q}(p)$. Namely, we set
$+1$ to represent the class of a small oriented circle around $p$
in $H_1(\rp^2\setminus\{p\};\Z)$.


\begin{defn}[Shumakovich-Viro, \cite{Shu}]\label{St-defn}
We set
$$\St(K)=-\sum\limits_p \ind_{K,q}(p)+
\ind^2_K(q)-\frac12,$$
where the sum is taken over all nodes $p\in K$.
\end{defn}

The number $\St(K)$ does not depend on the choice of the
base point $q\in K$ and stays invariant under both direct and
inverse self-tangency perestroikas. It changes by $\pm 1$
under the triple point perestroika and thus provides the
counterpart of Arnold's strangeness for generic immersions
of a circle into the projective plane (see Figure \ref{exampleK2}).

\begin{figure}[!ht]
\centerline{
\input{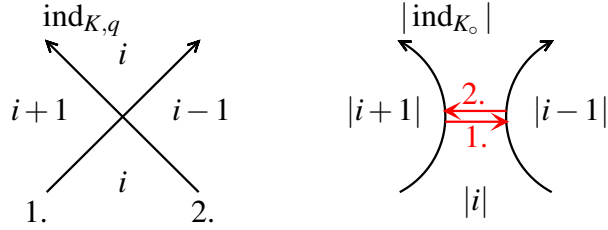}
}
\caption{The consistency of the two sign rules, namely $\text{sgn}(i) = - \sigma_q(I_p)$}
\label{signsagree}
\end{figure}

The relative interior of a $K_\circ$-membrane $I_p$ is contained in a single component
of $\rp^2\setminus K_\circ$. Thus $|\ind_{K_\circ}(I_p)|$ is well-defined.
Definition \ref{St-defn} may be rewritten as follows
in terms of the smoothing diagram $\Delta_K$:
\begin{equation}\label{St-Delta}
\St(K)=\St(\Delta_K)=\sum\limits_p \sigma_q(I_p)|\ind_{K_\circ}(I_p)|
+\ind^2_K(q)-\frac12.
\end{equation}
Here the base point $q$ can be placed anywhere on $K_\circ$
outside of the endpoints of the vanishing cycles.
To define the signs $\sigma_q(I_p)=\pm 1$ we trace
the diagram $\Delta_K$ starting from $q$ according to the orientation
of $K_\circ$ and jump to the other branch of $K_0$ at every
vanishing cycle $I_p$.
Note that whenever $|\ind_{K_\circ}(I_p)|\neq 0$ there is a
canonical orientation of $I_p$ from the arc of $K_\circ$
with a smaller $|\ind_{K_\circ}|$ to the arc with the
larger  $|\ind_{K_\circ}|$.
If the first jump at $I_p$ is in the direction of this
orientation, then we set $\sigma_q(I_p)=1$,
otherwise $\sigma_q(I_p)=-1$.
Figure \ref{signsagree} shows why this sign choice agrees
with the previous one. On the left hand side, we depicted
the local behaviour of $\ind_{K,q}$ around a node $p$.
Here, the choice of first and second branch induces the local orientation to
be counterclockwise, and this implies that $\ind_{K,q}$ decreases by 1 when
we cross a branch from left to right.
On the right hand side we depicted the two jumps between the branches of the local smoothing,
together with the values of $|\ind_{K_\circ}|$.
We want to compare the signs of $i = \ind_{K,q}(p)$ and $\sigma_q(I_p)$.
Assume $i > 0$, then $|i+1| > |i-1|$, and the first jump is from greater to lower value of $|\ind_{K_\circ}|$,
hence $\sigma_q(I_p)=-1$.
It follows $\text{sgn}(i) = - \sigma_q(I_p)$, which shows that
definition \ref{St-defn} and formula \eqref{St-Delta} are consistent.

\subsection{The immersion graph
\texorpdfstring{$\Gamma(K)$}{Gamma(K)}}
\label{immersiongraph}

To classify generic immersions
of a disjoint union of circles
with a given
number of nodes we can list their smoothing diagrams
(see Proposition \ref{prop-DeltaK}).
In turn, to exhaust such diagrams
it is useful to extract a certain graph
from the smoothing diagram and study
the properties of this graph.

Let $K \subset \rp^2$ be a generic immersion of a disjoint union of $k$ oriented circles.
We form a graph $\Gamma(K)$ whose vertices are the connected
components of $K_\circ$ and edges are the vanishing cycles of $\Delta_K$.
\begin{lemma}
If $l(K)=n(K)+k$, then the graph $\Gamma(K)$ is
a disjoint union of $k$  trees.
If $k=1$, then
$\Gamma(K)$ is connected
and $b_1(\Gamma(K))=n(K)+1-l(K) = 0 \mod 2$.
\end{lemma}
As usual, we denote with $b_1(\Gamma)$ the first Betti number
of a graph $\Gamma$, {\it i.e.}, the number of independent cycles of $\Gamma$.

\begin{proof}
If $l(K) = n(K) + k$, then
each vanishing cycle decreases the number
of connected components of the
smoothing diagram
exactly by one.
This implies
the first
statement
of the lemma.
If $K$ is connected, the graph $\Gamma(K)$ is also connected.
The modulo 2 congruence is
provided by
Lemma \ref{parity-Kcirc}.
\end{proof}

\begin{defn}\label{enhanced-Gamma}
The {\em immersion graph} of a generic immersion $K\subset\rp^2$
is the graph $\Gamma(K)$ enhanced with the following extra
information {\rm (}so that the smoothing diagram $\Delta_K$
can be uniquely recovered from the enhanced graph $\Gamma(K)${\rm )}.
The enhancement of $\Gamma(K)$
consists of
\begin{itemize}
\item {\em (1.\ oriented edges.)}
the orientation for {\em some}
of the edges of $\Gamma(K)$,
\item {\em (2.\ ribbon structure.)}
the ribbon structure
for $\Gamma(K)$,
\item {\em (3.\ two-coloring.)}
a two-coloring of its vertices and
\item {\em (4.\ projective enhancement.)}
a {\em projective enhancement} of $\Gamma(K)$
described below {\rm (}the treatment is different depending on $d(K)${\rm )}.
\end{itemize}
We describe these enhancements in details one by one.
\end{defn}

\noindent {\em 1. Oriented edges.}
Each edge $E$ of $\Gamma(K)$ connects two components of $K_\circ$.
If these components are ovals with disjoint interiors then
we do not orient $E$.
Also we do not orient $E$ if it is a loop-edge, {\it i.e.},
corresponds to a vanishing cycle connecting an oval to itself.

If the two components are nested ovals,
{\it i.e.}, the interior of one oval contains the other oval, then we
orient the edge in the direction from the interior oval to the
exterior one.

Recall that all connected components of $K_\circ$ except
possibly for a single one (which appears if $d(K)\neq 0$)
are ovals. If an edge $E\subset\Gamma(K)$
connects an oval to the pseudoline,
then we orient $E$ in the direction from the oval to the pseudoline.
Note that no vanishing cycle may connect a pseudoline in $\rp^2$
to itself by the orientation reasoning.

\noindent {\em 2. Ribbon structure.}
As there is some ambiguity in
what people call {\em ribbon graphs} we start
by quoting one of the most commonly accepted
definitions of an ribbon graph.
It is a finite graph $\Gamma$ with a choice of
cyclic order of adjacent half-edges for every vertex $v\in \Gamma$.

Recall that given such a structure on $\Gamma$
we may canonically reconstruct the {\em oriented} surface $S_\Gamma$ containing
the graph $\Gamma$ as its deformational retract as follows.
We take a copy of the closed oriented
disk $D_v$ for each vertex $v\in\Gamma$
of valence $\val(v)$.
Then we mark $\val(v)$
points at the boundary $\dd D_v$
so that the boundary orientation agrees with the cyclic
order from the ribbon structure.
Finally, we attach an orientation-preserving ribbon $R_E$ connecting
the disks $D_v$ and $D_{v'}$ at the corresponding
marked points for each edge $E$
connecting $v$ and $v'$.

\newcommand{\Cy}{\operatorname{Cy}}
\ignore{
We enhance $\Gamma(K)$ with a more general notion of
a ribbon structure (which may produce perhaps a non-orientable
surface in the way described above) according to the definition below.

\begin{defn}
A {\em ribbon structure} on a finite graph $\Gamma$ is a finite graph
enhanced with the following structure.
\begin{itemize}
\item
For each vertex $v$ we are given a
cyclic order of the edges adjacent to $v$
up to reversing this order.
We denote the set composed of these two
(opposite) cyclic orders with
$\Cy(v)$.
\item
Each edge $E$ connecting vertices
$v$ and $v'$ (including the case $v=v'$)
gives a preferred bijection
$$\rho_E:\Cy(v)\approx\Cy(v').$$
\end{itemize}
\end{defn}

As in the orientable case, any ribbon graph $\Gamma$
defines a surface $S_\Gamma$ containing
the graph $\Gamma$ as its deformational retract,
though this surface is not always orientable.

To construct $S_\Gamma$ we take a closed
disk $D_v$ for each vertex $v\in\Gamma$
of valence $\val(v)$.
Then we mark $\val(v)$
points at the boundary $\dd D_v$
in accordance with the orders from $\Cy(v)$
(two possible orientations of $D_v$ give two
cyclic orders from $\Cy(v)$).

Then we attach a ribbon $R_E$ connecting
$D_v$ and $D_{v'}$ at the corresponding
marked points for each edge $E$
connecting $v$ and $v'$.
Note that there are two distinct ways to attach
$R_E$ and each of these ways defines
a bijection between $\Cy(v)$ and $\Cy(v')$
by transporting a local orientation from
$D_v$ to $D_{v'}$ via $R_E$. We choose
the one specified by the ribbon structure on
$\Gamma$.
}

Recall that our immersion $K$ is oriented and so are all
components of $K_\circ$.
Because of this, the graph $\Gamma(K)$ admits a canonical ribbon structure
coming from the cyclic order of vanishing cycles
adjacent to a connected component $C_v\subset K_\circ$ associated to
a vertex $v\in\Gamma(K)$.

Nevertheless, the surface $S_\Gamma$ (that is conventionally associated to
a ribbon graph as we reviewed above) does not have a direct relation
to the topology of $K_\circ\subset\rp^2$.
Instead, our goal is to use additional structure that we define on $\Gamma(K)$
so that we may recover the pair $(\rp^2;\Delta(K))$.
We describe this procedure in the next subsection,
once we define the remaining part of the enhancement for $\Gamma(K)$.

\noindent {\em 3. Two-coloring.}

The vertices of $\Gamma$ are colored by two colors in the following way.
An oriented pseudoline $J\subset\rp^2$ defines
the orientation of $\rp^2\setminus J$.
An oriented
oval disjoint from $J$ is
{\it positive} (respectively, {\it negative})
if it defines in its interior the orientation opposite to (respectively, the orientation coinciding with)
the one given by the orientation of $\rp^2 \setminus J$.
If $d(K)$ is odd we choose $J=C_u$, the pseudoline
component of $K_\circ$.
If $d(K)$ is even, we
fix some pseudoline $J \subseteq \rp^2 \setminus K_\circ$.
We then set the vertices corresponding to all
positive ovals to be white vertices,
all other vertices (corresponding to negative ovals or the pseudoline)
are blue vertices.
If $d(K)$ is even, due to the choice of $J$ the colors are only well-defined up the following
flip. We say a vertex $v$ \emph{dominates} a vertex $v'$ if they can be connected via an oriented edge leading to $v$.
Let $v$ be maximal with respect to this partial order ({\it i.e.}, $v$ corresponds to an exterior oval).
Then the colors are well-defined up to flipping the colors of $v$ and all the vertices dominated by $v$ simultaneously.
  The colors are a convenient way to describe the following property of $\Gamma(K)$.
	For every cycle in $\Gamma(K)$, the number of non-oriented edges not contained in the root cluster
	is even.

Before we describe the projective enhancement,
let us collect some properties of the
first three enhancements satisfied by
immersion graphs $\Gamma=\Gamma(K)$.

\begin{defn}
A {\em cluster of vertices} is
a maximal connected subgraph
$\Gamma'\subset\Gamma$ such that all
its edges are unoriented
{\rm (}i.e., any other connected subgraph
 without oriented edges is either contained
 in $\Gamma'$ or disjoint from it{\rm )}.
\end{defn}

\begin{prop}[Oriented edges]
\label{or-edges}
If $\Gamma=\Gamma(K)$ is the immersion graph
of a generic immersion $K \subset \rp^2$,
then
all
oriented outgoing edges from a
cluster of vertices $A\subset\Gamma(K)$
lead to the same vertex $v\in\Gamma(K)$.

There is a unique cluster $U$ of vertices
of $\Gamma(K)$ such that there are no outgoing
edges from $U$.
In the case $d(K)\neq 0$ this cluster consists of the
vertex
associated to the pseudoline component of $K_\circ$.
\end{prop}

\begin{proof}
An oval $C_v\subset K_\circ$ might be contained inside several other ovals of $K_v$,
but all of them are nested. A vanishing cycle may connect $C_v$ only with the innermost oval
from this collection.
The second statement follows from the fact that $K$ and therefore $\Gamma$ are connected.
\end{proof}

If all
outgoing edges from a
cluster of vertices $A\subset\Gamma(K)$
lead to a vertex $v\in\Gamma(K)$,
we say that the cluster $A$ is {\it dominated by $v$}.
The unique cluster $U$ of $\Gamma(K)$ such that there are no outgoing edges from $U$
is called the {\it root cluster} of $\Gamma(K)$.

\begin{prop}
\label{coloring}
Let $E$ be an edge connecting
vertices $v, v'\in\Gamma$.
If $E$ is oriented, then $v$ and $v'$
are of the same color.
If $E$ is non-oriented and not contained in the root cluster,
then $v$ and $v'$ are of different color.
\end{prop}

\begin{proof}
%
Let $J$ denote the pseudoline used to define the colors.
In both cases of the statement, we may assume that the vanishing cycle $I$
corresponding to $E$ does not intersect $J$. The statement than follows
from the coherence of $I$ (see \ref{def-coherent-membrane}).
\end{proof}

\newcommand{\Ad}{\operatorname{Ad}}
For a vertex $v$, we denote with $\Ad(v)$
the set of all
edges adjacent to $v$.
The sets $\Ad^+(v)$  and $\Ad^-(v)$
are the subsets of $\Ad(v)$ formed by the edges
oriented away from $v$ and towards $v$, respectively.

Let $A$ be a cluster of vertices. Build a surface
$\Sigma_A$ as follows. Take
a disjoint union $\Delta$
of oriented discs $D_v$
over all $v\in A$.
For each vertex $v \in A$,
mark
points,
indexed by $\Ad(v) \setminus \Ad^-(v)$,
on $\dd D_v$ in the cyclic order
provide by
the ribbon structure of $A$. For each edge
connecting $v,v'\in A$ we
add a ribbon to $\Delta$ connecting small neighborhoods
of the corresponding marked points at $\dd\Delta$ so that
the ribbon {\em disrespects} the orientations
of $D_v$ and $D_{v'}$.
Denote the resulting surface with $\Sigma_A$.

\begin{prop}\label{planarity}
Suppose that $\Gamma=\Gamma(K)$ is an immersion graph.
Let $A$ be a cluster which is not the root cluster.
Then the surface $\Sigma_A$ is homeomorphic to a sphere with holes.
Moreover, all the remaining
marked points
{\rm (}indexed by $\bigcup_v \Ad^+(v)${\rm )}
lie on the same boundary component of $\Sigma_A$.
\end{prop}

The latter component is called the \emph{exterior} boundary component.

\begin{proof}
The surface $\Sigma_A$ is homeomorphic to
a regular neighborhood $N$
of the part of the smoothing diagram formed by the ovals
and vanishing cycles from $A$.
All remaining marks sit on the exterior boundary component of $N$
(the boundary of the non-oriented component in $\rp^2\setminus N$).
\end{proof}

\begin{remark} \label{combprop}
  The previous proposition can also be expressed more combinatorially.
	Let $A$ be a cluster of $\Gamma$ which is not the root cluster.
	A \emph{ribbon cycle} in $A$ is an oriented closed path $p$ in $A$ such
	at each white (resp.\ blue) vertex $v$ the outgoing edge is the
	successor (resp.\ predecessor) of the incoming edge
	(according to the cyclic order given by the ribbon structure).
	The oriented edges $e \in \Ad^+(v)$ lying between the incoming
	and outgoing edge are said to \emph{lie on $p$}. Let
	$V,E,B$ be the number of vertices of $A$, edges of $A$, ribbon cycles of $A$, respectively.
	Then the proposition can be reformulated as
	\[
	  V - E + B = 2
	\]
	and the condition that all edges in $\bigcup_{v\in A} \Ad^+(v)$ lie on the same ribbon cycle.
\end{remark}

Let us finally discuss the fourth enhancement.

\noindent{\em 4. Projective enhancement.}

If $d(K) = 0$, consider the surface $\Sigma_A$,
where $A$ is the root cluster. If
$\Sigma_A$ is oriented, the projective enhancement
is a choice of boundary component of $\Sigma_A$,
namely the exterior one (see the proof of
Proposition \ref{planarity}). If $\Sigma_A$
is non-oriented, no extra information is needed.

If $d(K)=1$, {\it i.e.}, $[K]\neq 0\in H_1(\rp^2)$,
there is a pseudoline component of $K_\circ$.
We treat the corresponding vertex $u\in\Gamma(K)$
as the {\em root vertex}.
In figures, we
draw the root vertex in a special way, as
a horizontal line.

As any other vertex of $\Gamma(K)$ the vertex $u$ comes with
a natural cyclic order on the adjacent edges.
This order is
the cyclic order
in which the vanishing cycles appear when we travel along
the pseudoline $C_u\subset K_\circ$.

There is however an extra data we can extract from the way
how the vanishing cycles are attached to $C_u$.
Even though the component $C_u\subset\rp^2$ is one-sided,
locally it has two sides and we can check whether the two
consecutive vanishing cycles come to $C_u$ from the same side,
or from opposite sides.
In the latter case we place a cross on the corresponding
arc of the circle corresponding to the vertex $u\in\Gamma(K)$.
Since $C_u$ is one-sided, the total number of crosses on
the circle corresponding to $u$ must be odd.
This information is the projective enhancement in the case
$d(K)=1$.

If $v\neq u$ is not a root vertex,
we define the a cyclic order on $\Ad^-(v)$
induced by the cyclic order on $\Ad(v)$.
For the root vertex $u$ we define the cyclic
order on $\Ad^-(u)$ by requiring that
the edge following $E\in\Ad^-(u)$ is the first edge
after $E$ (in the order defined by the orientation
of the pseudoline $C_u$) that
is separated from $E$
by an even number of crosses.

Let $A$ be a cluster, and let $\Ad^+(A) = \bigcup_v \Ad^+(v)$
denote the edges adjacent to a cluster $A$
and oriented in the direction outgoing from $A$.
(Note that $\Ad^+(U)=\emptyset$ for the root
cluster $U$.)
Due to Proposition \ref{planarity}, the orientation
on the exterior boundary component
of $\Sigma_A$ induces a
cyclic order on $\Ad^+(A)$.

\begin{prop}
\label{clusters}
Let $\Gamma=\Gamma(K)$ be an immersion graph.
For any cluster $A\subset\Gamma$ dominated by a vertex $w$
the cyclic order on $\Ad^+(A)\subset\Ad^-(w)$
agrees with that on $\Ad^-(w)$.
Furthermore, the sets $\Ad^+(A)$ for different clusters
dominated by $w$ are unlinked,
i.e., if $A$ and $B$ are two clusters dominated by $w$,
there is a segment in $\Ad^-(w)$ {\rm (}according to the cyclic order{\rm )}
that contains $\Ad^+(A)$ and is disjoint from $\Ad^+(B)$.
\end{prop}

\begin{proof}
The proposition follows since every
component of $\rp^2\setminus K_\circ$ adjacent to $w$ and $A$
is homeomorphic to an open disk with punctures.
\end{proof}

\begin{definition}
  \label{def:enhancedgraph}
  A graph $\Gamma$ enhanced with the a partial orientation
	of its edges, a ribbon structure, a two-coloring and
	a projective enhancement as defined in paragraphs
	1. -- 4. above and satisfying the properties
	of propositions \ref{or-edges}, \ref{coloring},
	\ref{planarity} and \ref{clusters} is called an
	\emph{enhanced graph}.
\end{definition}


We conclude these considerations
with the following elementary statement.

\begin{prop} \label{criteria-odddegree}
Let $\Gamma$ be an enhanced graph for odd $d$.
Then there exists a generic
immersion $K\subset\rp^2$ of a collection of disjoint
oriented circles with $\Gamma=\Gamma(K)$.
The isotopy type of $K$ is determined by $\Gamma$.
\end{prop}
\begin{proof}
We construct the smoothing diagram
in $\rp^2$ by representing the root vertex $u$
with a pseudoline and the other vertices with ovals.
Propositions \ref{or-edges}, \ref{clusters}
and \ref{planarity} ensure
existence of
such diagram while
Proposition \ref{coloring} ensures that the orientations
are compatible.
Then we ``unsmooth" $K_\circ$, {\it i.e.}, replace each
vanishing cycle with a node as in Figure \ref{or-smoothing}.
Topological uniqueness is inductive.
\end{proof}

\begin{remark}
  \label{rem:deven}
	In the case $d=0$, only few changes are necessary.
	Proposition \ref{criteria-odddegree} still holds if we add
	the requirement that the surface $\Sigma_A$ for the root cluster $A$
	is homeomorphic to either a sphere with holes or $\rp^2$ with holes.
	(In terms of \ref{combprop}, we want $V - E + B = 1$ or $2$.)
	Then this surface can be embedded in $\rp^2$ by gluing discs to all
	boundary components, except possibly to the boundary component
	specified by the projective enhancement, to which we glue a M\"obius strip
	instead. The construction then continues as in the odd case.
\end{remark}

\ignore{
Indeed, all oriented edges incoming to $v$
correspond to one of these sides.
We call these edges {\em inner}.
All other edges adjacent to $v$ (which may be outgoing or non-oriented)
correspond to the other side.
We call these edges {\em outer}.
Note that according
to our convention for partial orientation of $\Gamma(K)$ all
the edges adjacent to $u$ are oriented and incoming (with respect to $u$).

\subsection{Properties of \texorpdfstring{$\Gamma(K)$}{Gamma(K)}}
Here we collect formal properties of the enhanced graph $\Gamma(K)$
so that it can come from a smoothing diagram of immersed collection of circles
(which we assume throughout the subsection).
\begin{defn}
A {\em cluster of vertices} is
a maximal connected subgraph
$\Gamma'\subset\Gamma$ such that all
its edges are unoriented
(i.e. any other connected subgraph
 without oriented edges is either contained
 in $\Gamma'$ or disjoint from it).
\end{defn}
Note that if a cluster of vertices contains a root vertex then
it must consist solely of the root vertex (as all edges adjacent to
the root vertex are oriented).
Otherwise, the cluster consists of the ovals.
Since no pair of these ovals forms a nest,
all the ovals from the cluster are adjacent
from the outside
to the same connected component of $\rp^2\setminus K_\circ$.
We obtain the following proposition.

\begin{prop}[Oriented edges]
\label{or-edges} 
All oriented and outgoing edges from a
cluster of vertices $A\subset\Gamma(K)$
lead to the same vertex $v\in\Gamma(K)$.
We call such a cluster $A$ {\em dominated by $v$}.

There is a unique cluster $U$ of vertices
of $\Gamma(K)$ such that there are no outgoing vertices from $U$.
In the case $d(K)\neq 0$ this cluster consists of the root vertex of $\Gamma(K)$.
We call this cluster {\em the root cluster}.
\end{prop}

To describe other properties of $\Gamma(K)$ it is convenient to
color its vertices into two colors.
If $d(K)\neq 0$ then the smoothing $K_\circ$ has a pseudoline component $C_u$
(corresponding to the root vertex $u$). Note that the orientation of $C_u$
(induced by the orientation of $K$) induces the orientation of
$\rp^2\setminus C_u$. Likewise, the orientation of each oval $C_v$
(corresponding to a non-root vertex $v$) induces the orientation of
the interior of $C_v$.

\begin{defn}[Coloring]\label{def-coloring}
We say that an oval $C_v$ is {\em negative}
(i.e. $C_u$ and $C_v$ form a negative pair)
and color the vertex $v$ {\em black}
if the orientations of the interior of $C_v$ induced by the orientations
of $C_u$ and $C_v$ coincide.
Otherwise, we say that $C_v$ is {\em positive}
(i.e. $C_u$ and $C_v$ form a negative pair)
and color the vertex $v$ {\em white}.
\end{defn}

If $d(K)=0$ then all connected components of $K_\circ$
are ovals. Because of non-orientability of $\rp^2$ we need
to make some auxiliary choice to define colors on the vertices
of $\Gamma(K)$. Let $C_u$ be a choice of an oriented pseudoline
disjoint from $K_\circ$ (even though the graph $\Gamma(K)$
does not have any root vertex $u$ in this case).
Once we have made this choice, any oval of $K_\circ$
acquires a sign and any vertex of $\Gamma(K)$ acquires a color
according to Definition \ref{def-coloring}.
\begin{prop}
If $d(K)=0$ then any cluster of vertices except for the root cluster
is bipartite according to the coloring defined by Definition \ref{def-coloring}.
An edge of the root cluster connects two vertices of the same color
if and only if the corresponding vanishing cycle intersects $C_u$.
\end{prop}
\begin{proof}
The proposition follows from the fact that the pseudoline $C_u$
represents the first Stiefel-Whitney of $\rp^2$.
\end{proof}

Clearly, in the case $d(K)=0$ we have ambiguity in coloring that comes
from the choice of $C_u$.
Namely, we can change the color of any vertex $v$ of the root cluster
by isotoping the reference pseudoline $C_u$ around the oval $C_v$
(topologically we may think of the oval $C_v$ to go across the pseudoline $C_u$).
This results in changing the color of $v$ as well as
the color of all vertices in the clusters dominated by $v$.

Let us define the (simplicial) 1-cochain $W\in C^1(\Gamma(K);\Z_2)$
as follows.
If an edge $E\subset\Gamma(K)$ connects vertices of different colors
we set $W(E)=0$. Otherwise we set $W(E)=1\in\Z_2$.

\begin{prop}[Coloring in the case $d(K)=0$]
The 1-cochain $W$ represents the 1-cohomology class
$w_1(\Gamma(K))\in H^1(\Gamma(K);\Z_2)$ from the enhancement
of Definition \ref{enhanced-Gamma}.

In particular, all clusters of $\Gamma(K)$ except for the root cluster
form bipartite graphs.
\end{prop}
\begin{proof}
We have the cochain $W$ defined by $C_u$ while $C_u$ represents
the first Stiefel-Whitney class of $\rp^2$ as it is a pseudoline.
\end{proof}
Similarly, we have the following proposition
for the odd degree case (when the root is a single vertex $u\in\Gamma(K)$).
\begin{prop}[Coloring in the case $d(K)\neq 0$]
Any cluster of $\Gamma(K)$ in the case $d(K)\neq 0$ is bipartite.
\end{prop}

\newcommand{\col}{\operatorname{color}}
We can use the coloring of $\Gamma(K)$ together with
its ribbon structure 
to recover the surface $S_{\Gamma(K)}^{\col}$ that can be obtained
from the surface $S_{\Gamma(K)}$ by twisting those ribbons
that correspond to the edges connecting vertices of opposite color.
Namely, we adopt the following definition.
\begin{defn}[Surface $S^{\col}_{\Gamma}$]
\label{def-coloribbon}
Given a colored ribbon graph $\Gamma$ we associate
to it the surface $S^{\col}_{\Gamma}$ in the following way.
For each vertex $v\in\Gamma$
of valence $\val(v)$ we take
an oriented disk $D_v$ with
$\val(v)$ marked points on the boundary
$\dd D_v$ so that each marked point corresponds
to an edge adjacent to $v$ and so that
the cyclic order induced by the orientation
of $\dd D_v$ agrees with the ribbon structure.

Consider the disjoint union of such disks
over all vertices of $\Gamma$.
Now for each edge $E\subset\Gamma$ we attach
a ribbon $R_E$ at the marked point corresponding
to its edges.
If the endpoints of $E$ have different
colors then we attach the ribbon $R_E$
so that
it respects the orientation of the disks
it joins (i.e. so that the union of these
        disks and $R_E$ admits the orientation
        agreeing with the orientation of the disks).
If the endpoints of $R_E$ are of the same color
then we attach the ribbon $R_E$ otherwise
(i.e. so that an orientation of $R_E$
     agrees with the orientation of the disks
     at one of the endpoints, but not the other).
\end{defn}


We use Definition \ref{def-coloribbon} to identify
the properties of the enhanced graph $\Gamma(K)$
that ensures that $\Delta(K)$ can be embedded to $\rp^2$.
To do this we associate a ribbon graph $A_{\Gamma'}$
to every cluster $\Gamma'\subset\Gamma(K)$.
If $\Gamma'$ is a root cluster we set $A_{\Gamma'}=\Gamma'$.
Otherwise, $\Gamma'$ is dominated by a vertex $v_{\Gamma'}$
by Proposition \ref{or-edges} and we set $A_{\Gamma'}=\Gamma'\cup\{v_{\Gamma'}\}$.
The coloring of $A_{\Gamma'}$ is induced by the coloring
of $\Gamma(K)$.
The ribbon structure of $A_{\Gamma'}$ is the choice
of cyclic order of the adjacent edges at its vertices.
For all vertices of $A_{\Gamma'}$ except for the
root vertex $u$ we induce this cyclic order from that of
the ribbon structure of $\Gamma(K)$.
For the root vertex $u$ (which exists if $d(K)\neq 0$)
we define the cyclic order of the adjacent edges according
to the order they appear at the boundary of the regular
neighborhood of the pseudoline corresponding to $u$
(recall that such neighborhood is diffeomorphic to a M\"obius band).

\begin{prop}[Planarity of extended clusters]
If $\Gamma'\subset\Gamma$ is a non-root cluster then
the surface $S^{\col}_{A_{\Gamma'}}$ can be uniquely embedded
to the sphere $S^2$.
If $\Gamma'\subset\Gamma$ is a root cluster then
the surface $S^{\col}_{A_{\Gamma'}}$ can be uniquely embedded
to the projective plane $\rp^2$.
\end{prop}
\begin{proof}
If $\Gamma'$ is not the root cluster then it is composed of an oval $O_{\Gamma'}$
(corresponding to the dominant vertex of $\Gamma'$) and all ovals $O$
that lie directly inside $O_{\Gamma'}$, i.e. such that $O$ and $O_{\Gamma'}$
are adjacent to the same connected component of $\rp^2\setminus K_\circ$.
Topologically this component is a sphere $S^2$ with punctures. The surface
$S^{\col}_{\Gamma'}$ is embedded to this $S^2$ (after we compactify it
by adding back the punctures). 

If $\Gamma'$ is the root cluster and $d(K)\neq 0$ then it consists of a single
point. If $d(K)=0$ then $\Gamma'$ is composed of ovals that are adjacent to
the unique connected component of $\rp^2\setminus K_\circ$.
Topologically this component is $\rp^2$ with some punctures.
We get an embedding of $S^{\col}_{\Gamma'}$ to $\rp^2$ in a similar way as before.

Note that these embeddings are unique up to to a homeomorphism of the target
($S^2$ or $\rp^2$) since we can recover the target by gluing disks to
all connected components of $S^{\col}_{\Gamma'}$.
\end{proof}

...
}
\ignore{
\begin{defn}
Two pairs of edges $\{e_A,e'_A\}$ and $\{e_B,e'_B\}$
all adjacent to the same non-root vertex $v\in\Gamma(K)$ are
called {\em linked} if they appear as $e_A,e_B,e'_A,e'_B$
or $e_A,e'_B,e'_A,e_B$
according to $\Cy(v)$ (note that for this property
it does not matter which of the two opposite cyclic order at $v$
to consider). Otherwise they are called {\em unlinked}.
\end{defn}
\begin{prop}[Planarity inside all ovals]
\label{v-planarity}
If $A$ and $B$ are two different cluster of vertices dominated
by the same non-root vertex $v$ and $\{e_A,e'_A\}$ (resp. $\{e_B,e'_B\}$)
is a pair of edges adjacent to the cluster $A$ (resp. B)
and the vertex $v$ then  $\{e_A,e'_A\}$ and $\{e_B,e'_B\}$
are unlinked.
\end{prop}
\begin{proof}
Linking of the corresponding pairs of edges implies
non-zero intersection of the parts of the vanishing cycle
diagrams corresponding to the clusters $A$ and $B$.
\end{proof}

If $d(K)=1$ then there exists the root vertex $u\in\Gamma(K)$.
There is a modification of Proposition \ref{v-planarity} for
the clusters of vertices dominated by $u$.
\begin{defn}
Two pairs of edges $\{e_A,e'_A\}$ and $\{e_B,e'_B\}$
all adjacent to the same root vertex $u\in\Gamma(K)$ are
called {\em linked} if the following two conditions hold
simultaneously.
\begin{itemize}
\item
The edges appear as $e_A,e_B,e'_A,e'_B$
or $e_A,e'_B,e'_A,e_B$
according to $\Cy(u)$.
\item
Let us consider the four segments of $u$ given by
the the four edges $e_A,e_B,e'_A,e'_B$.
Then exactly one of these segments contains
an even number of crosses. \todo{}
\end{itemize}
Otherwise the two pairs of edges $\{e_A,e'_A\}$ and $\{e_B,e'_B\}$
are called {\em unlinked}.
\end{defn}
\begin{prop}[Planarity at the root vertex]
\label{U-planarity}
Suppose that $d(K)\neq 0$. If $A$ and $B$
are two different cluster of vertices dominated
by the root vertex $u$ and $\{e_A,e'_A\}$ (resp. $\{e_B,e'_B\}$)
is a pair of edges adjacent to the cluster $A$ (resp. B)
and the vertex $u$ then  $\{e_A,e'_A\}$ and $\{e_B,e'_B\}$
are unlinked.
\end{prop}
\begin{proof}
Linking of the corresponding pairs of edges implies
non-zero intersection of the parts of the vanishing cycle
diagrams corresponding to the clusters $A$ and $B$.
\end{proof}

\begin{prop}
\label{noloopd1}
If $d(K)\neq 0$ then $\Gamma(K)$ contains no loop edges.
All of its unoriented edges connect vertices of different color
while all of its oriented vertices connect vertices of the same color.
\todo{colors?}
\end{prop}
\begin{proof}
If vanishing cycle connects an oval $K_\circ$ to itself
and is disjoint from the pseudoline then
the union of the oval and the cycle can be isotoped to
the orientable affine part $\R^2$ of $\rp^2$.
We get a contradiction with the orientation of $K_\circ$
at the endpoints of a vanishing cycle.

Similarly the pseudoline $J\subset K_\circ$ cannot be connected
to itself by a vanishing cycle as $\rp^2\setminus J\approx \R^2$.
Orientability of $\rp^2\setminus J$ also implies the condition
on colors of the endpoints of an edge according to its orientation.
\end{proof}

If $d(K)=0$ we may consider a certain planarity property
of the unique cluster $U$ of nest-depth 0.
This property involved the characteristic class
$w_1\in H^1(\Gamma(K);\Z_2)$.
Let us consider the continuous retraction
\begin{equation}
\label{U-retraction}
\rho_U:\Gamma(K)\to U
\end{equation}
\todo{} defined by sending every vertex $v$ with $\nu\delta(v) \ > 0$
to the unique vertex in the root cluster that dominates $v$.
Proposition \ref{or-edges}
ensures that this map is well-defined.
\begin{prop}
The characteristic class $w_1\in H^1(\Gamma(K);\Z_2)$
is induced by its restriction to $U$ through $\rho_U$,
i.e.
$$w_1=\rho_U^*(w_1|_U).$$
\end{prop}
\begin{proof}
The proposition follows from the observation that the interior
of any oval is contractible and thus the restriction of $w_1$ to it
is trivial.
\end{proof}

Furthermore, we may express the characteristic class $w_1$ through a
collection of edges (which is however not uniquely defined).
Let $W$ be a collection of open edges of $U\subset\Gamma(K)$.
We identify $W$ with the union of open edges in this collection
which can be seen as an open subset $W\subset U$.
Then there is the unique element $\alpha_W\in H^1(U,U\setminus W;\Z_2)$
such that $\alpha_W(E)=1\pmod2$ for every edge $E$ from the union $W$.
\begin{defn}
\label{W-collection}
A collection of open edges
is called {\em characteristic}
if $w_1(\Gamma_K)|_U$ is the image of $\alpha_W$ under the homomorphism
$$H^1(U,U\setminus W;\Z_2)\to H^1(U;\Z_2)$$
induced by the inclusion $U\subset(U,U\setminus W)$.
\end{defn}
The following proposition follows directly from Definition
\ref{W-collection}.
\begin{prop}
A characteristic collection exists.
Any characteristic collection $W$ may be used to recover
$w_1\in H^1(\Gamma(K);\Z_2)$.
Namely, the value of $w_1$ on any cycle $Z\subset\Gamma(K)$
takes a non-zero value if and only if $Z$ contains odd number of edges
from $W$.
\end{prop}
\begin{prop}\label{prop-bipartite}
For any characteristic collection $W$
the graph $\Gamma(K)\setminus W$ is bipartite. \todo{oriented edges?}
In the same time the edges of $W$
connect the vertices of the same color
(for any bipartite coloring of
$\Gamma(K)\setminus W$),
and have nest-depth zero.
\end{prop}
\begin{proof}
We can specify a local orientation of $\rp^2$ near the union
of the ovals of $K_\circ$ with the vanishing cycles not
corresponding to $W$. Then the ovals of $K_\circ$ split
into positively (counterclockwise) and negatively (clockwise)
oriented. If a vanishing cycle $I$ connects two ovals
from the same cluster then the two ovals
must have the opposite orientations.
Otherwise, one of the ovals sits inside another
and they have the same orientation.

We color the vertices corresponding to positively oriented ovals
with even $\nu\delta$ as well as those
corresponding to negatively oriented ovals
with odd $\nu\delta$ in one color.
We color all the other vertices in the other color.
\end{proof}

\begin{prop}
\label{samevertexloop}
Any loop-edge (an edge connecting a vertex to itself) of $\Gamma(K)$
is contained in any characteristic collection $W$.
Furthermore, all loop-edges of $\Gamma(K)$ must be attached
to the same vertex.
\end{prop}
\begin{proof}
A loop-edge connects an oval with itself,
thus it has to reverse orientation.
Suppose we have two loop-edges on distinct vertices. They are connected
in the complement of a characteristic collection $W$.
Thus the corresponding vanishing cycles must intersect.
\end{proof}

\todo{$U\setminus W$ can be disconnected}
The previous lemma may be generalized as follows.
Let $S(U\setminus W)$ be the compact surface with boundary
obtained as a small regular neighborhood of $U\setminus W$
in $\rp^2$.
Note that $U\setminus W$ is connected and so is $S(U\setminus W)$.
\todo{}

The following condition may be considered as the counterpart
of the planarity condition at $U$ in the case when $d(K)=0$.
\begin{prop}
\label{U-linking}
All vanishing cycles in a characteristic collection $W$
connect points from the same connected component
$C^W_\infty\subset\dd S(U\setminus W)$.
The endpoints of any two vanishing cycles from $W$
are linked in $C^W_\infty\approx S^1$.
This means that an arc connecting
two endpoints of one vanishing cycle in $A$
contains a single endpoint of the other vanishing cycle.
\end{prop}
\begin{proof}
There is a single nonorientable component
of $\rp^2\setminus S(U\setminus W)$.
This component must be a M\"obius band as $U$ is connected.
All vanishing cycles from $W$ must be contained in this component.
Since they are disjoint in this M\"obius band
and reverse local orientation in $S(U\setminus W)$,
their endpoints are linked.
\end{proof}
}

\ignore{
We conclude this treatment by the following straightforward theorem.
Recall that a {\em semiorientation} is an orientation treated up
to its reversal. This notion makes sense for disconnected orientable
manifolds (such as disjoint unions of several circles).
Namely, for any two local orientations at different
points of such a manifold a semiorientation tells if they are coherent.
Clearly, a connected orientable manifold (such as a circle $S^1$)
has a unique semiorientation.

\begin{thm}\label{thm-graph}
Any pair consisting of a connected finite
graph $\Gamma$ with $n$ edges and $m$ vertices
and a degree $d\in\Z_2 = H_1(\rp^2;\Z_2)$ which is
enhanced as in Definition \ref{enhanced-Gamma} and satisfying to all
properties of Propositions \ref{or-edges}--\ref{U-linking}
encodes a generic immersion $K\subset\rp^2$
(which we identify with its image in $\rp^2$)
of a semioriented
disjoint union of $k$ circles. Here $k\le n-m+2$ and $k\equiv n+m\pmod2$.

In particular, if $\Gamma$ is a tree (i.e. $n-m+2=1$) then
it always corresponds to a generic immersion of a circle.
\end{thm}
\begin{proof}
We start by recovering the smoothing diagram $\Delta(K)$
near the cluster $U$ of nest-depth zero.

If $d=1$ then we associate to the root vertex $u\in\Gamma$
a non-trivially embedded circle $C_u\subset\rp^2$ (so that
$C_u$ is isotopic to $\rp^1\subset\rp^2$).
We represent the edges of $\Gamma$ adjacent to $u$
(which are all oriented towards $u$ according to Proposition \ref{or-edges})
by disjoint small intervals orthogonal to $C_u$ at one of their endpoints.
We choose the order of their appearance at $C_u$ to
agree with $\Cy(u)$ and choose (locally) the side of the
adjacency of the consecutive intervals to $C_u$ to be the same
or opposite according to the projective enhancement of $\Gamma$.
Proposition \ref{U-planarity} allows to represent the vertices
of nest-depth one by ovals adjacent to the other side of the intervals
as well as to insert the vanishing cycles corresponding to
non-oriented edges from the clusters they form.
Inductively we proceed to the ovals of higher nest-depth.

If $d=0$ then we choose a characteristic collection $W$.
Then we represent the subgraph $U\setminus W$ by a collection
of non-nested ovals connected with vanishing cycles in the
affine part $\R^2\subset\rp^2$.
Proposition \ref{U-linking} allows us to extend this representation
to $W$ in agreement with (local) orientation of the ovals.
The rest of the reconstruction is inductive
by the nest-depth of the vertices as in the case of
$d=1$.

In both cases we get the smoothing diagram $\Delta$ consisting
of oriented circles as well as vanishing cycles which are
in agreement with their orientation so that if we replace
each vanishing cycle with two (transversally intersecting)
oriented intervals we get a generic immersion of an oriented
1-manifold $K$ once we choose an arbitrary orientation at a
component corresponding to any vertex of $\Gamma$.

The normalization $\tilde K$ of the immersed 1-manifold $K$
is diffeomorphic to a disjoint union of $k$ circles.
The manifold $\tilde K$ is obtained from the disjoint union
of $n$ circles by (oriented) connected sum corresponding
to the $m$ vanishing cycles. Performing each connected sum
increases or decreases the number of components by 1.
Thus $k\equiv n+m$. Connected sum performed along
a maximal tree in $\Gamma$ must be connected. Thus
$k\le n-m+2$.
\end{proof}

\begin{corollary}
\label{coro-graph}
Any generic immersion of a circle into $\rp^2$ is encoded
by a unique finite connected graph enhanced
as in Definition \ref{enhanced-Gamma} with $m\equiv n+1\pmod 2$
(though some of such enhanced graphs with $m<n+1$ do not correspond
to immersions of a single circle).
\end{corollary}

\subsection{Drawing of the graph \texorpdfstring{$\Gamma(K)$}{Gamma(K)}}
\label{draw-graph}
Corollary \ref{coro-graph} allows us to use graphs $\Gamma(K)$
to encode generic immersions of a circle into $\rp^2$.
Recall that the result $K_\circ$ of smoothing of $K$
(according to one of its orientation) is a smooth multicomponent
projective curve that is the disjoint union of ovals if $d=0$
or the disjoint union of some ovals and the unique non-oval curve
(called {\em pseudoline}) if $d\neq 0$.
The pseudoline $J\subset K_0$ is isotopic to a projective line.

We choose the affine part $\R^2=\rp^2\setminus\rp^1$
so that it contains all ovals of $K_\circ$ and so that
$J\setminus\R^2$ consists of a single point (in the case $d\neq 0$).
We choose an orientation on $\R^2$ as well as an orientation
on $K_\circ$ compatible with an orientation on $K$.

We denote a vertex corresponding to an oval of $K_\circ$
with a small circle which can be either black or white,
depending whether the orientation on it agrees or disagrees
with that on $\R^2$. The pseudoline is denoted with a horizontal
line above all ovals oriented from left to right.

The edges of $\Gamma(K)$ are intervals (vanishing cycles)
from the diagram $\Delta(K)$.
In our drawing they are edges connecting circles to circles
or circles to the upper horizontal line.
Note that by Proposition \ref{noloopd1}
there cannot be edges connecting the horizontal line to itself.

We draw edges on the plane so that the cyclic order of their appearance
at each circle (or the horizontal line) agrees with the ribbon structure of $\Gamma(K)$.
Namely, the appearance of vanishing cycles adjacent to a component of $K_\circ$
in the direction of a chosen orientation of $K$ has to be left to right
for the horizontal line, clockwise for black circles and counterclockwise
for white circles. We color all circles adjacent to the horizontal line black. \todo{}

There are crosses on the horizontal lines indicating change of the side of adjacency
of the vanishing cycles. If two vanishing cycles $I_1$ and $I_2$ are adjacent to $J$
locally at different sides of $J$ (resp. the same side)
when we follow on $J$ from $I_1$ to $I_2$
according to the chosen orientation of $K$ and the edge representing $I_1$
is to the left from the one representing $I_2$ then they have to be separated
with odd number of crosses (resp. even number of crosses).
As $J$ is one-sided the total number of crosses is odd.

Consider ambiguities in this representation. If there is a horizontal line
(i.e. if $d(K)\neq 0$) then there is a choice of a point where we cut
the pseudoline $J$ in order to obtain this line
as well as in placing crosses on that horizontal line.
Edges adjacent to the horizontal line divide it into intervals.
Only the parity of the number of crosses on each interval matters
for reconstructing $K$, so we may assume that there is one cross or none
on each such interval.
Also we adopt the convention that there is
no cross at the rightmost interval
(as they are glued into the same interval on $J$).
Graphs $\Gamma(K)$ are equivalent if they differ by cyclic shifts of edges adjacent to the horizontal line,
cf. Figure \ref{cyclic-shift}. Note that as a result of such shift we may get two crosses on the same
interval. According to our remark above we may erase any even number of crosses on the same interval. \todo{}
\begin{figure}[ht]
\centerline{
\includegraphics[height=22mm]{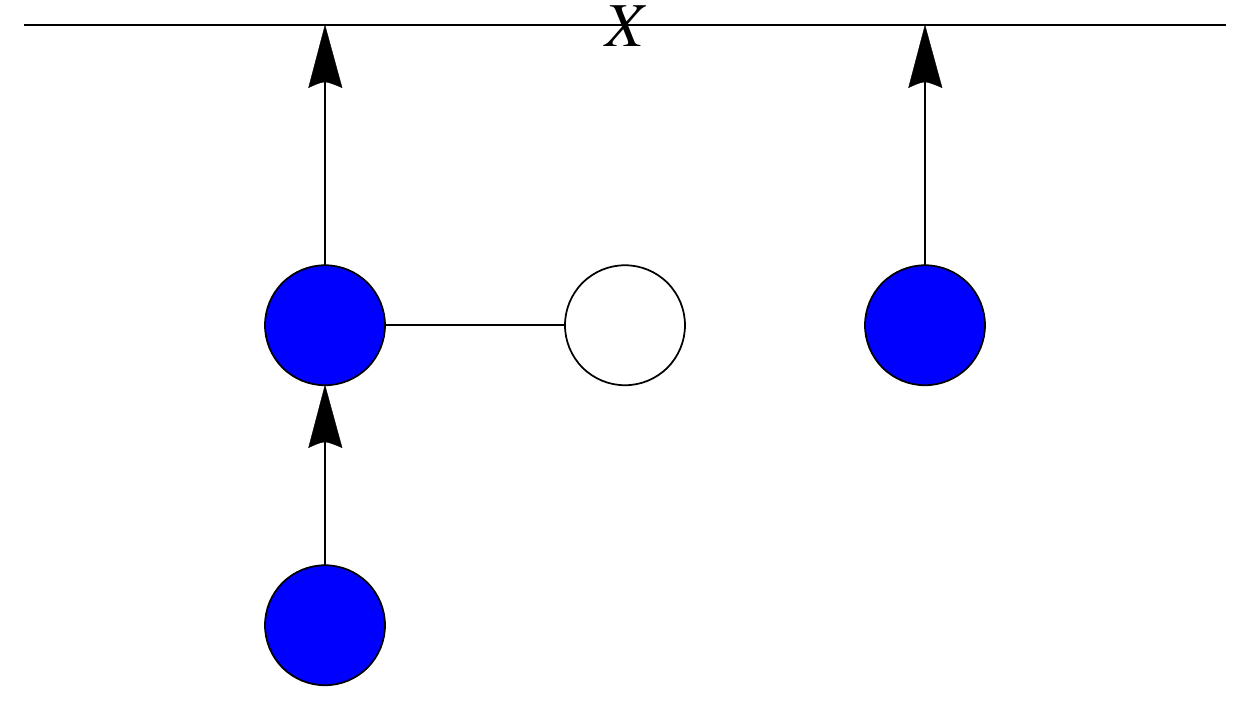}\hspace{28mm}
\includegraphics[height=22mm]{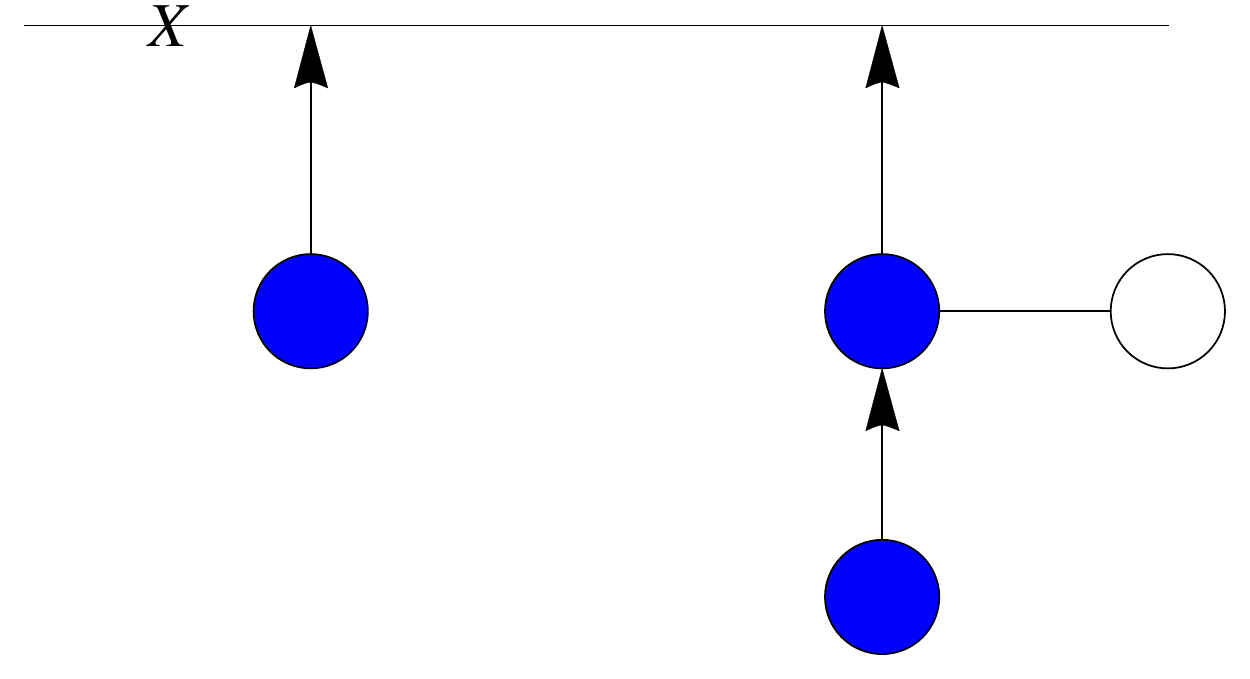}
}
\caption{Cyclic shift along the pseudoline $J$}
\label{cyclic-shift}
\end{figure}

If $d(K)=0$ there is no horizontal line, but there is a characteristic collection
of unoriented edges that connect vertices of the same color
(in our drawings we mark such edges with the "$\infty$" sign).
The ambiguity in the choice of a characteristic collection
comes from coloring of the vertices in the root cluster.
If we change the color of one of these vertices then
any unoriented edge connecting it to a different vertex will
become characteristic if it was not characteristic and will cease
to be characteristic if it was already characteristic, cf Figure \ref{char-collection}.
\begin{figure}[ht]
\centerline{
\includegraphics[height=25mm]{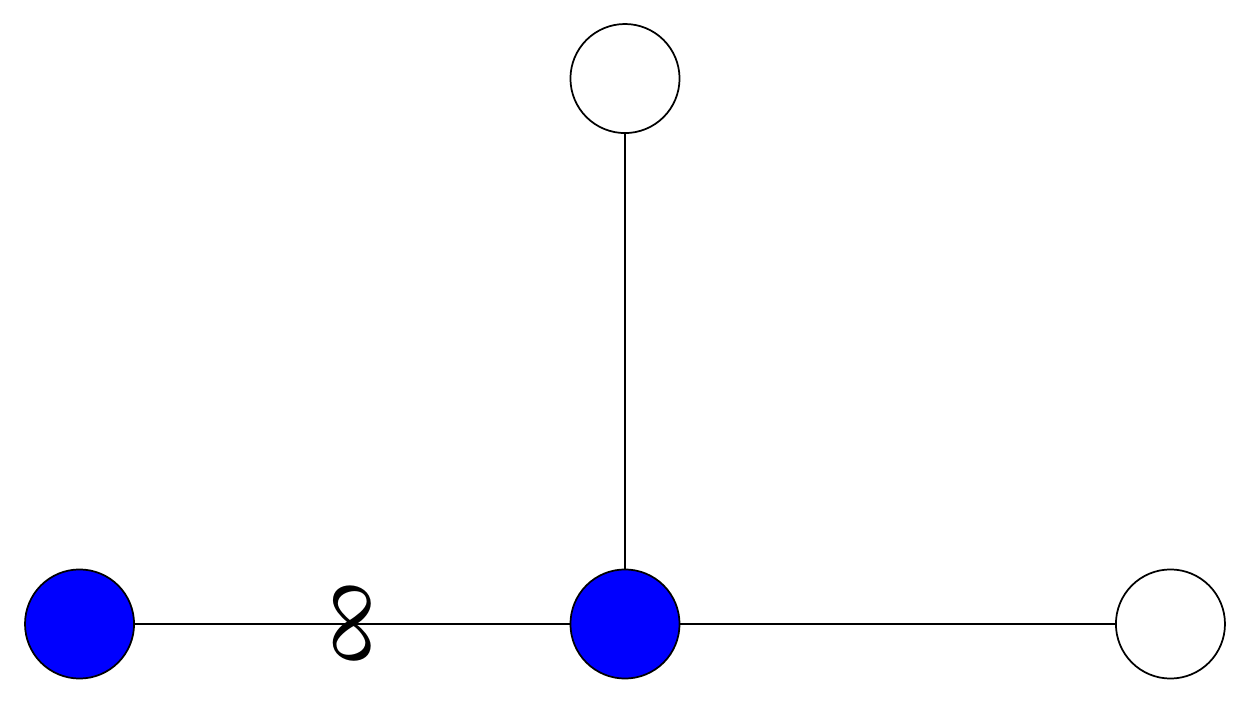}\hspace{28mm}
\includegraphics[height=25mm]{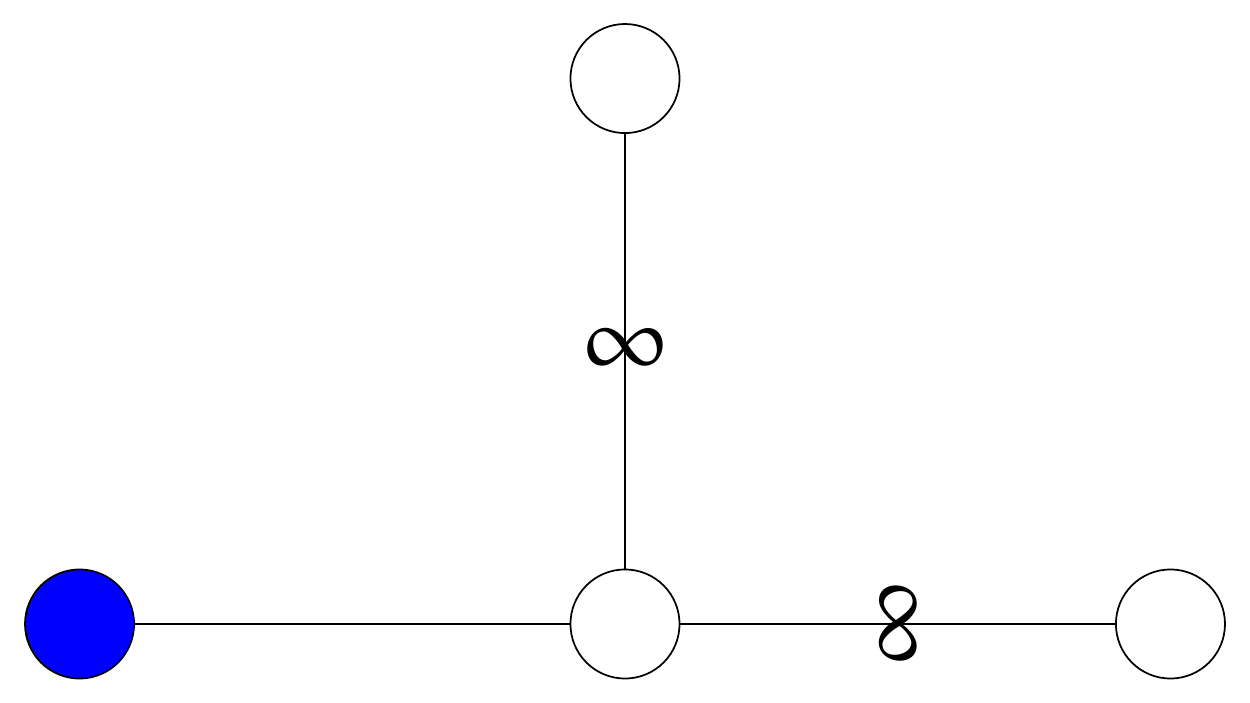}
}
\caption{Change of a characteristic collection}
\label{char-collection}
\end{figure}

We summarize these considerations in the following statement.
\begin{prop}
The diagram representing
a generic immersion of an oriented circle
(or a collection of oriented circle)
is unique up to the ambiguities described above.
\end{prop}

\subsection{Reidemeister moves in terms of \texorpdfstring{$\Gamma(K)$}{Gamma(K)}}
Let us now describe the two Reidemeister moves from Figures \ref{j+},
\ref{j-} and \ref{reidemeister3} in terms of this enhanced graph.

\noindent {\em The second Reidemeister move.} Recall that we have
distinguished between direct and inverse self-tangency perestroikas.
In terms of the graph $\Gamma(K)$ there are two options for each of these two
cases: when the two branches belong to distinct or to the same connected component
of $K_\circ$.

Figure \ref{R2-1} shows the corresponding moves for $\Gamma(K)$.
First to rows correspond the case when the two branches belong to distinct ovals.
Note that the two ovals may be of the same depth (i.e. dominated by the same
vertex of $\Gamma(K)$) or one oval may dominate the other.
This determines whether the two new edges are oriented or not (and respectively
whether the two vertices are of the same color).
Note that if $d=0$ and both ovals are from the root cluster we may also have
appearance of unoriented edges marked by $\infty$ between ovals of the same color.
These situations are depicted in the first three rows of Figure \ref{R2-1}.

The fourth row of Figure \ref{R2-1} represents
the case when one of the branch is an oval
while the other is a pseudoline.
Finally, the last row corresponds to the case when both branches belong to the same
oval. In this case $d(K)=0$ and the oval has to be from the root cluster.
In our figures the letters ($a, b, c$)
stand not necessarily for one, but for a collection of adjacent edges
connecting the depicted piece to the rest of the graph.
}

\begin{figure}[ht]
\begin{minipage}{0.45\textwidth}
\centerline{
\includegraphics[width=22mm]{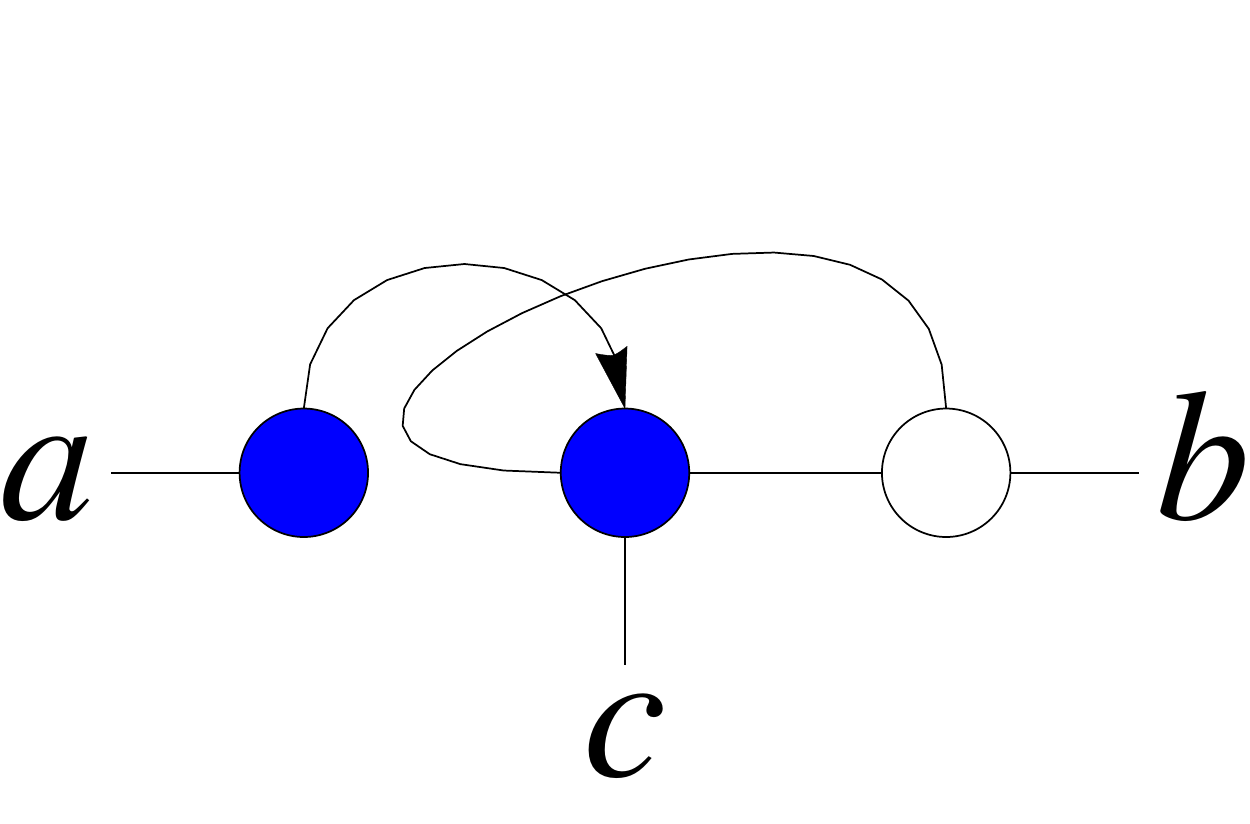}\hspace{10mm}
\includegraphics[width=22mm]{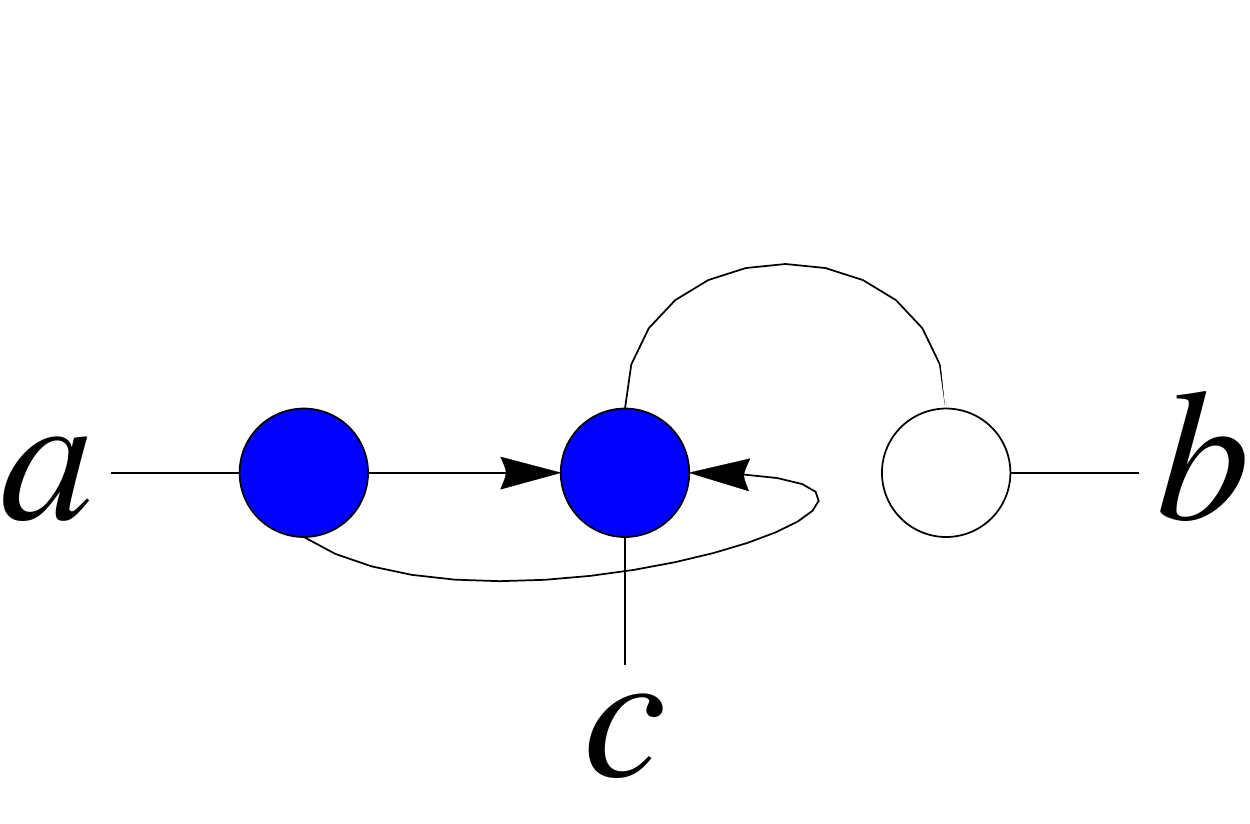}
}
\vspace{0mm}
\centerline{
\includegraphics[width=22mm]{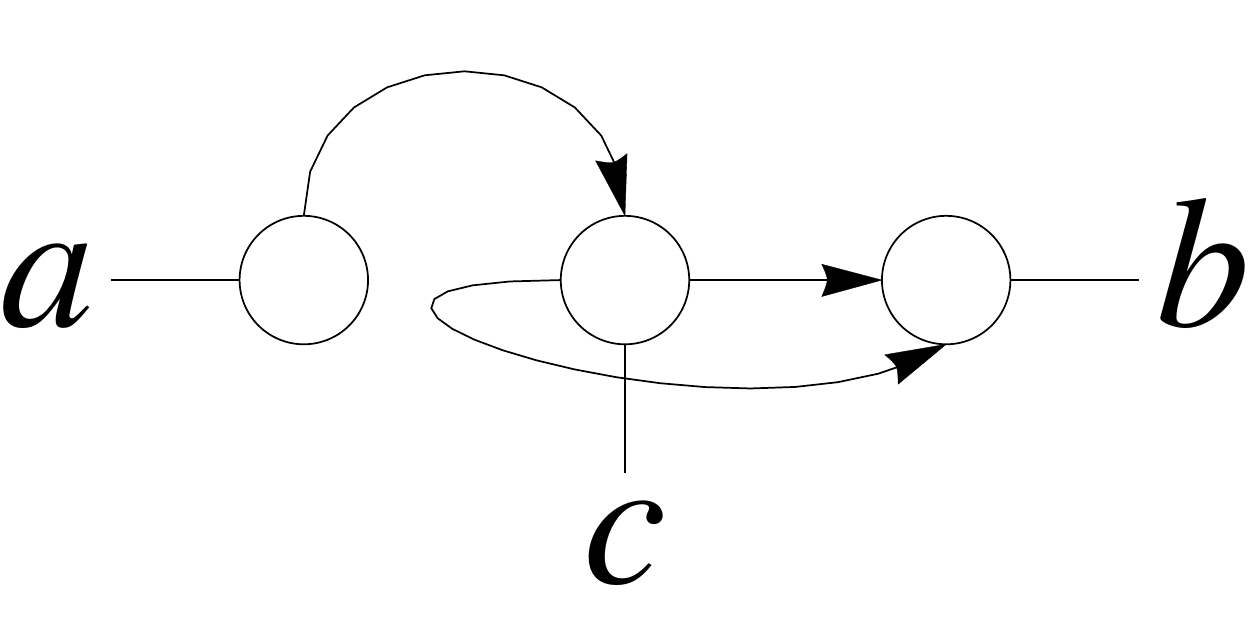}\hspace{10mm}
\includegraphics[width=22mm]{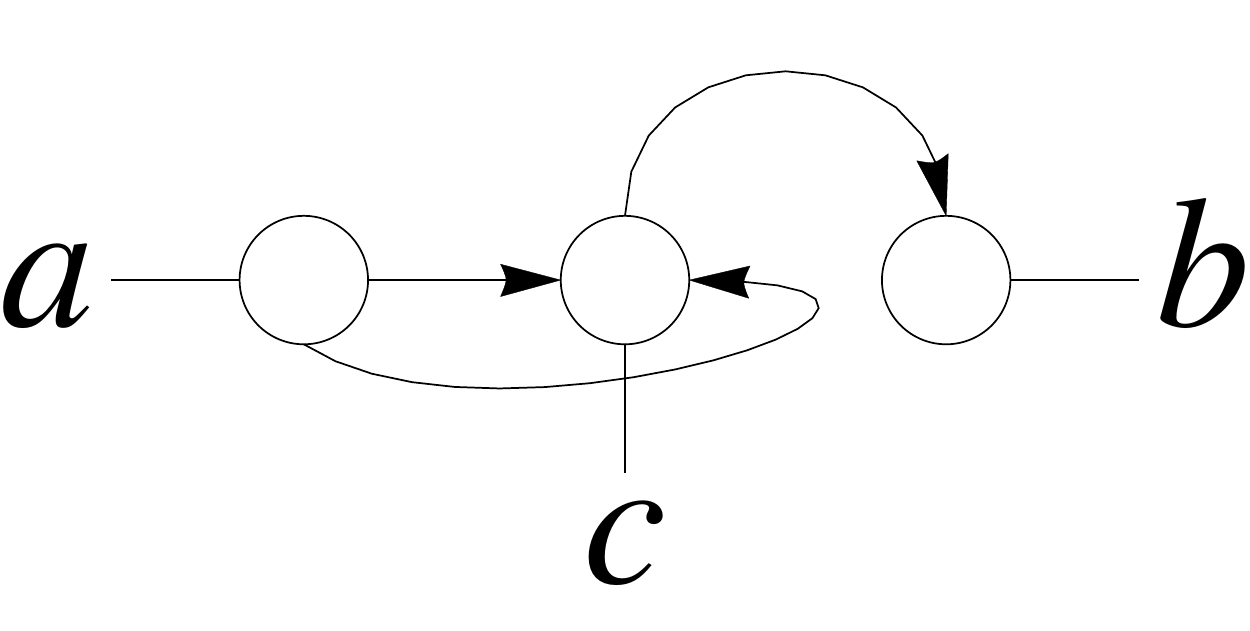}
}
\vspace{5mm}
\centerline{
\includegraphics[height=12mm]{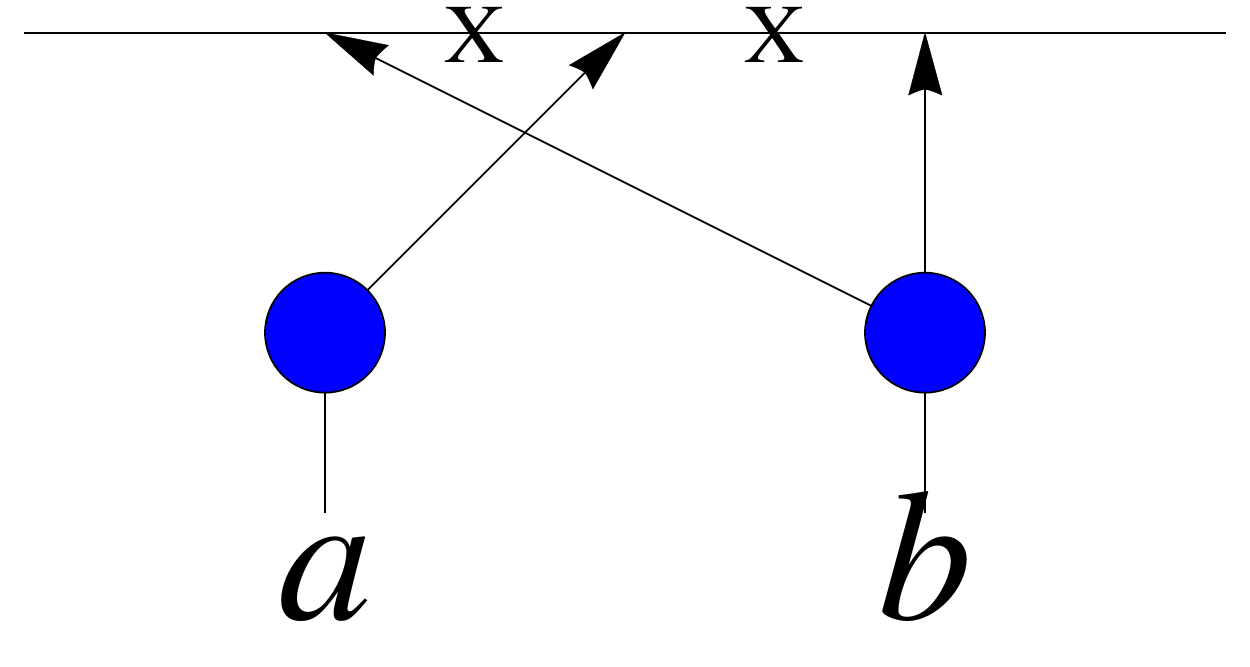}\hspace{10mm}
\includegraphics[height=12mm]{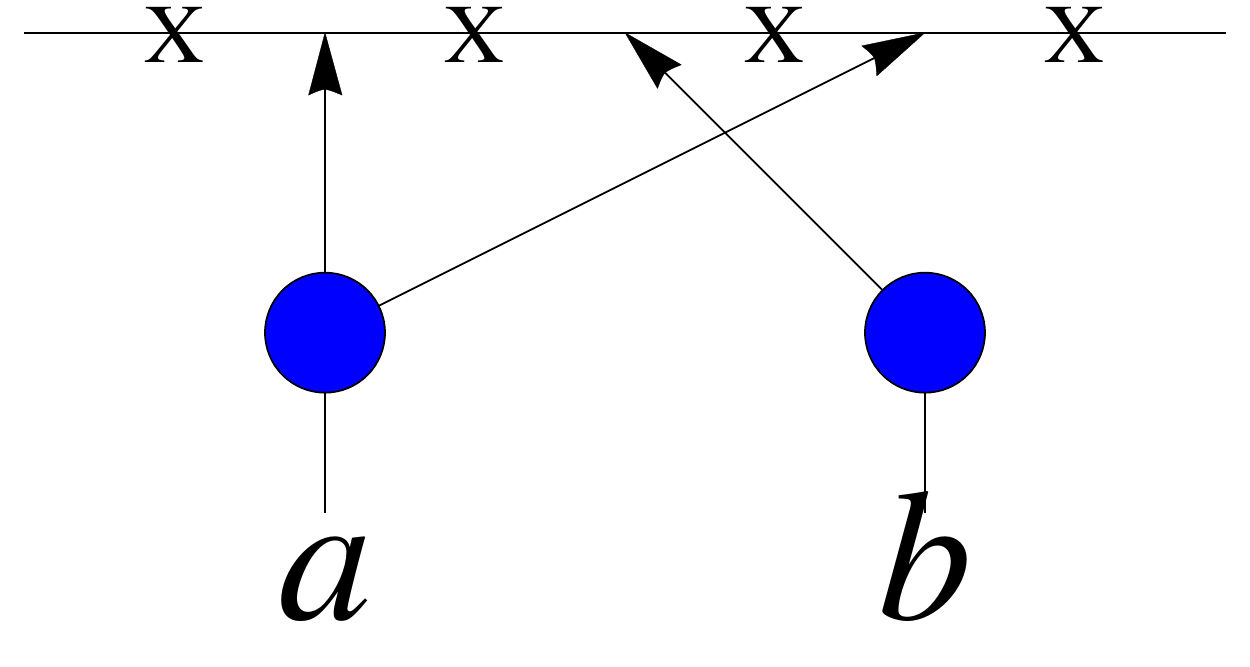}
}
\vspace{5mm}
\centerline{
\includegraphics[height=12mm]{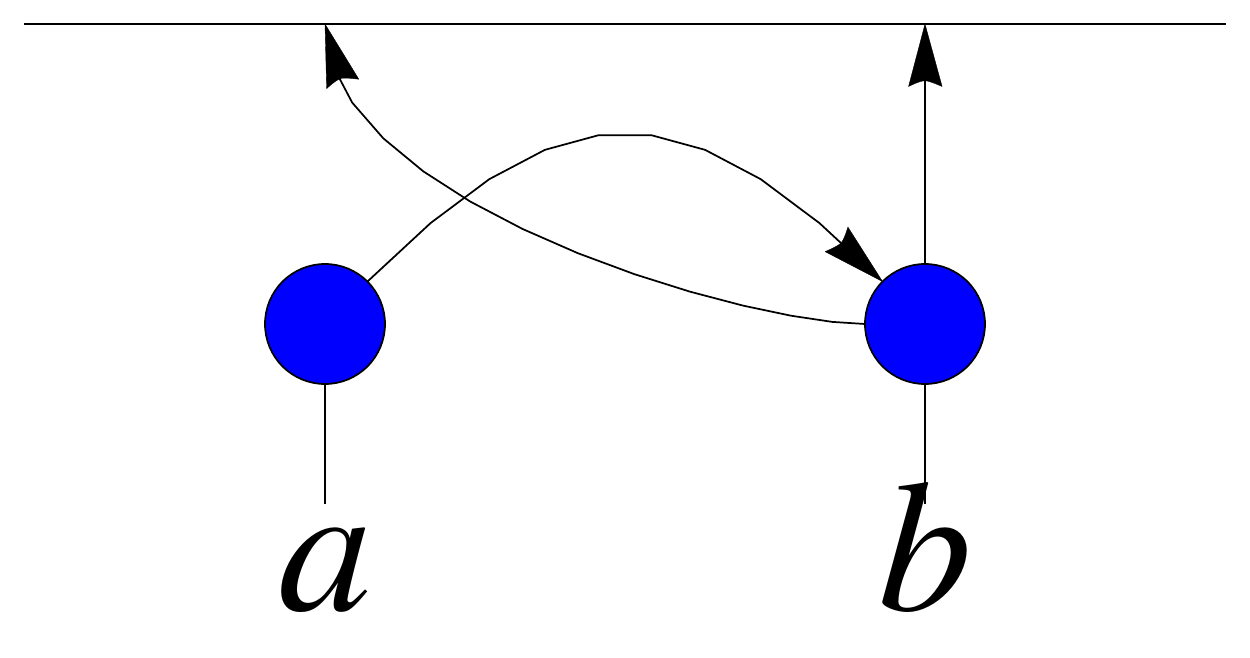}\hspace{10mm}
\includegraphics[height=12mm]{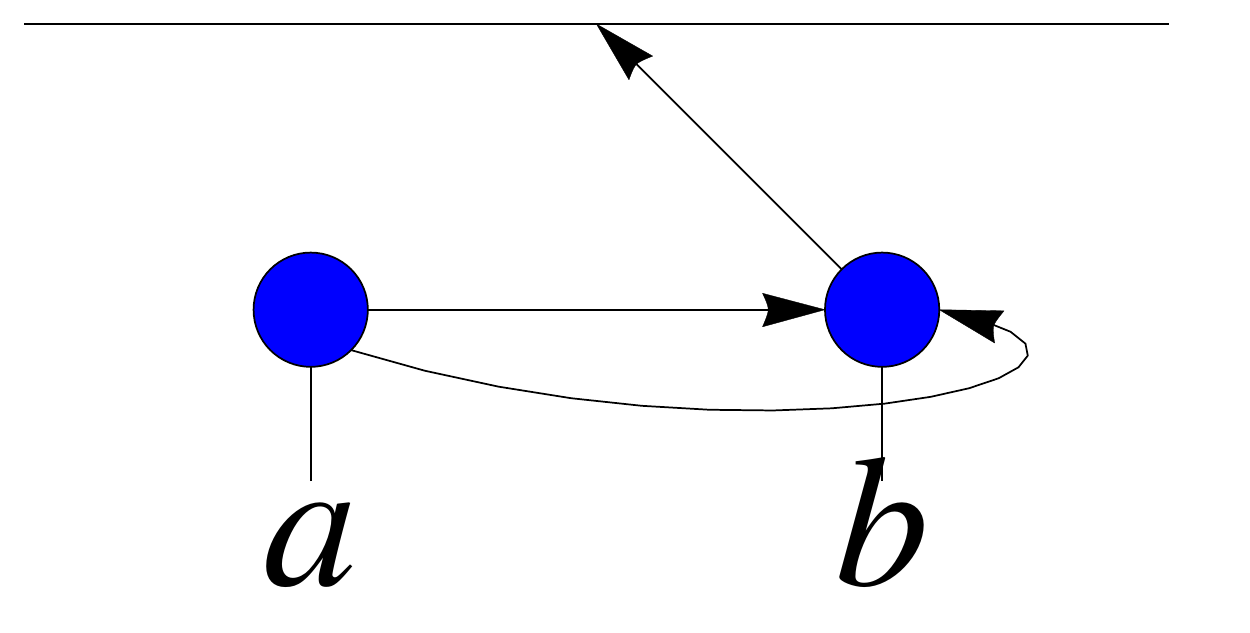}
}
\centerline{\vspace{2.5mm}
\includegraphics[height=19mm]{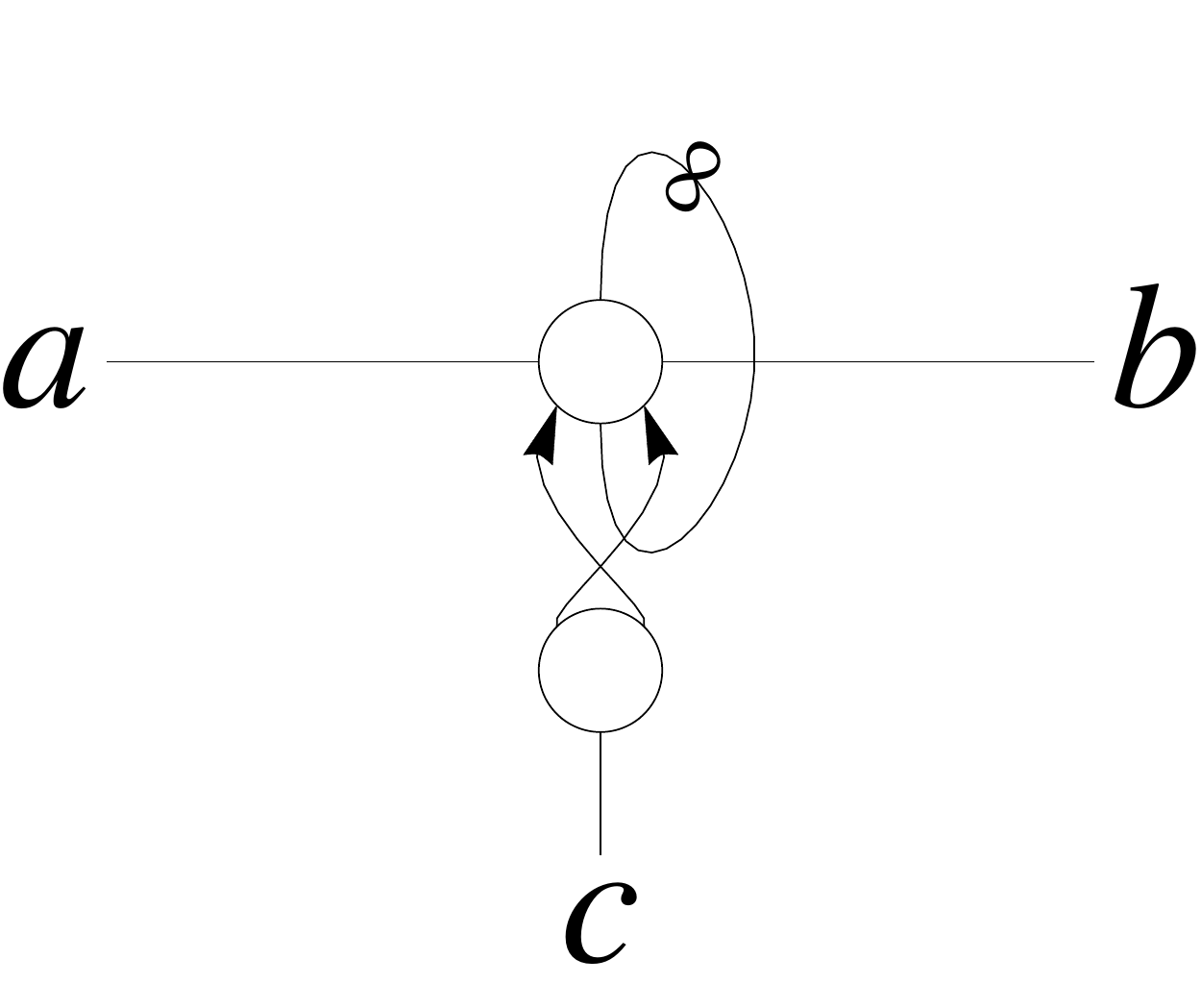}\hspace{10mm}
\includegraphics[height=19mm]{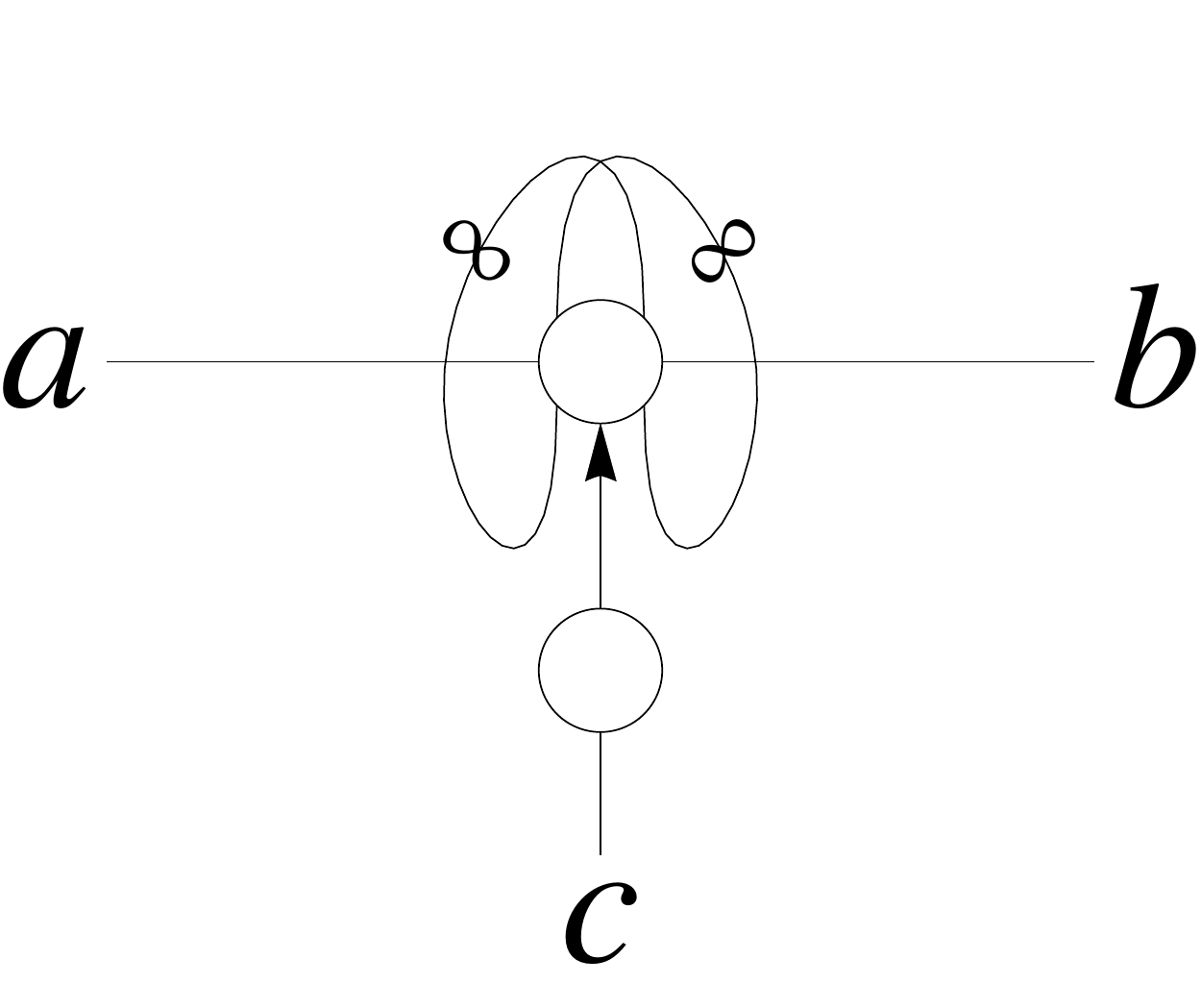}
}
\vspace{2.5mm}
\centerline{
\includegraphics[height=10mm]{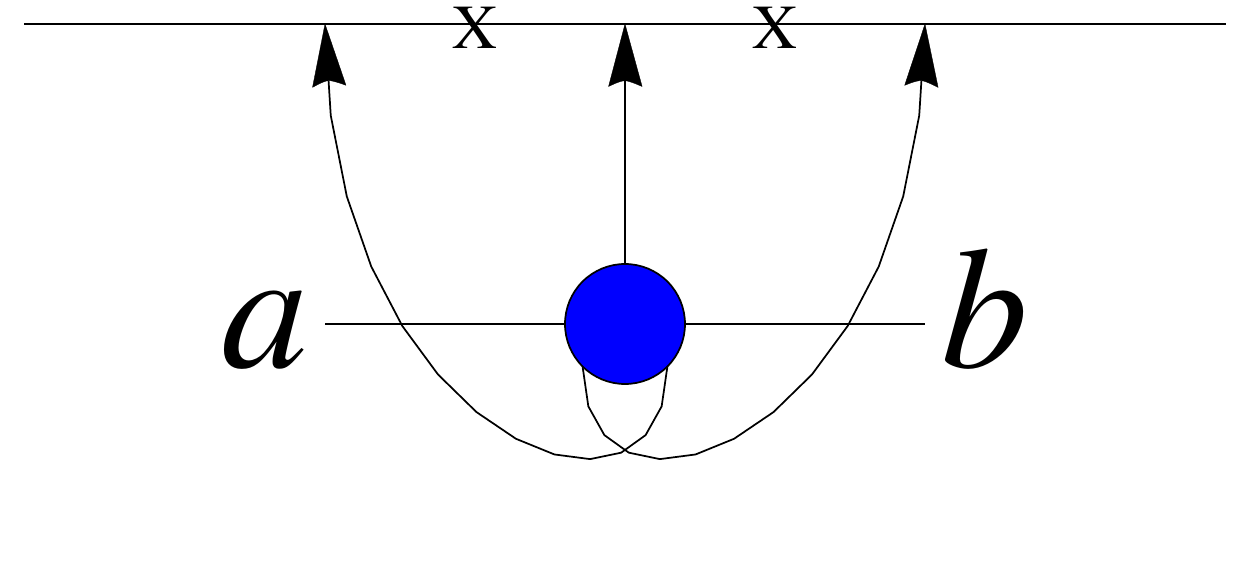}\hspace{10mm}
\includegraphics[height=10mm]{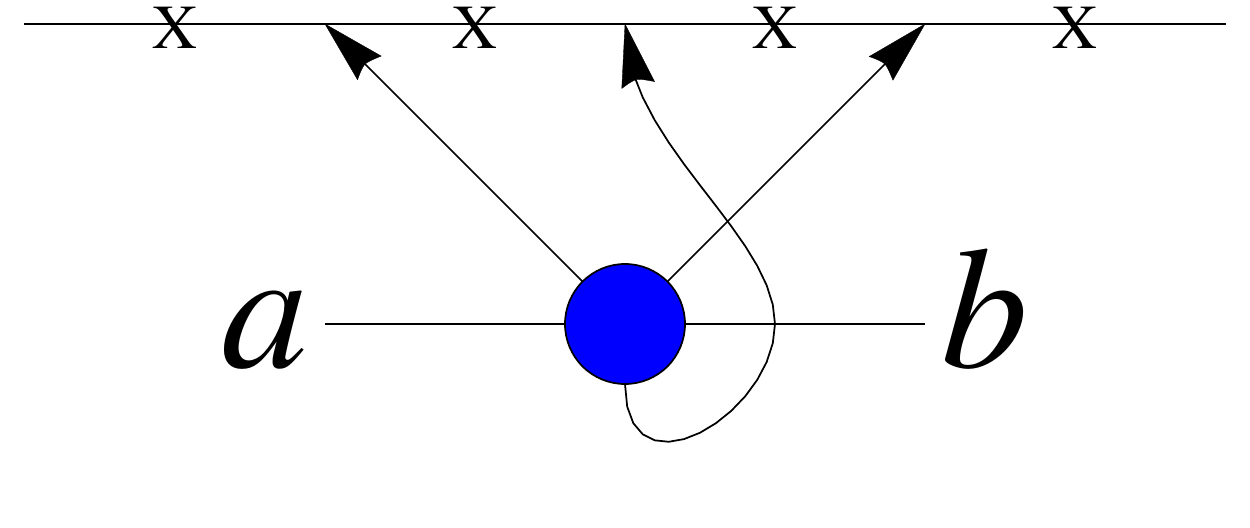}
}
\caption{Weak triple point perestroikas $\Gamma(K)$
\label{R3-1}}
\end{minipage}
\hspace{10pt}
\begin{minipage}{0.45\textwidth}
\centerline{
\includegraphics[width=20mm]{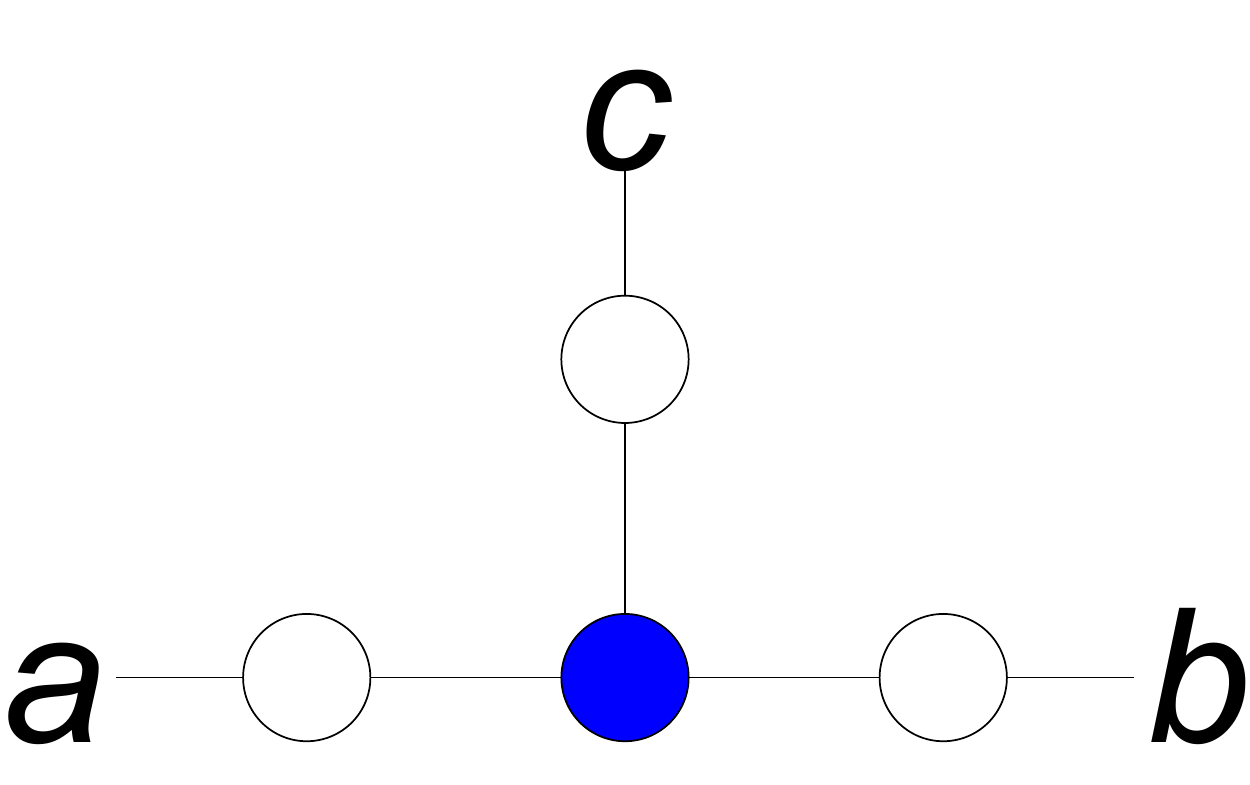}\hspace{5mm}
\includegraphics[width=20mm]{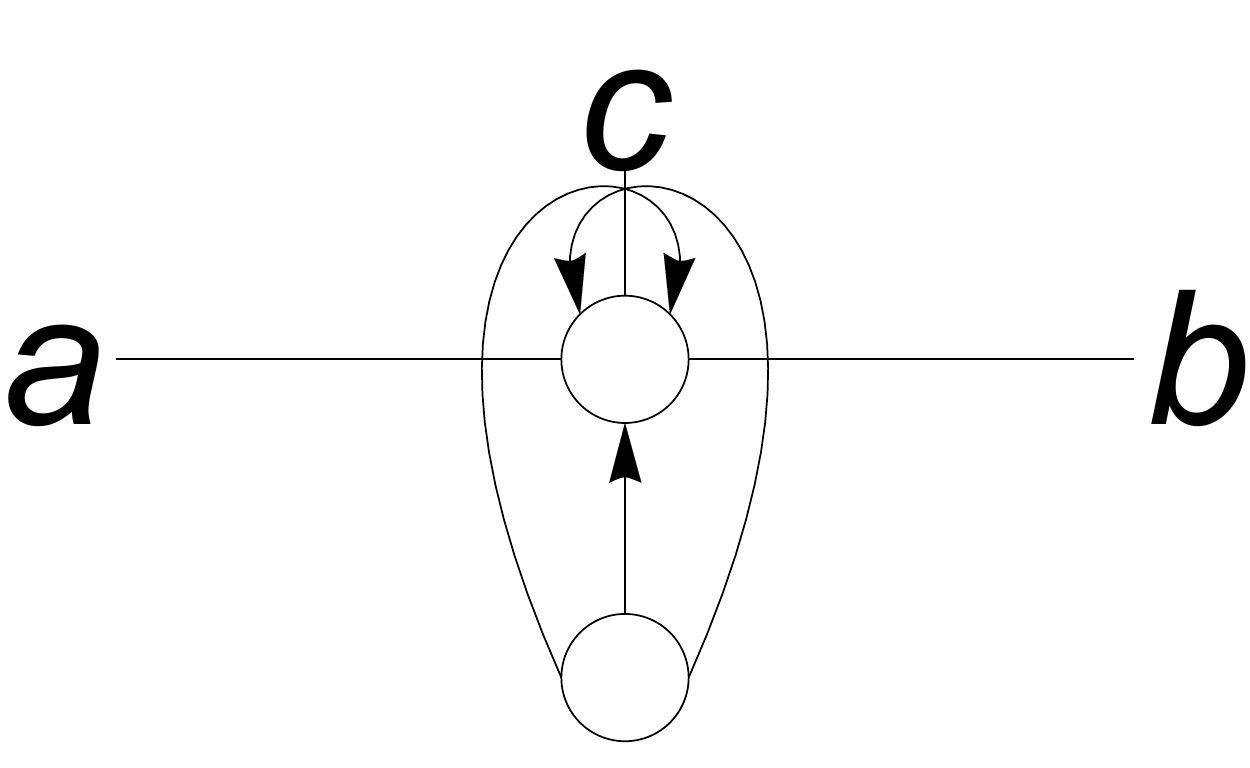}
}
\vspace{0mm}
\centerline{
\includegraphics[width=20mm]{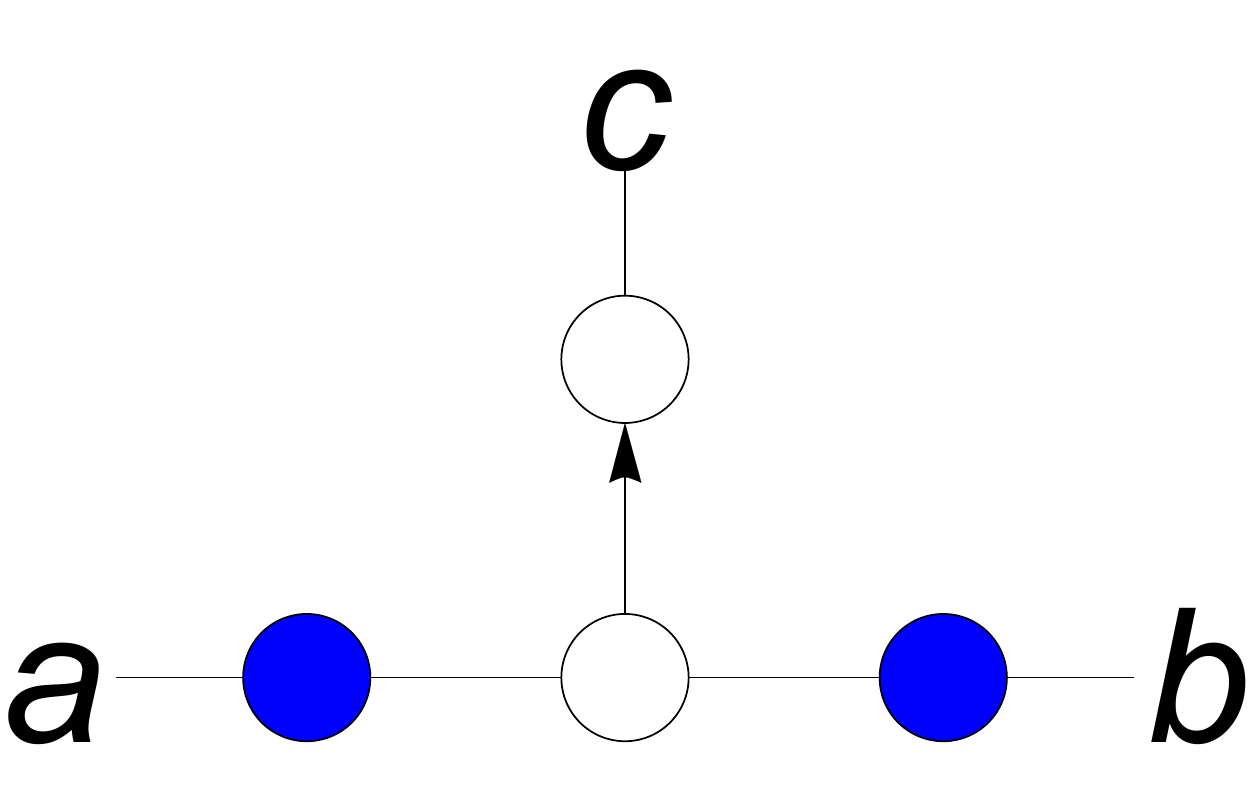}\hspace{5mm}
\includegraphics[width=20mm]{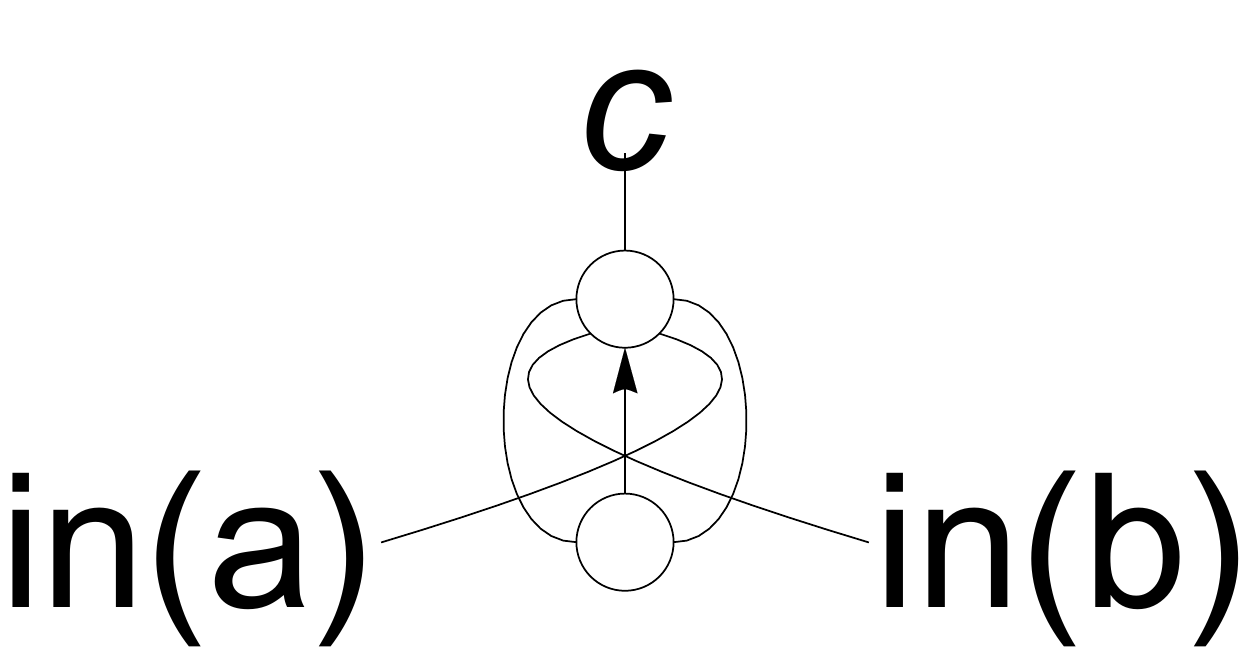}
}
\vspace{5mm}
\centerline{
\includegraphics[height=11mm]{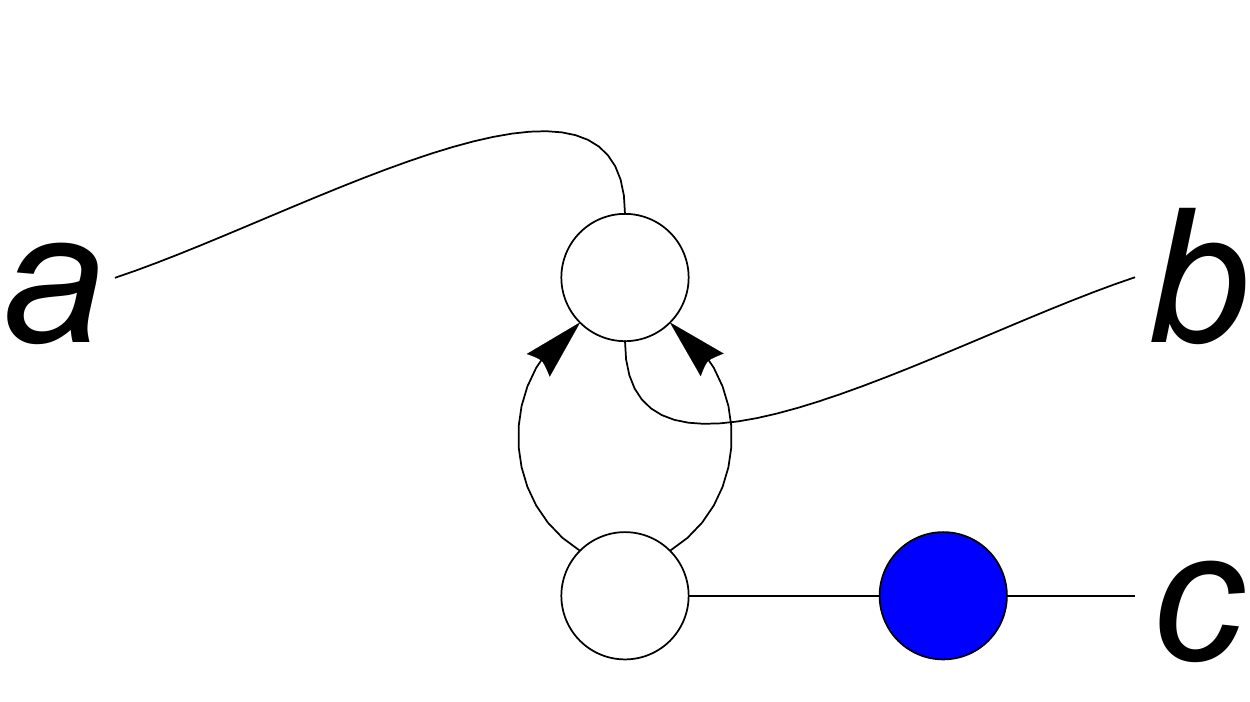}\hspace{5mm}
\includegraphics[height=10mm]{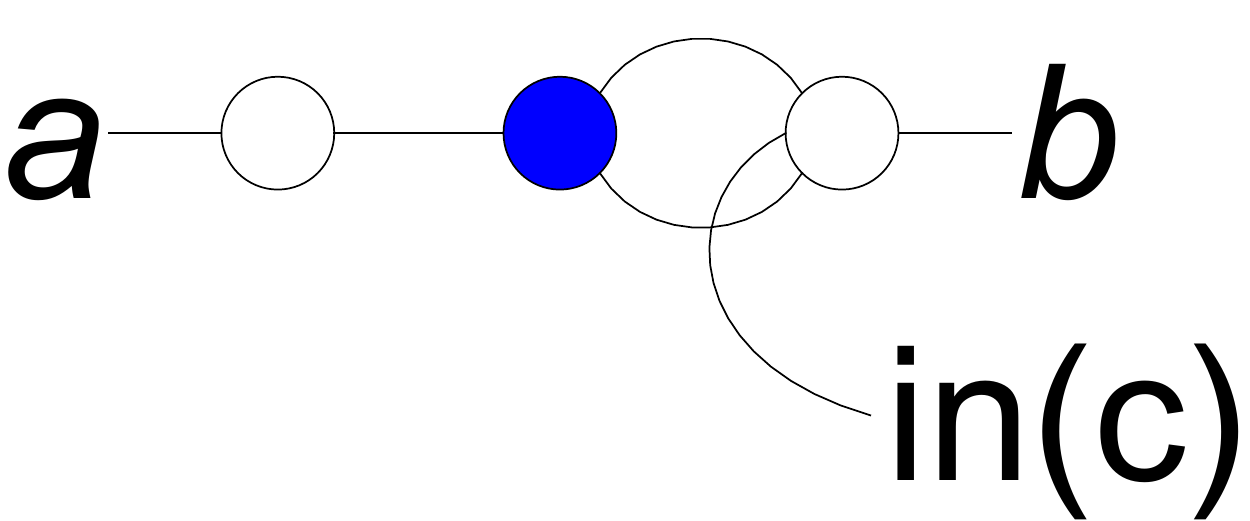}
}
\vspace{5mm}
\centerline{
\includegraphics[height=11mm]{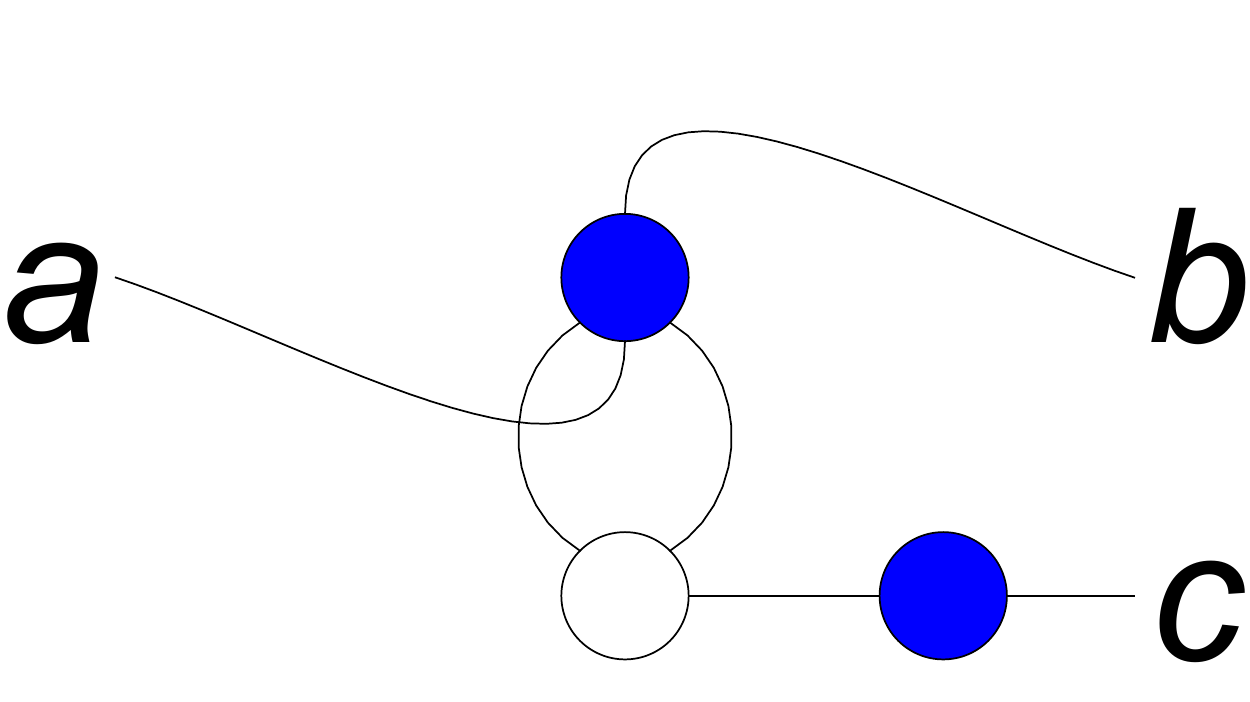}\hspace{5mm}
\includegraphics[height=10mm]{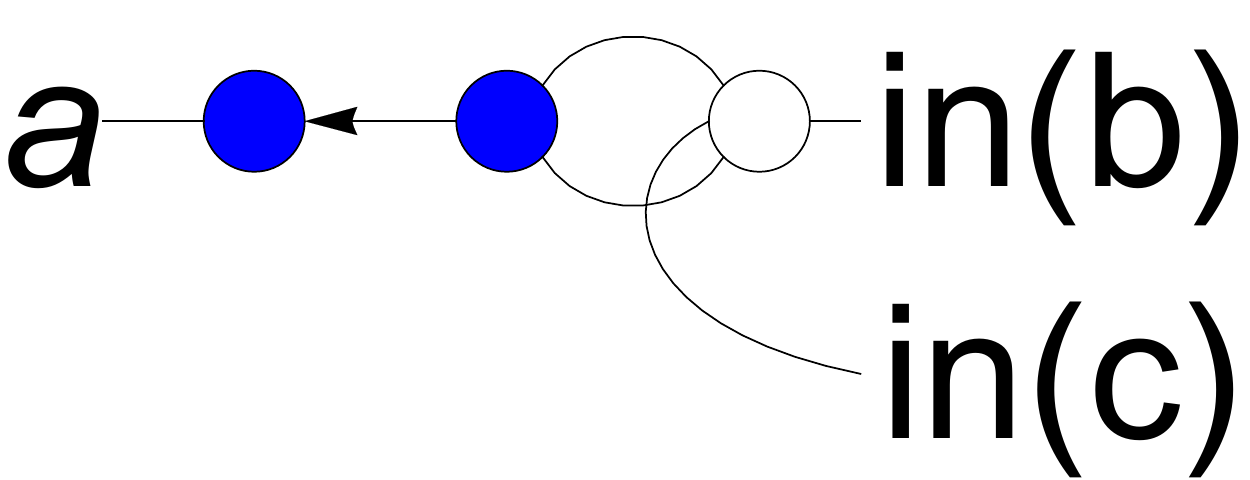}
}
\vspace{11mm}
\centerline{
\includegraphics[height=8mm]{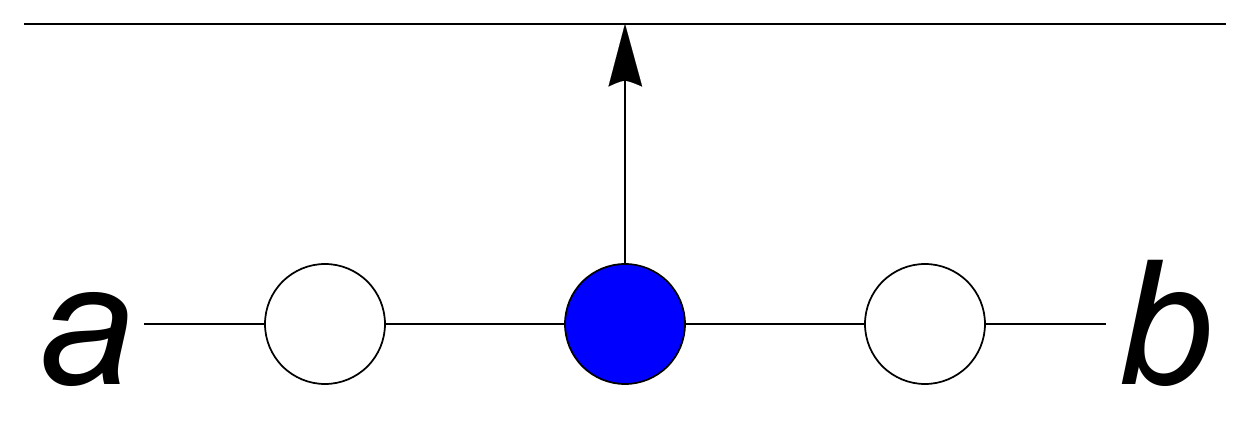}\hspace{5mm}
\includegraphics[height=8mm]{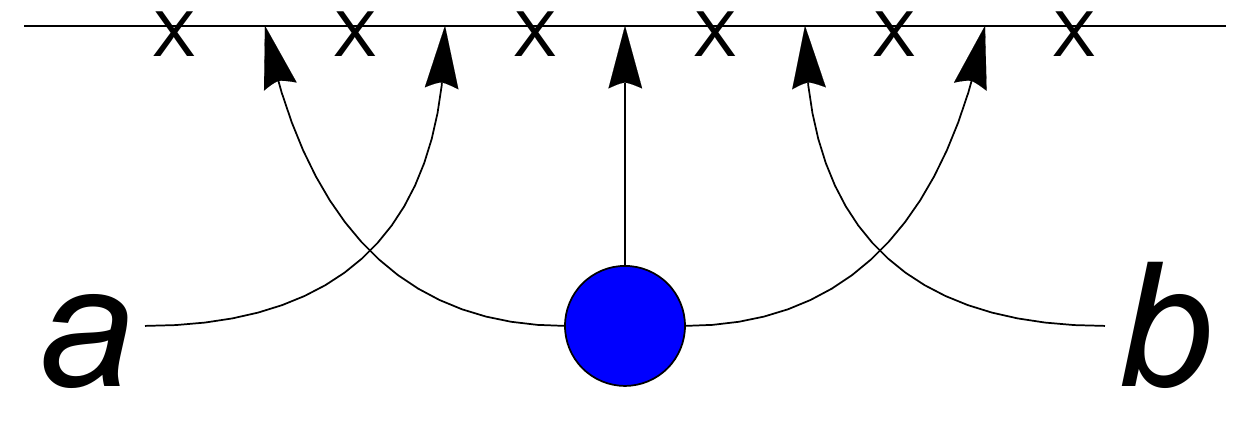}
}
\vspace{11mm}
\centerline{
\includegraphics[height=8mm]{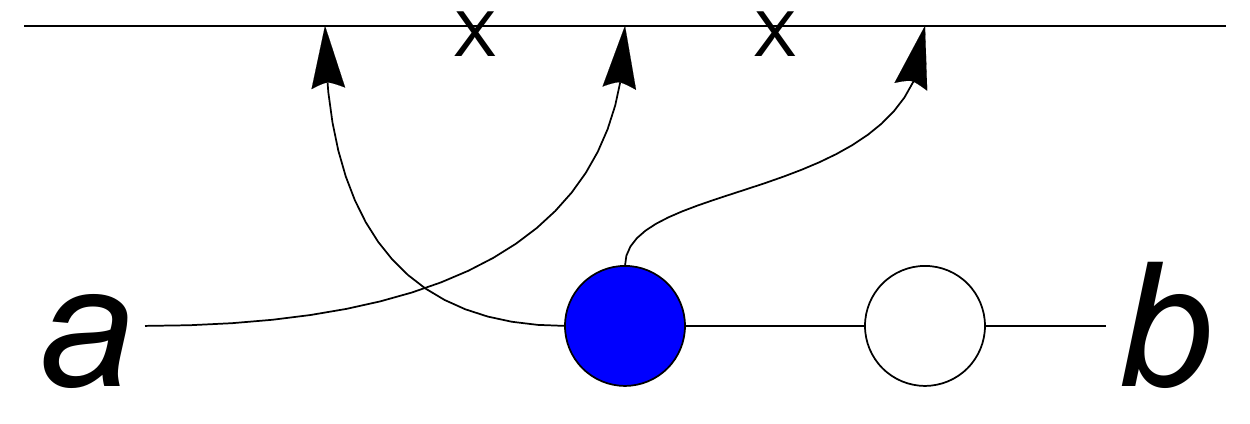}\hspace{5mm}
\includegraphics[height=8mm]{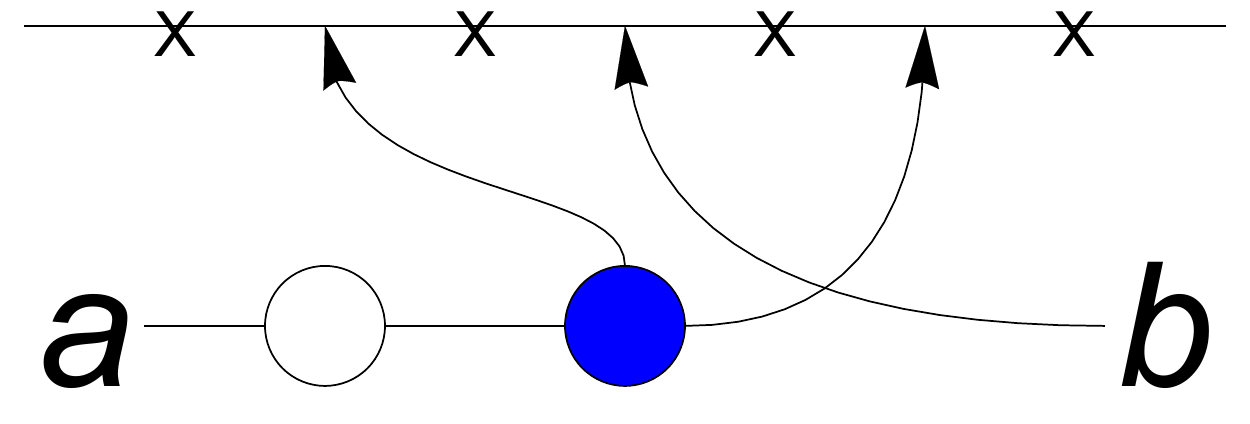}
}
\vspace{11mm}
\centerline{
\includegraphics[height=9mm]{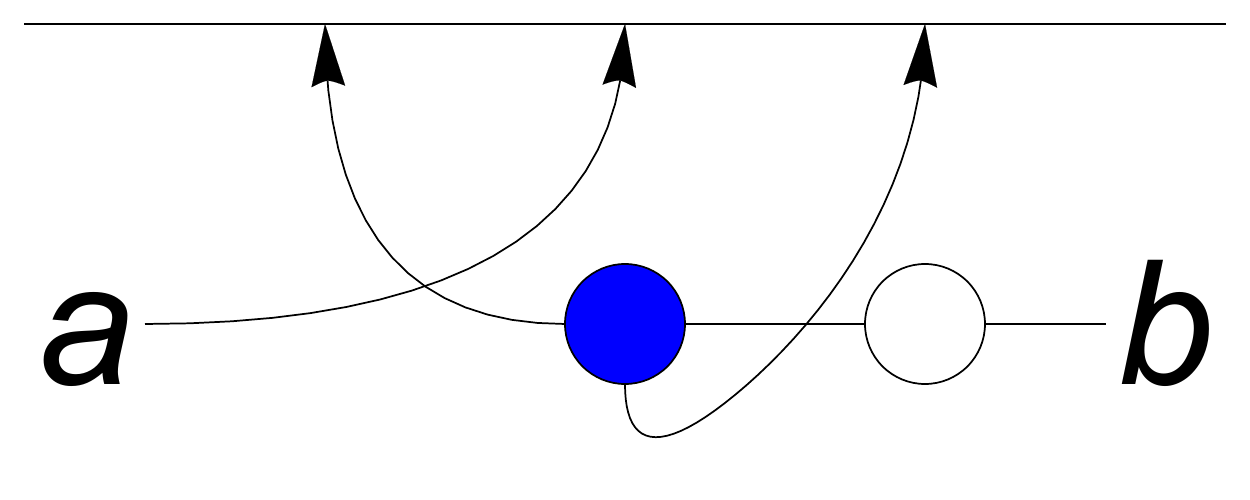}\hspace{5mm}
\includegraphics[height=9mm]{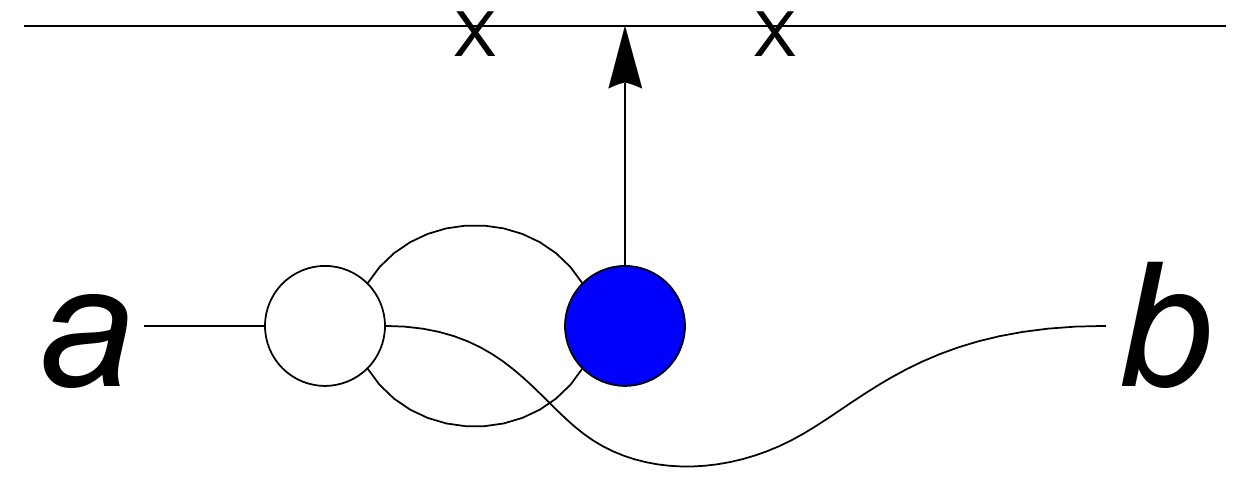}
}
\vspace{8mm}
\centerline{
\includegraphics[height=8mm]{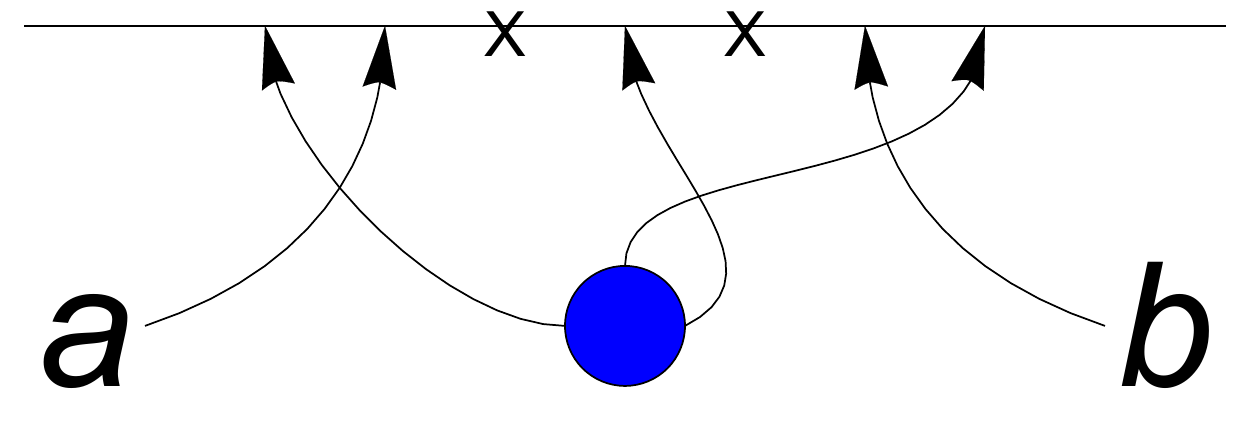}\hspace{5mm}
\includegraphics[height=8mm]{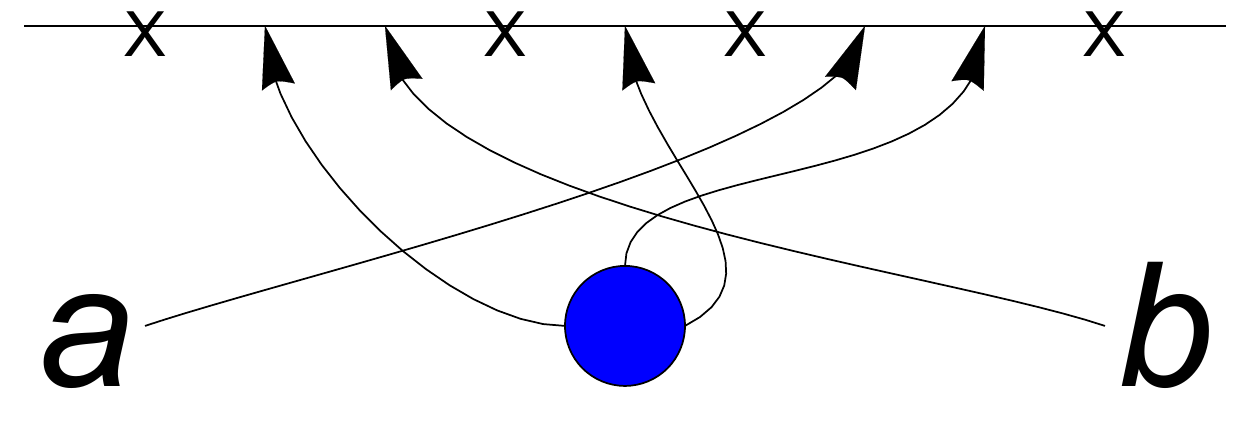}
}
\caption{Strong triple point perestroikas $\Gamma(K)$
\label{R3-2}}
\end{minipage}
\end{figure}

\begin{remark}
Figures \ref{R2-1}, \ref{R2-2}, \ref{R3-1}, \ref{R3-2}
show the Reidemeister moves in terms of the graph $\Gamma(K)$.
Here and in the following, the edges in the root cluster connecting
vertices of the same color are marked with the symbol $\infty$
{\rm (}which expresses that the corresponding vanishing cycles intersect
the pseudoline $J$ ``at infinity'' which we chose in order to define
the colors{\rm )}.
\end{remark}

\begin{figure}[ht]
\begin{minipage}{0.45\textwidth}
\centerline{
\includegraphics[height=3mm]{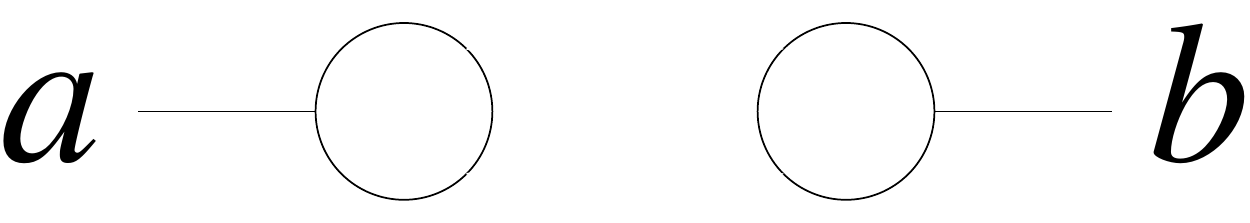}\hspace{18mm}
\includegraphics[height=3mm]{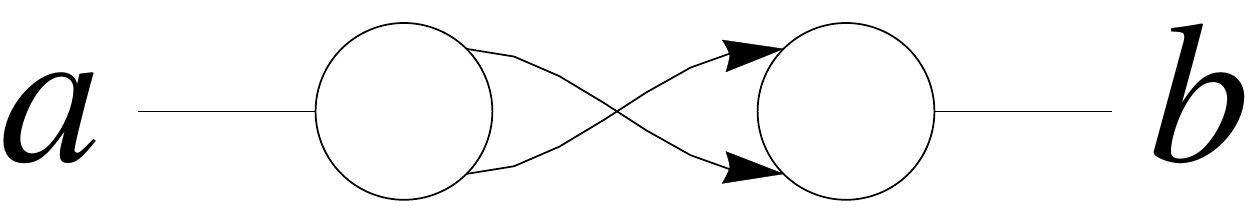}
}
\vspace{7mm}
\centerline{
\includegraphics[height=9mm]{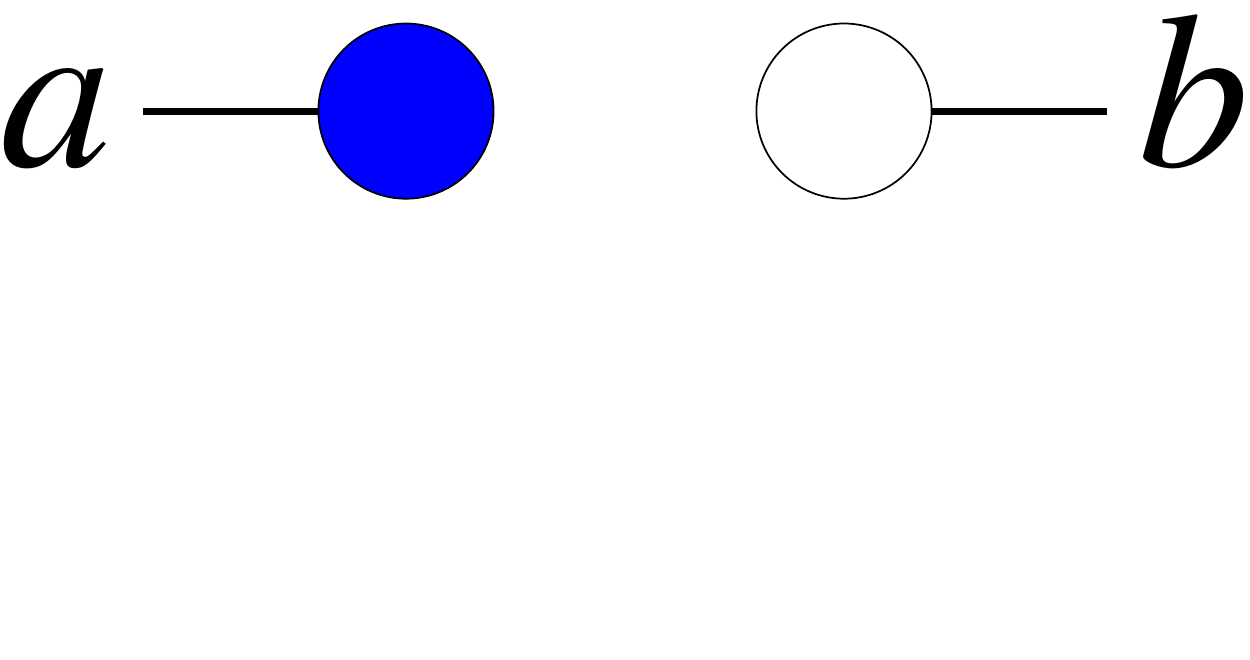}\hspace{18mm}
\includegraphics[height=9mm]{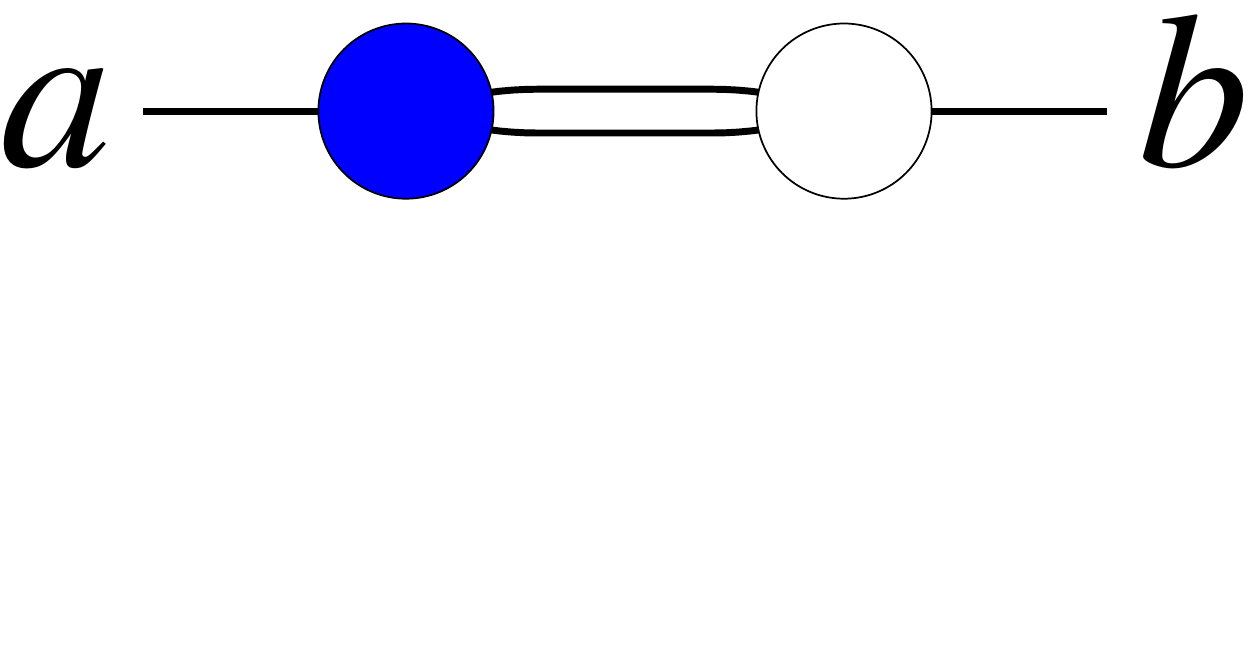}
}
\vspace{0mm}
\centerline{
\includegraphics[height=6.5mm]{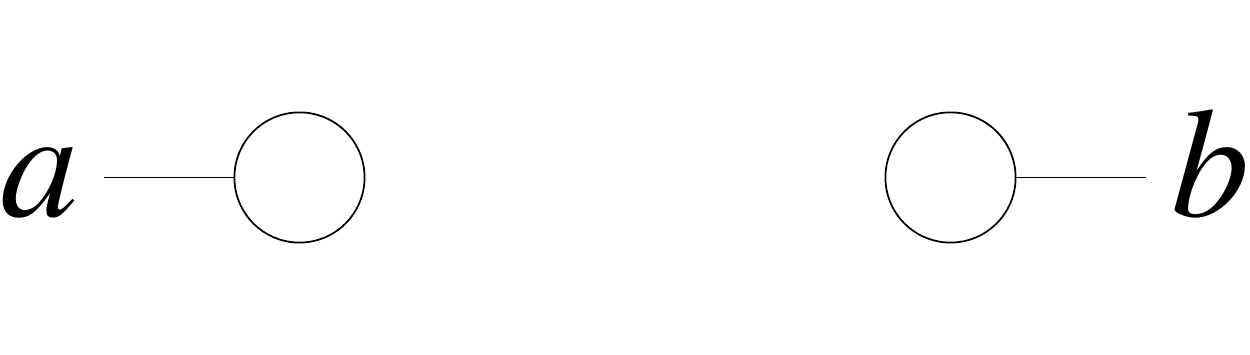}\hspace{11mm}
\includegraphics[height=6.5mm]{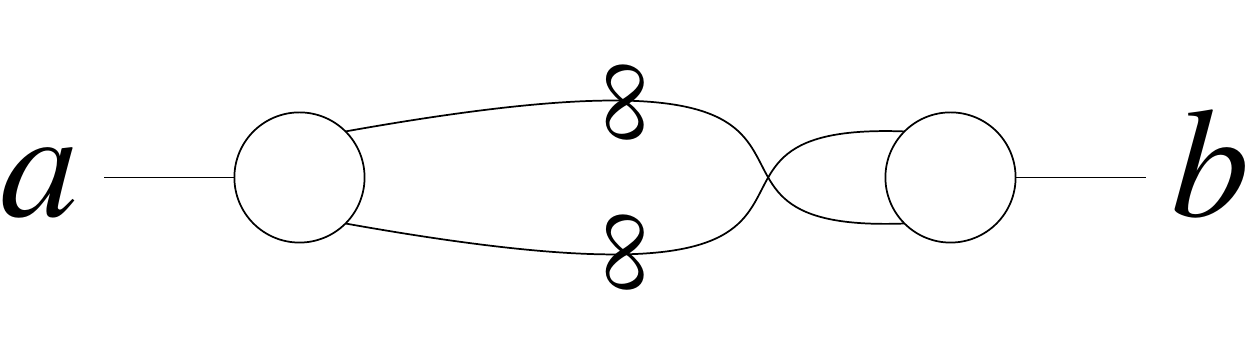}
}
\vspace{10mm}
\centerline{
\includegraphics[height=13mm]{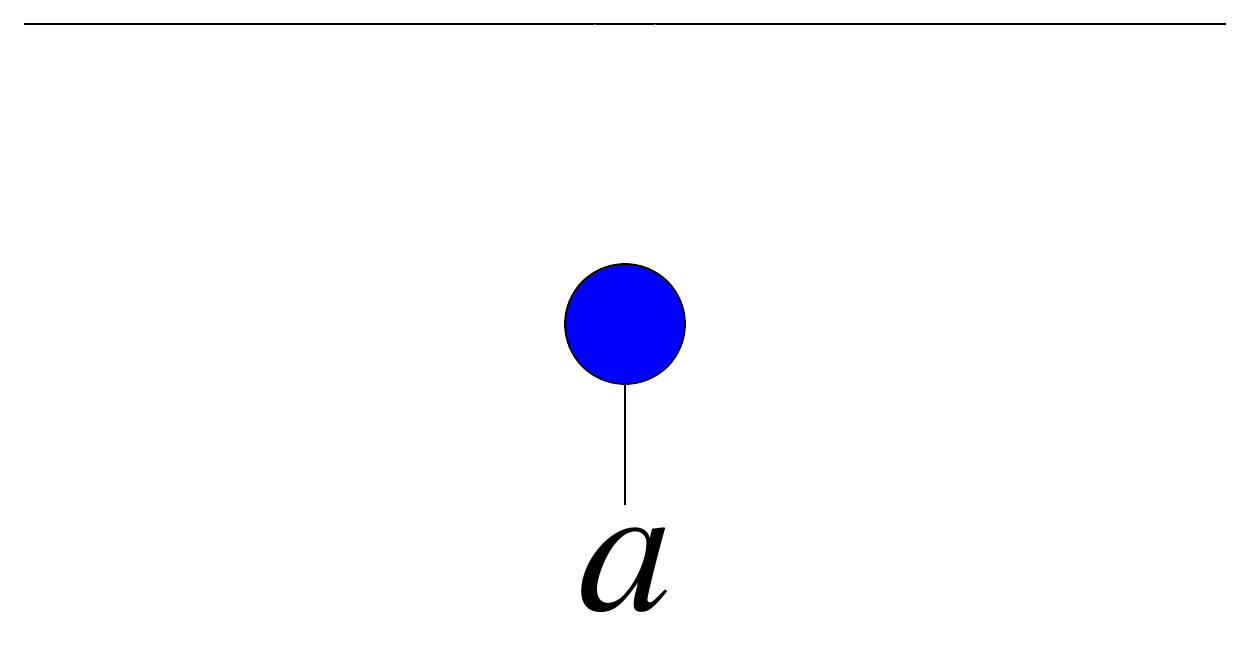}\hspace{8mm}
\includegraphics[height=13mm]{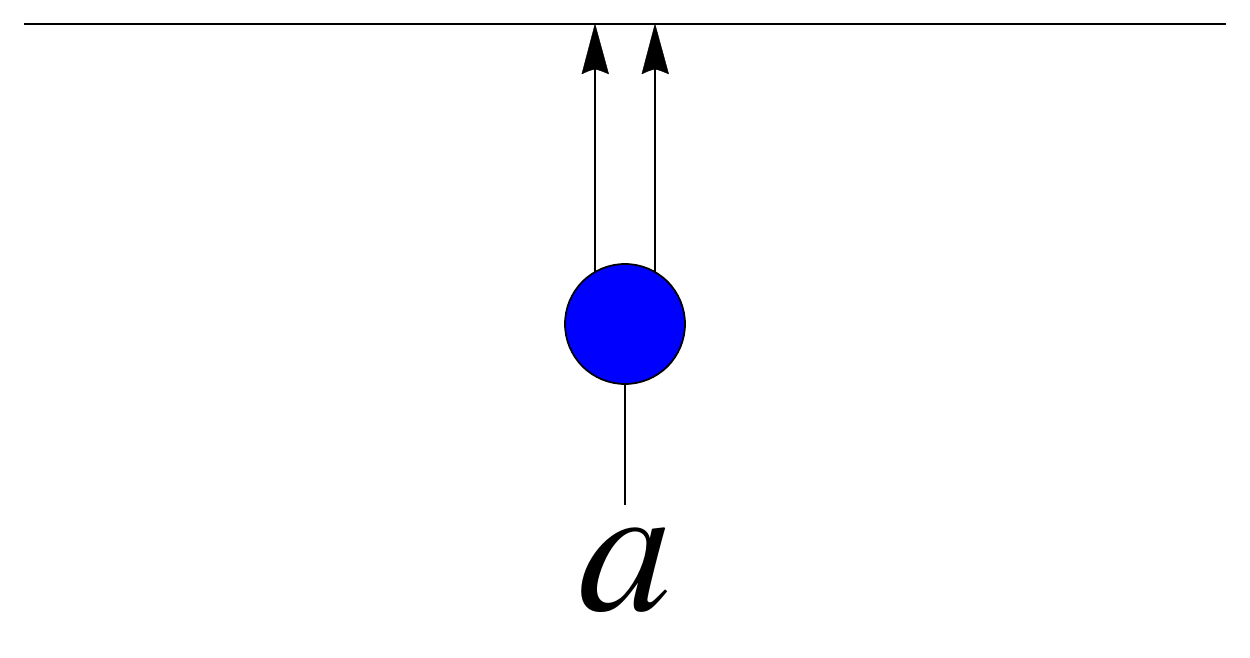}
}

\vspace{0mm}
\centerline{\hspace{-5mm}
\includegraphics[height=31mm]{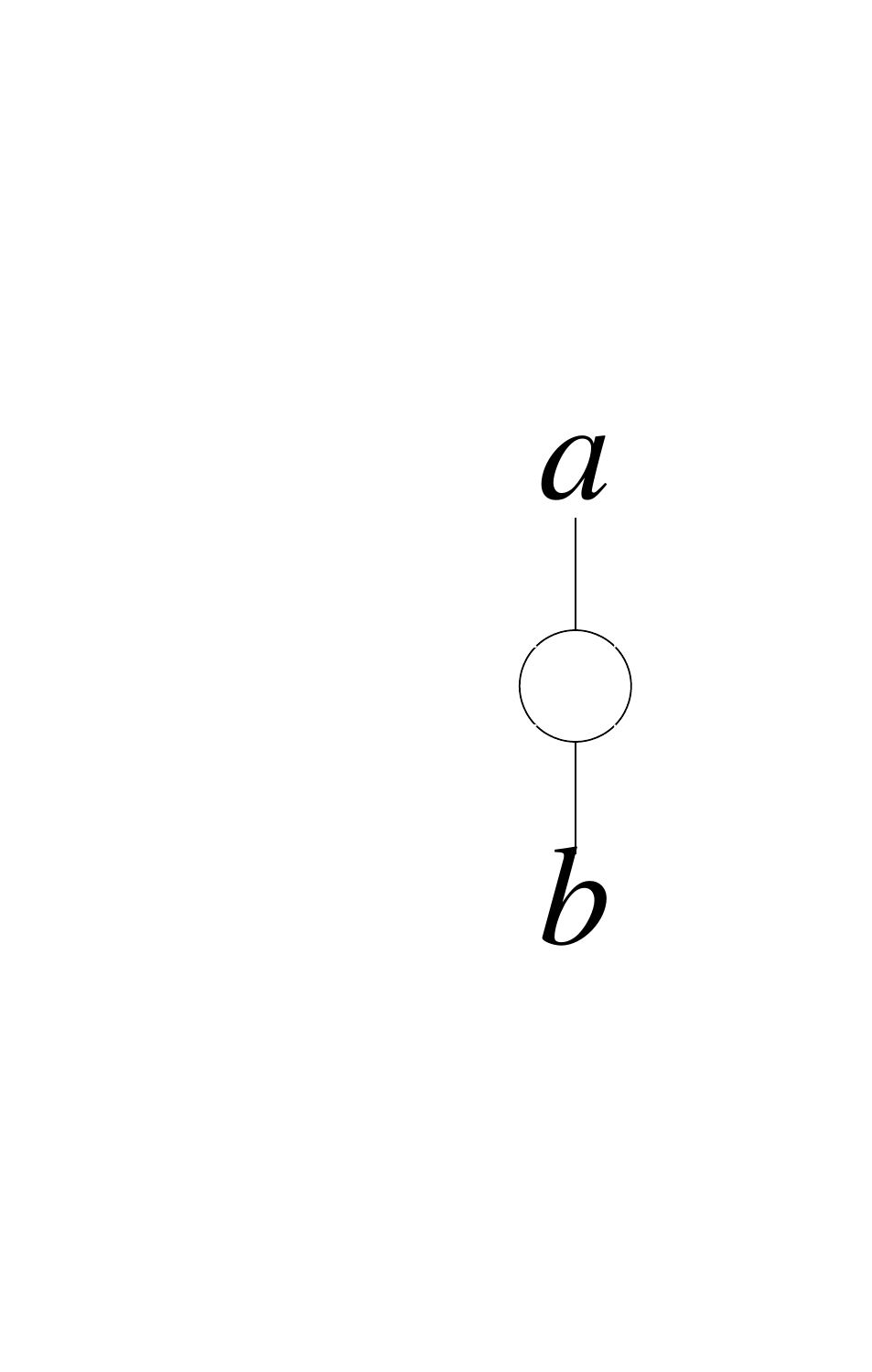}\hspace{20mm}
\includegraphics[height=31mm]{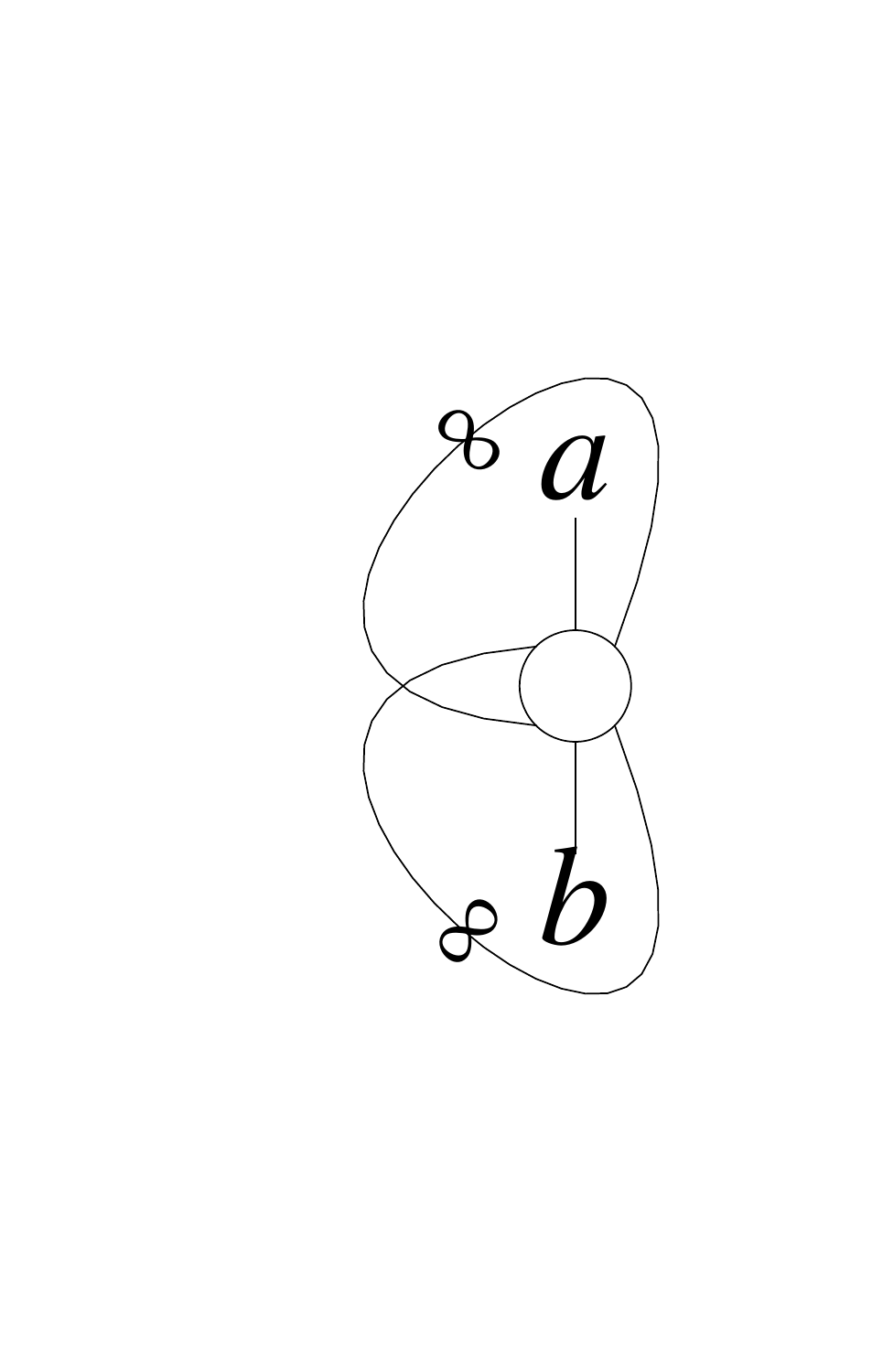}
}
\caption{Direct self-tangency perestroika for $\Gamma(K)$}
\label{R2-1}
\end{minipage}
\hspace{0.10\textwidth}
\begin{minipage}{0.45\textwidth}
\centerline{
\includegraphics[height=3mm]{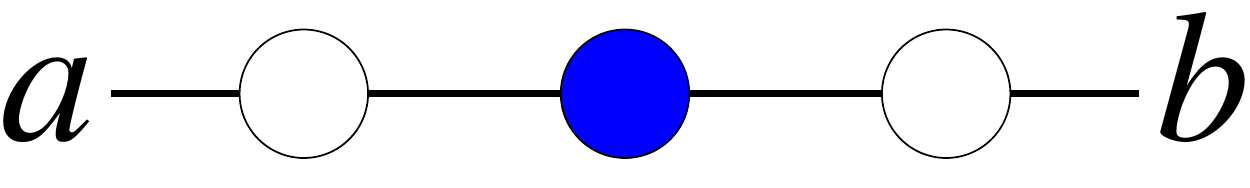}\hspace{19mm}
\includegraphics[height=3mm]{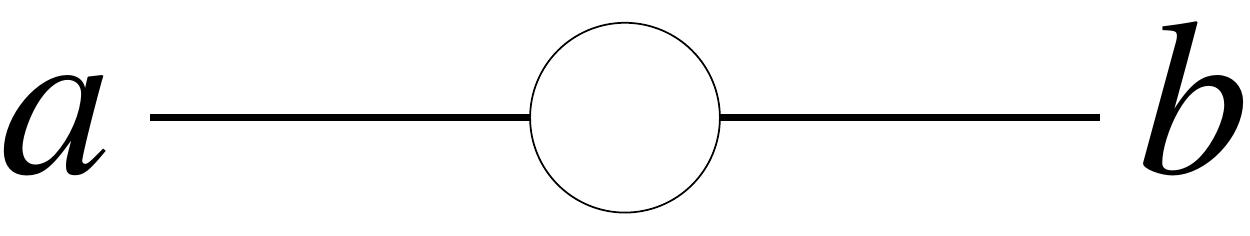}
}
\vspace{-.5mm}
\centerline{
\includegraphics[height=9mm]{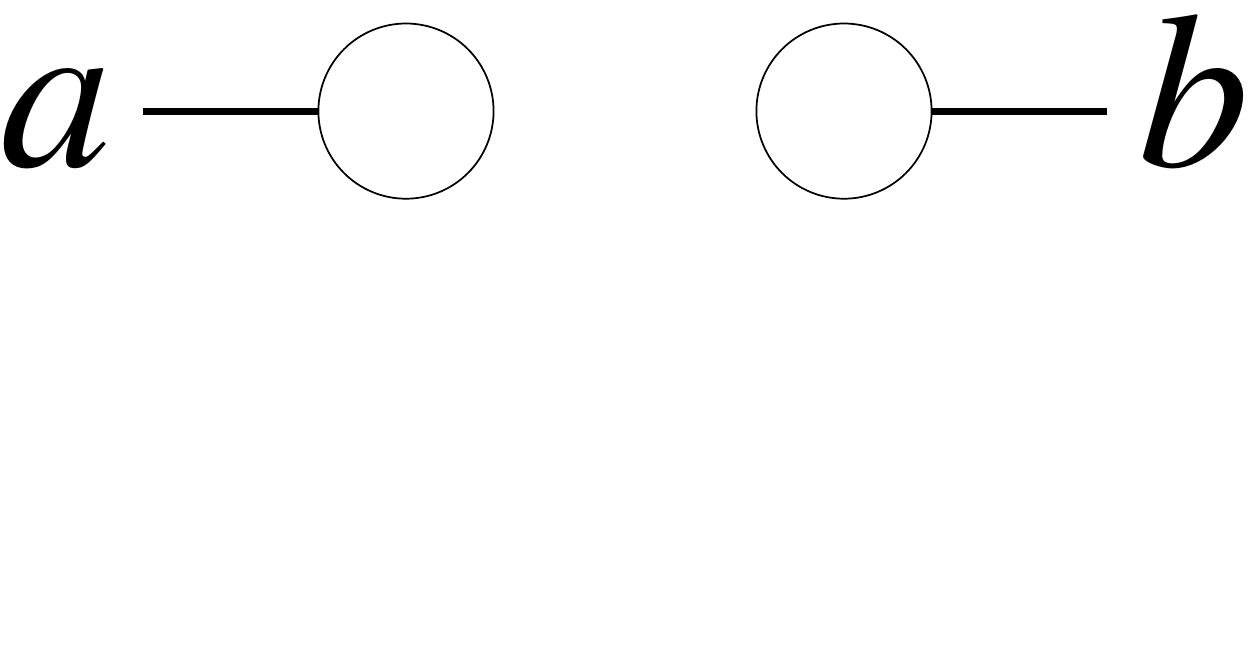}\hspace{20mm}
\includegraphics[height=20mm]{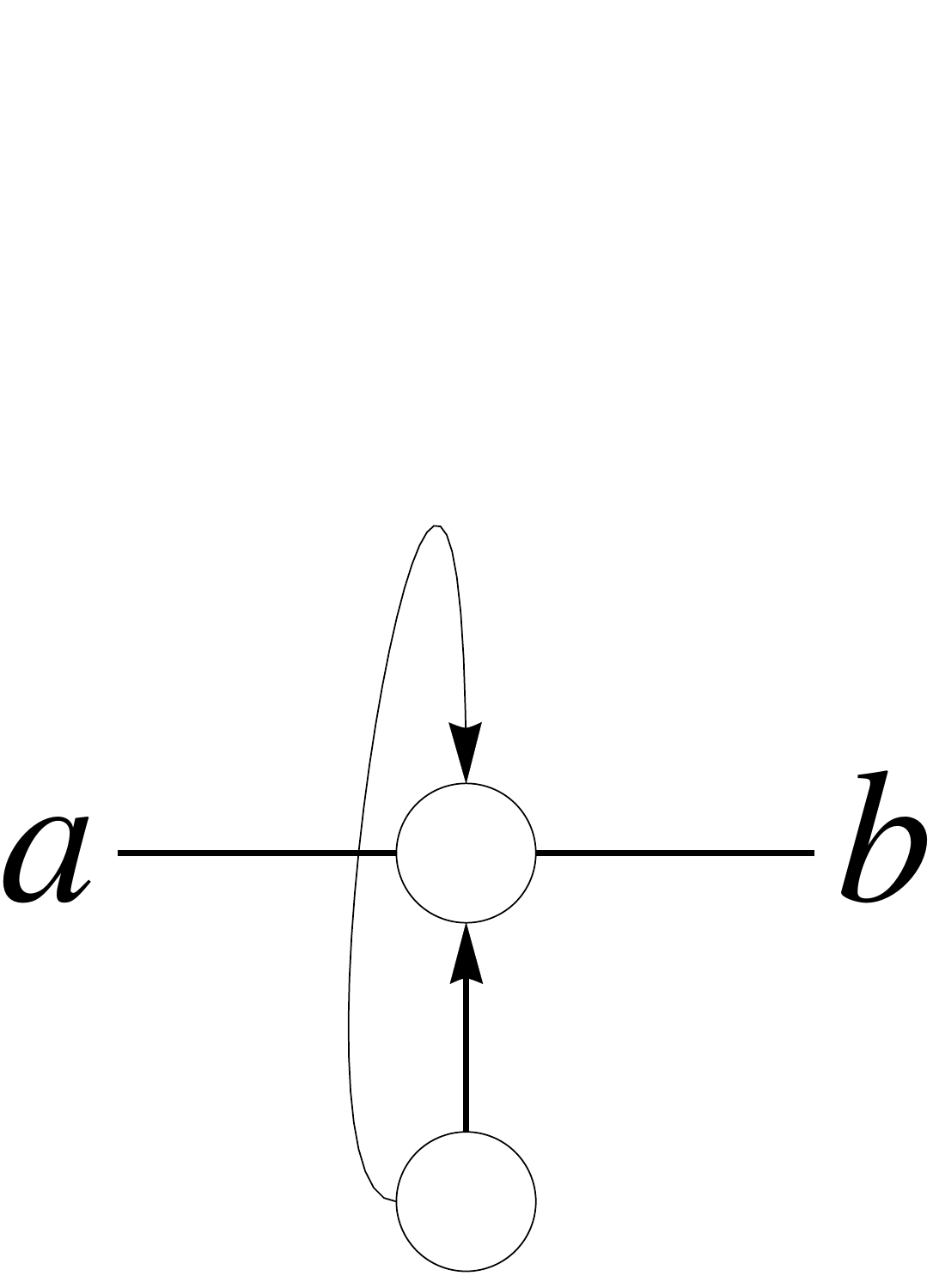}
}
\vspace{3mm}
\centerline{\hspace{5mm}
\includegraphics[height=6mm]{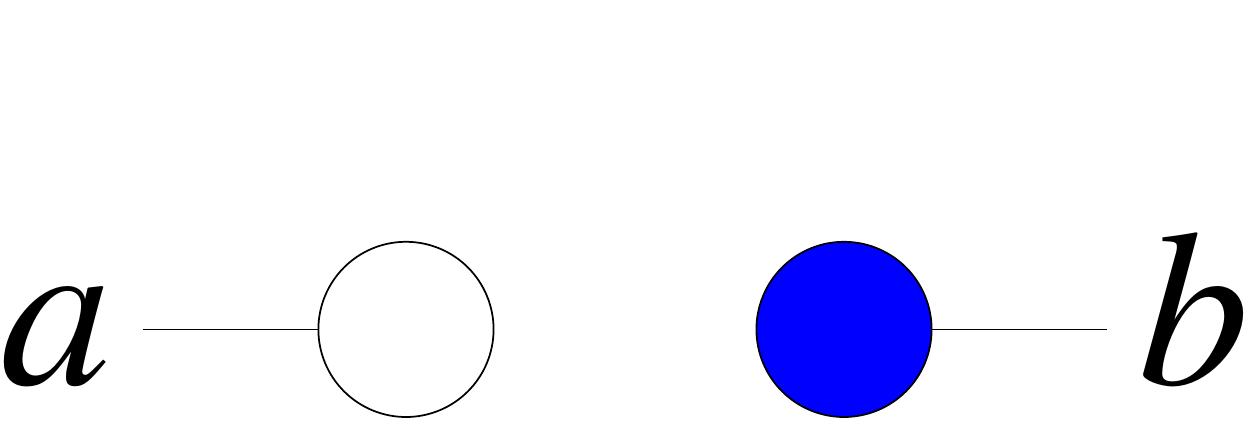}\hspace{15mm}
\includegraphics[height=6mm]{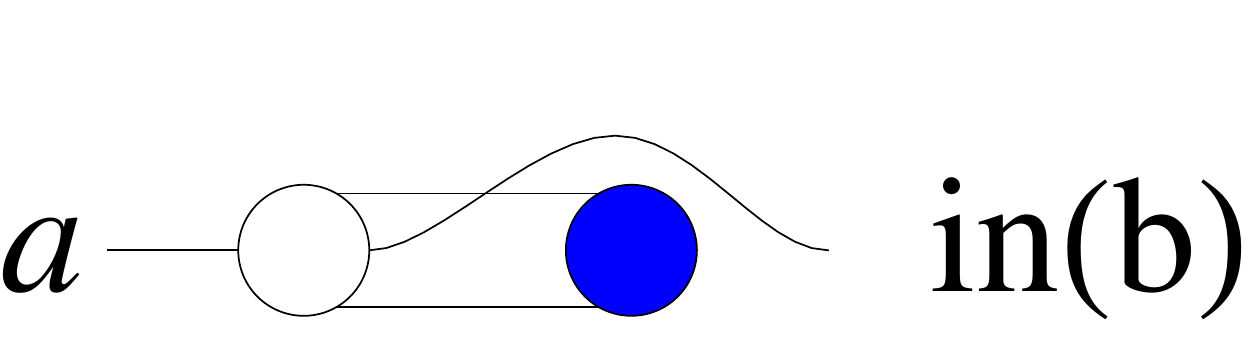}
}
\vspace{10mm}
\centerline{\hspace{5mm}
\includegraphics[height=11mm]{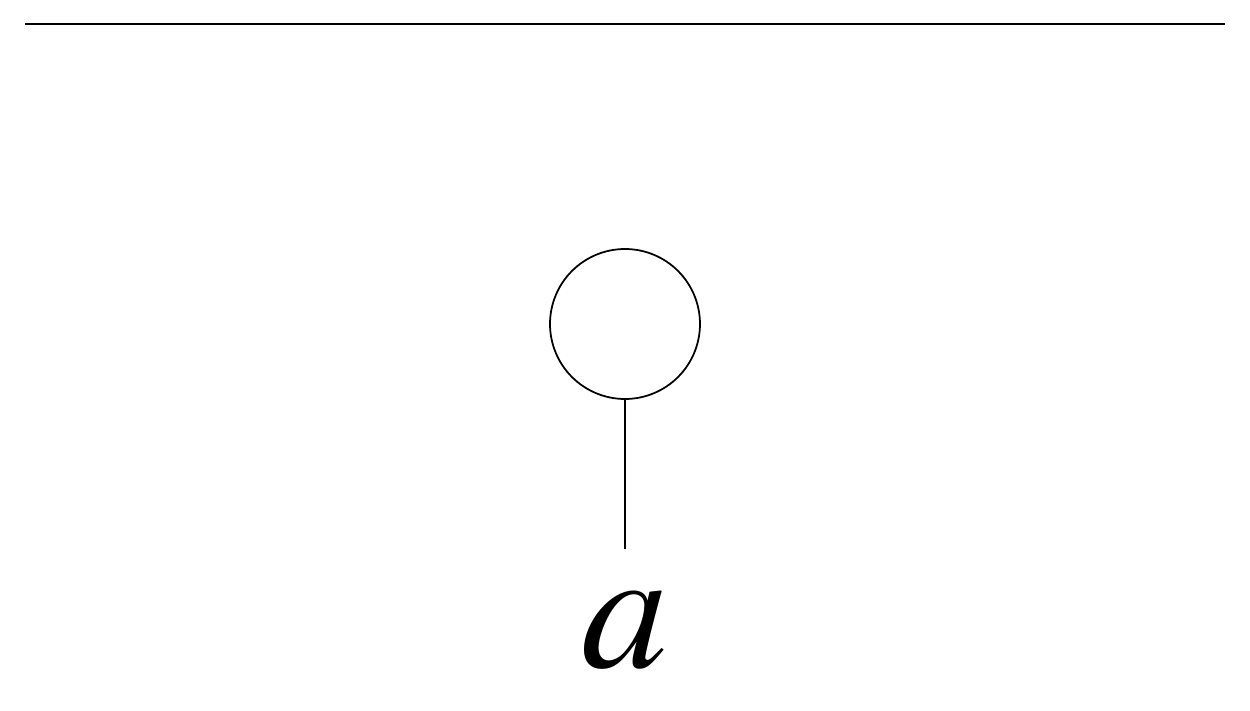}\hspace{15mm}
\includegraphics[height=11mm]{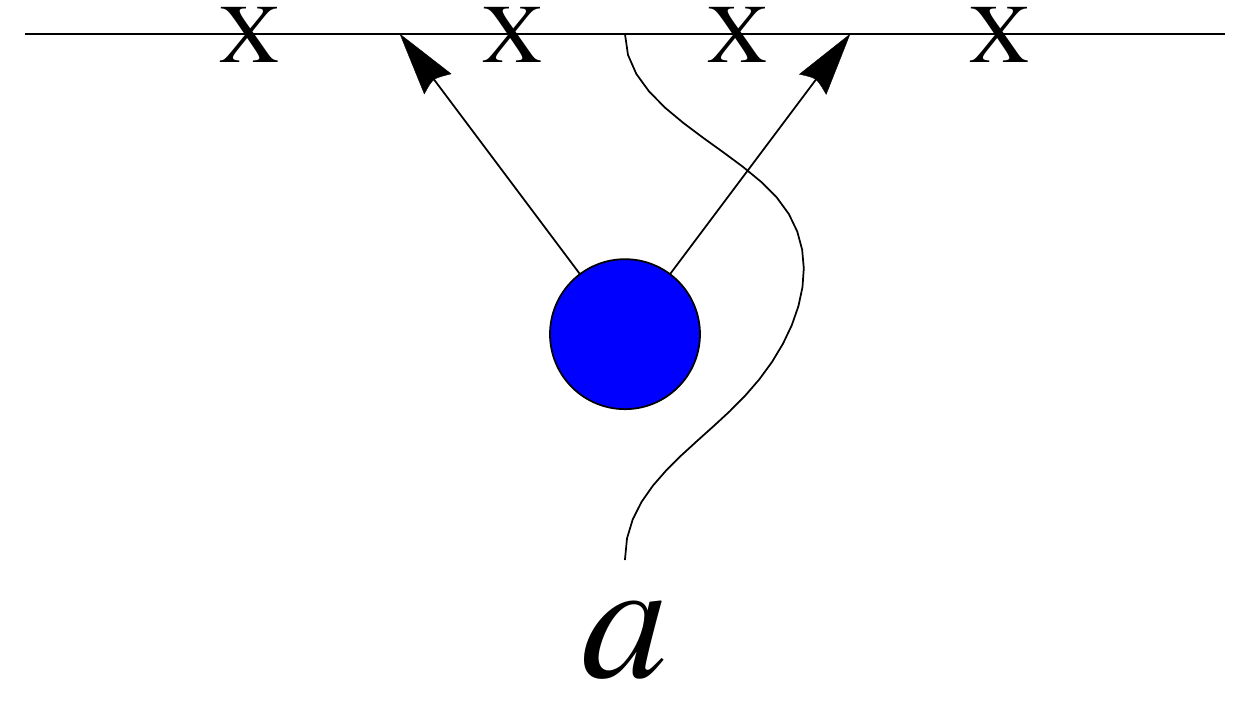}
}
\vspace{6mm}
\centerline{\hspace{7mm}
\includegraphics[height=13mm]{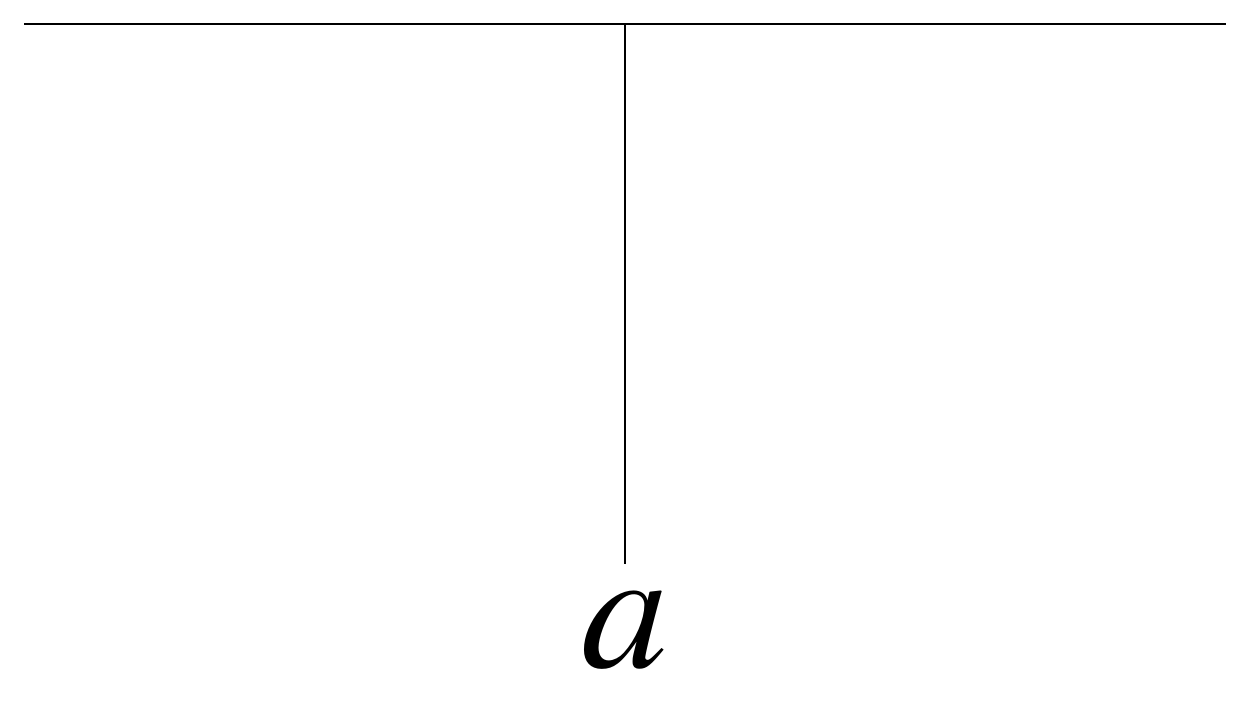}\hspace{12mm}
\includegraphics[height=13mm]{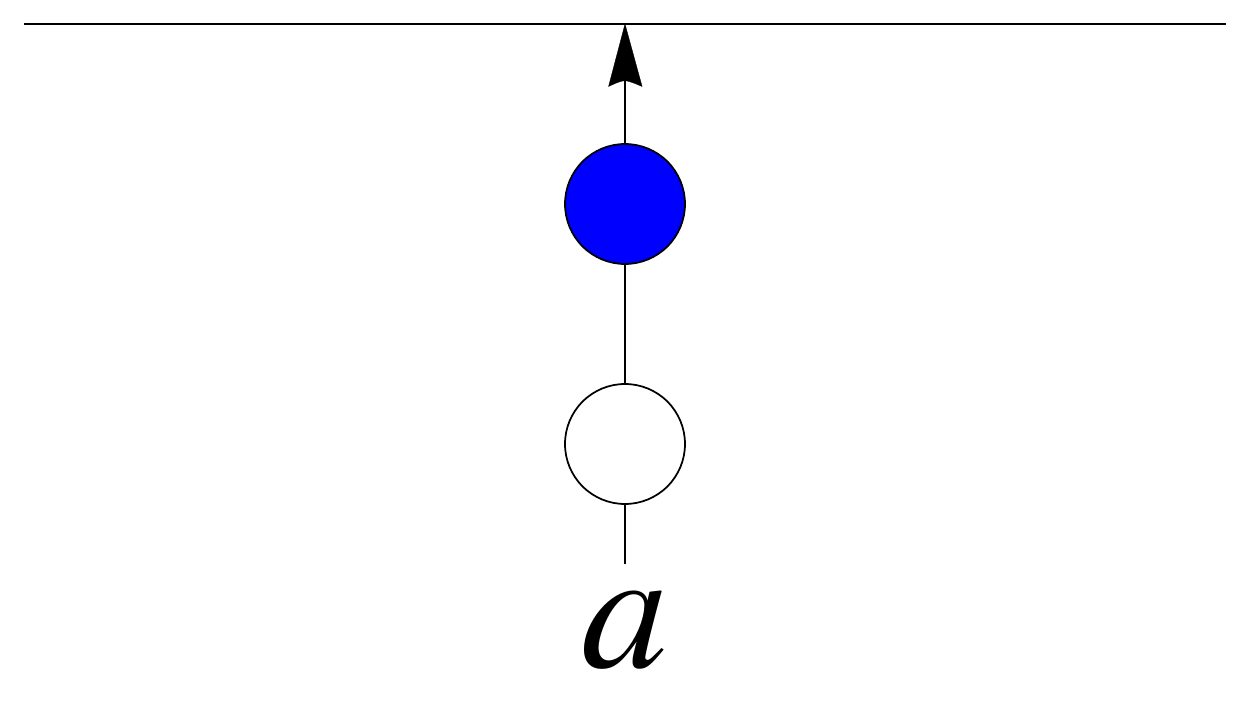}
}
\vspace{8mm}
\caption{Inverse self-tangency perestroikas for $\Gamma(K)$
\label{R2-2}}
\end{minipage}
\end{figure}

%

The following conventions are adopted.
Each letter $a$, $b$, $c$ stands for a sequence of edges (one after another in the cyclic order of the ribbon graph)
adjacent to a vertex of $\Gamma(K)$. These edges may be oriented  or unoriented.
The symbol $"in(a)"$ stands for the inversion of the sequence $a$, {\it i.e.}, inserting it in the reverse order.
Such an inversion happens every time the sequence gets attached to a vertex of different color after the move.

Also, some moves in Figure \ref{R3-2} change a sequence ($a$ or $b$) from being adjacent to a non-root vertex
(depicted by a blue or white disk) to the root vertex (depicted by the line). Recall that all edges adjacent to the root
vertex are oriented towards it. If an edge in a sequence $a$ (or $b$) adjacent to a non-root vertex was unoriented,
then it becomes oriented after such a move.
If it was already oriented, then it remains oriented, but we add two crosses,
one on each side of its adjacency to the line. Recall that this means that the corresponding vanishing cycle is attached
from the other (local) side of the non-contractible component of the smoothing diagram.

The third and the fifth move in Figure \ref{R2-1},
as well as the fifth move in Figure \ref{R3-1},
involve ovals of the same color connected
with non-oriented edges. These moves are only applicable
to the even $d$ case as the relative interior
of the corresponding vanishing
cycles must intersect the auxiliary pseudoline $J\subset\rp^2$
while being disjoint from $K_\circ$.

\subsection{Classification of generic immersions
of a circle with small \texorpdfstring{$n$}{n}}
If $n(K) = 0$, then our immersion is an embedding.
If $d(K) = 1 \mod 2$,
then the embedding is isotopic to the standard embedding
$\rp^1\subset\rp^2$.
If $d(K) = 0 \mod 2$, then $K$ must be the standard embedding of a circle
into $\rp^2$ as an oval (the one that bounds an embedded disk).

Immersion graphs provide an exhausting way to classify all
immersions of a circle with a given number $n$ of nodes.
To do that we may list all connected graphs which have $n$ edges
and $m$ vertices with $m = n + 1 \mod 2$ enhanced
as in Definition \ref{enhanced-Gamma}
and extract those that correspond to immersions of a circle.

\begin{figure}[ht]
\begin{minipage}{0.45\textwidth}
\centerline{
\includegraphics[width=31mm]{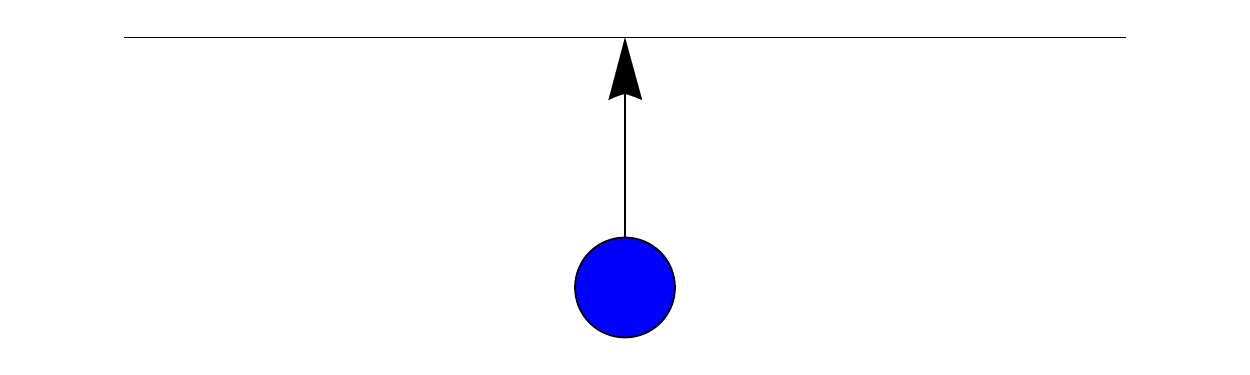}
\hspace{2mm}
\includegraphics[width=30mm]{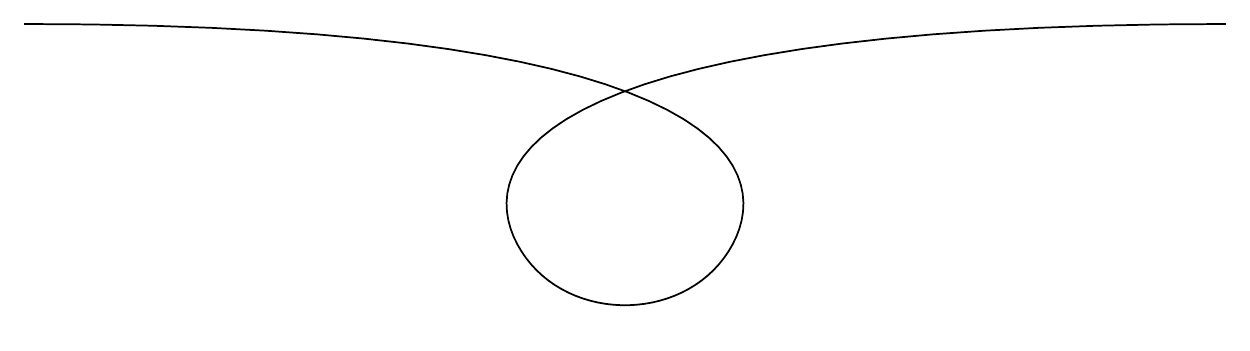}
}
\caption{The graph and the corresponding immersion for
$n=1$, $d=1$
\label{n1classd1}}
\vspace{5mm}
\centerline{
\includegraphics[width=39mm]{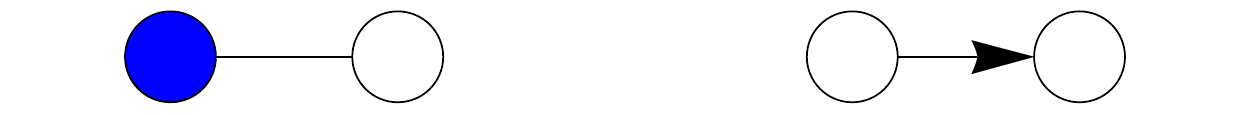}
}
\vspace{5mm}

\centerline{
\includegraphics[width=39mm]{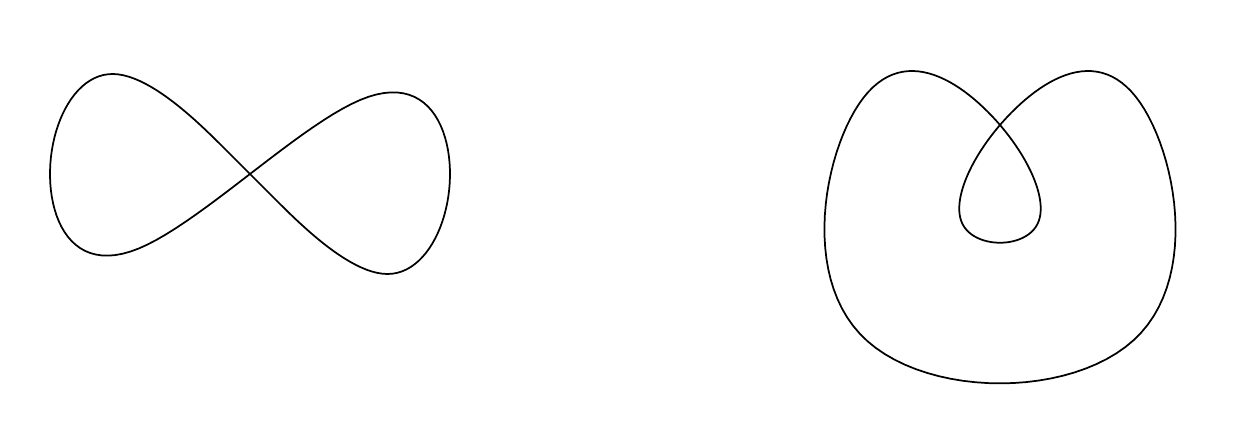}
}
\caption{Graphs and corresponding immersions for $n=1$, $d=0$
\label{n1classd0}}
\end{minipage}
\hspace{0.10\textwidth}
\begin{minipage}{0.45\textwidth}
\centerline{
\includegraphics[width=31mm]{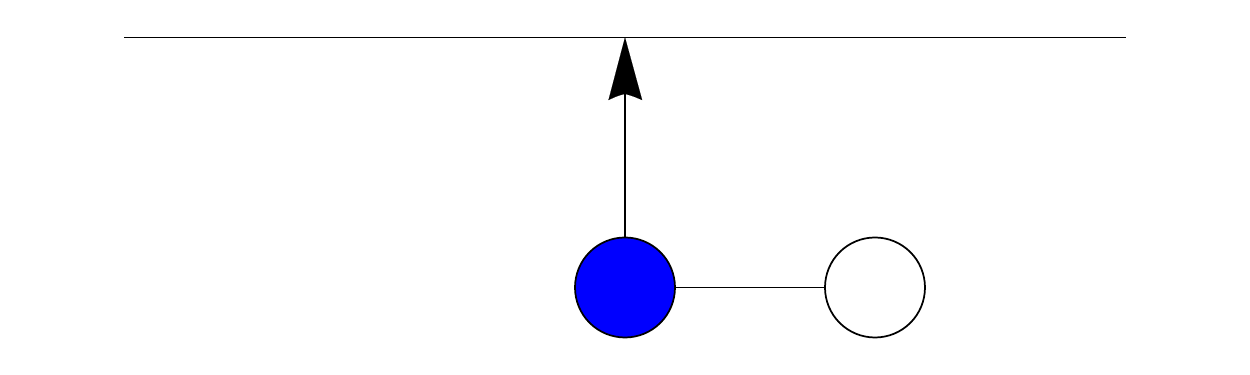}
\hspace{5mm}
\includegraphics[width=30mm]{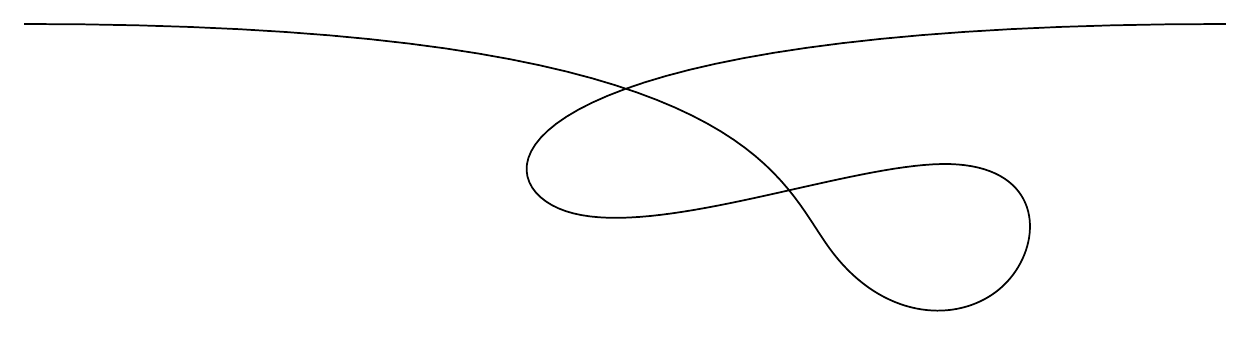}}
\vspace{10mm}
\centerline{
\includegraphics[width=31mm]{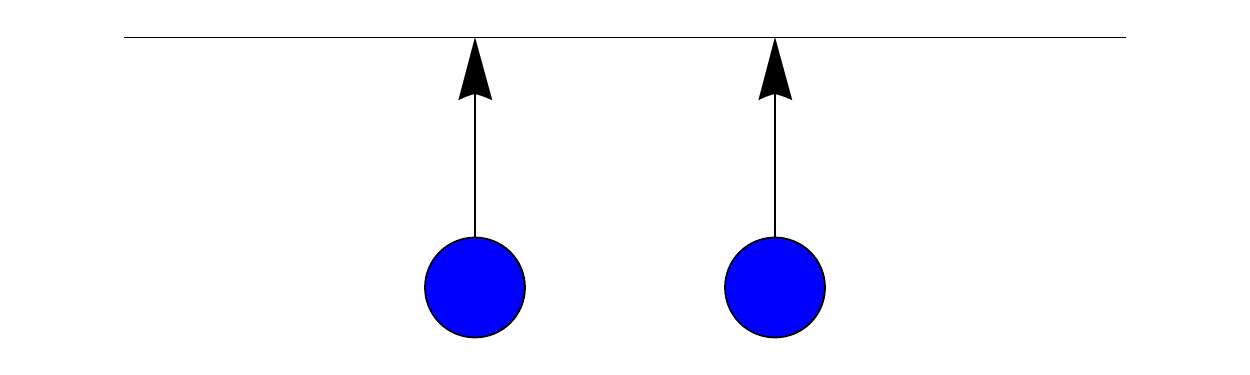}
\hspace{5mm}
\includegraphics[width=30mm]{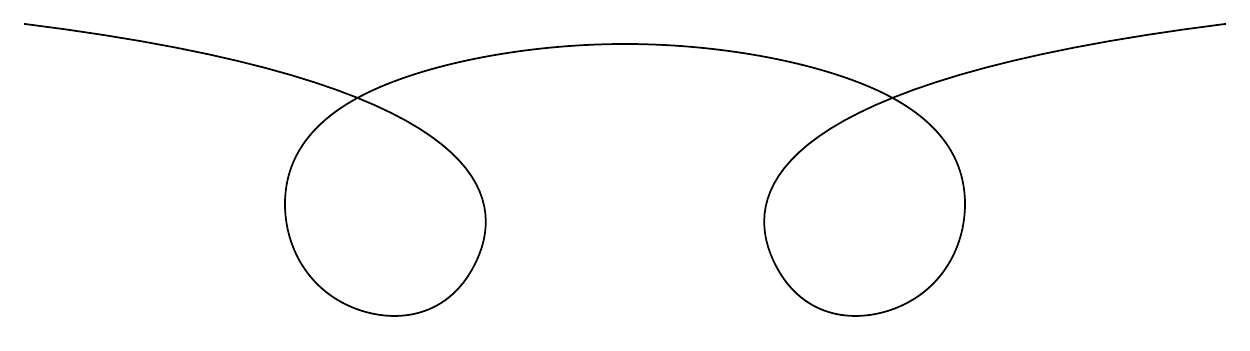}}
\vspace{10mm}
\centerline{
\includegraphics[width=31mm]{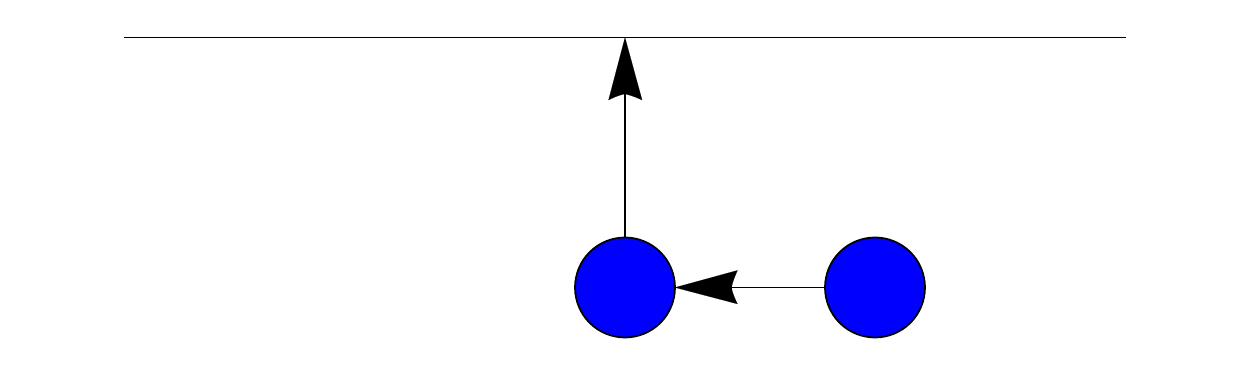}
\hspace{5mm}
\includegraphics[trim = 0mm 20mm 0mm 0mm, clip,width=30mm]{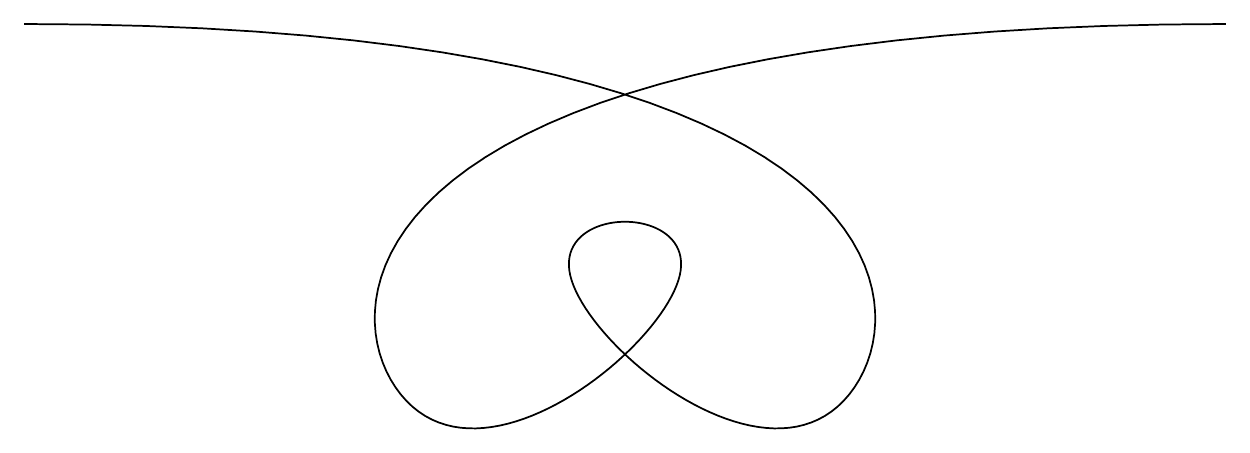}}
\vspace{2mm}
\caption{Graphs and corresponding immersions for
$n=2$, $d=1$
\label{n2classd1}}
\end{minipage}
\end{figure}

If $n=1$, then
the graph $\Gamma(K)$ must have two vertices (as the number of vertices is not greater
than two and has the same parity).
For $d = 1 \mod 2$ there is a unique
graph $\Gamma(K)$, see Figure \ref{n1classd1}.
For $d = 0 \mod 2$ there are two cases,
see Figure \ref{n1classd0} for the graphs and corresponding immersions.

In the case $n=2$, $d=1$ the graph $\Gamma(K)$
must have three vertices (as the root vertex can not be adjacent
to itself if $\Gamma$ comes from an immersion to $\rp^2$).
There are three possibilities depicted (along with
the corresponding immersions) at Figure \ref{n2classd1}.

\begin{figure}[ht]
\centerline{\hspace{-2mm}
\includegraphics[width=18mm]{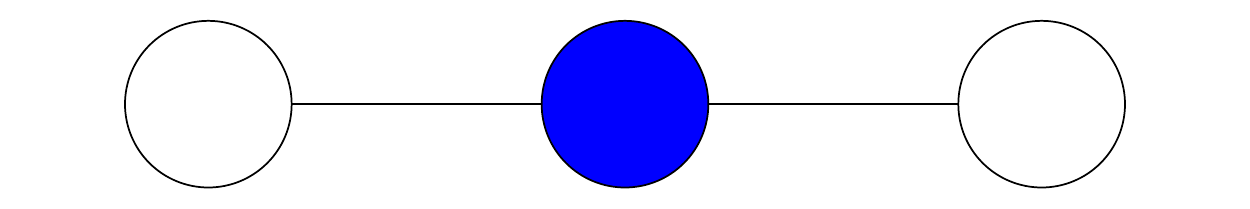}
\hspace{3mm}
\includegraphics[width=18mm]{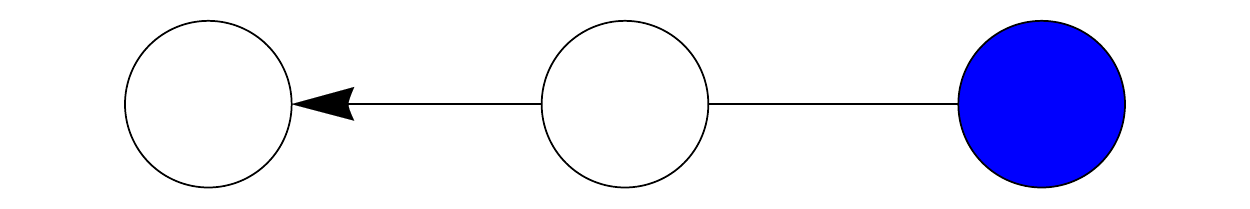}
\hspace{3mm}
\includegraphics[width=18mm]{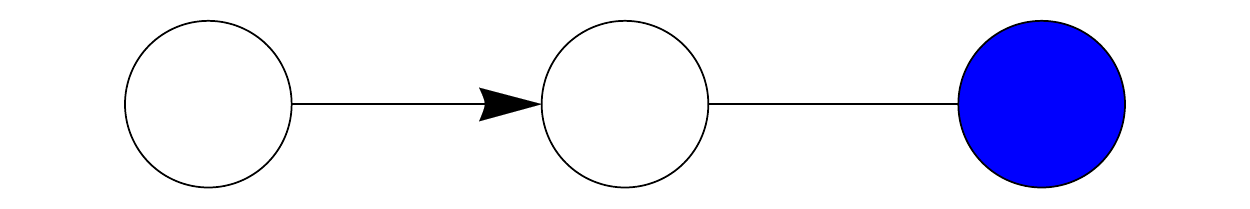}
\hspace{3mm}
\includegraphics[width=18mm]{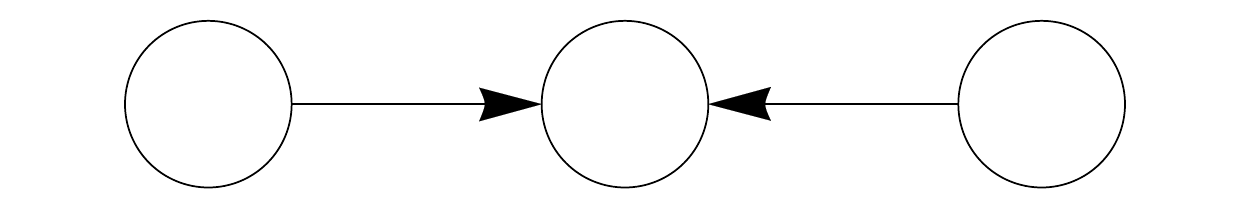}
\hspace{3mm}
\includegraphics[width=18mm]{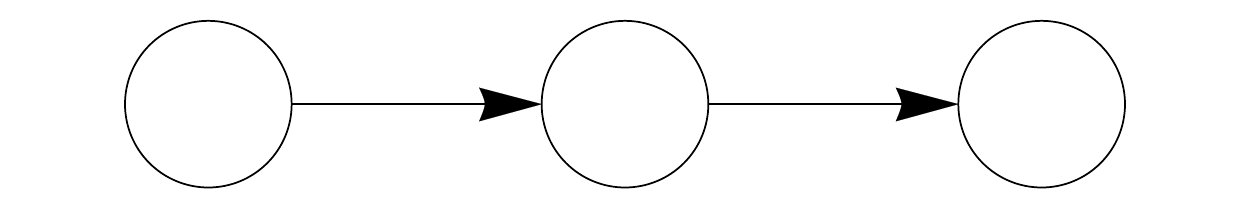}
\hspace{8mm}
\includegraphics[trim = 20mm 23mm 0mm 0mm, clip,width=11mm]{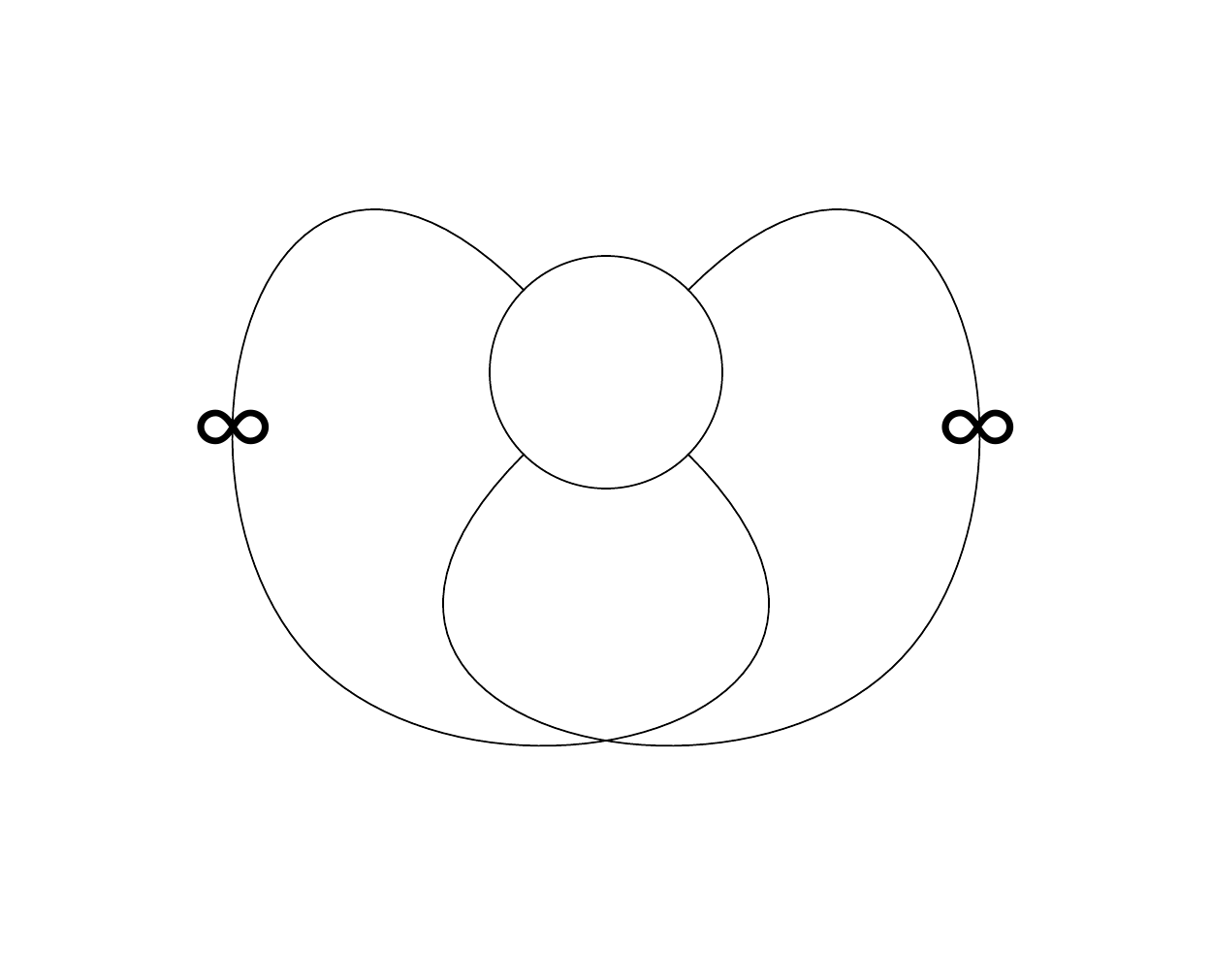}
}
\vspace{10mm}
\centerline{
\includegraphics[width=18mm]{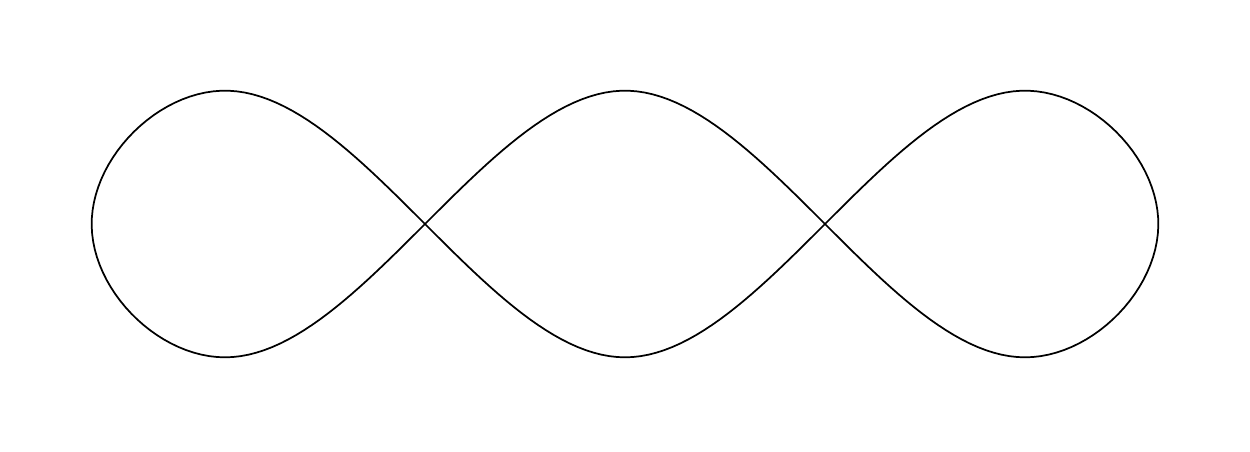}
\hspace{3mm}
\includegraphics[width=18mm]{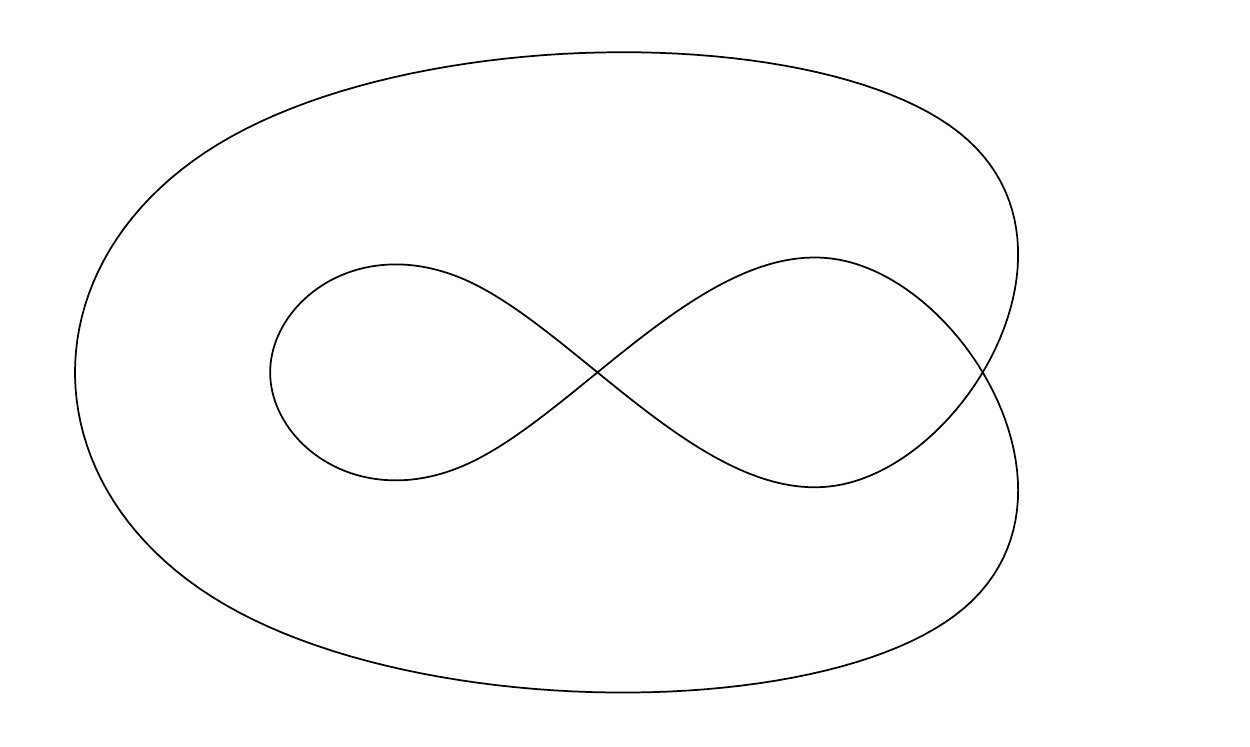}
\hspace{3mm}
\includegraphics[width=18mm]{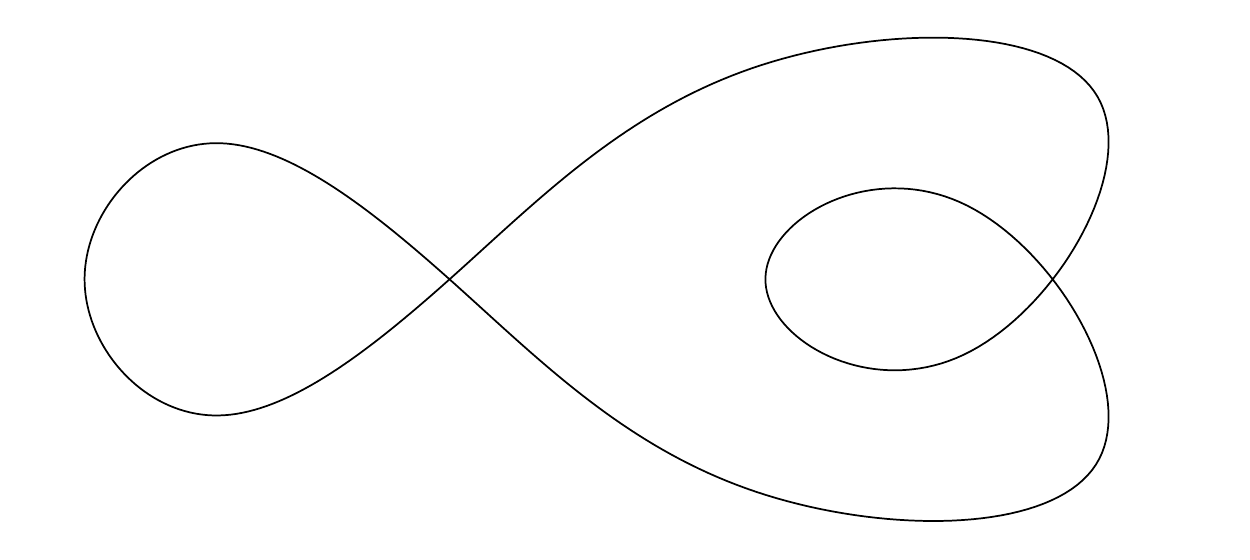}
\hspace{3mm}
\includegraphics[width=18mm]{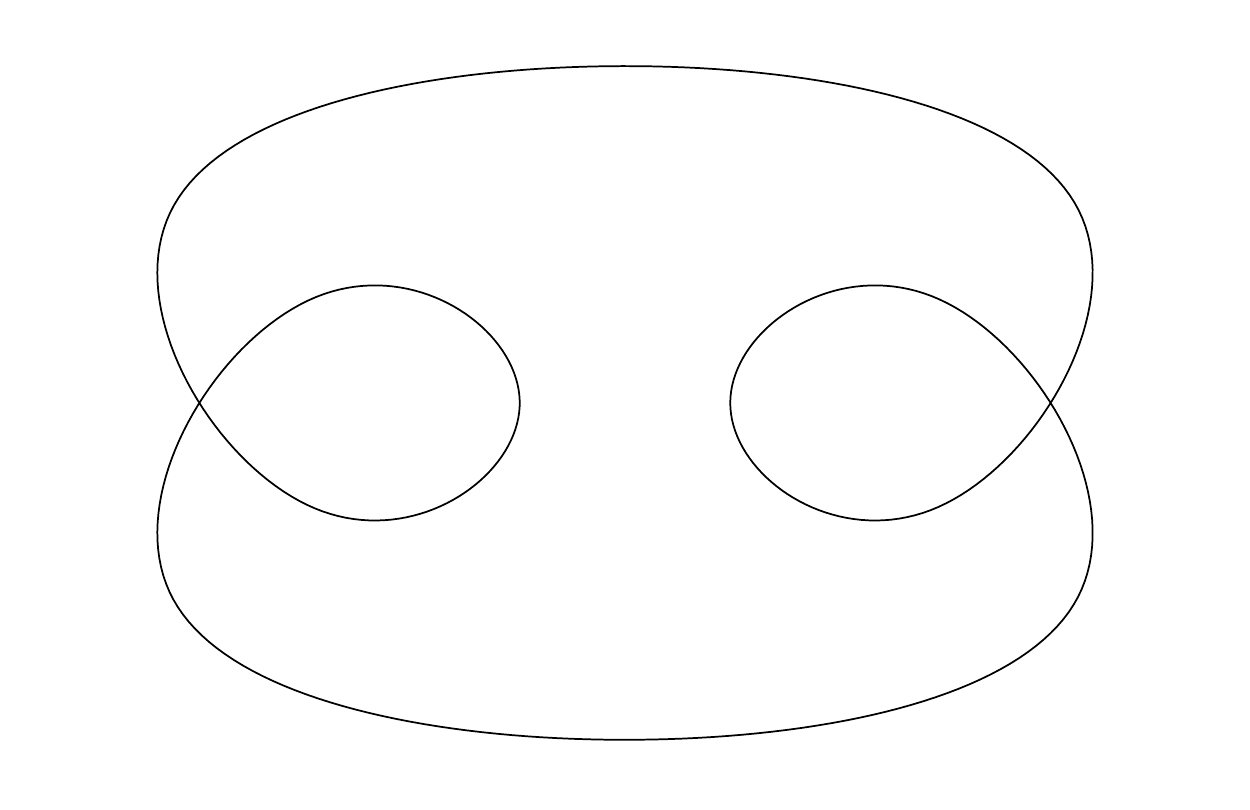}
\hspace{3mm}
\includegraphics[width=18mm]{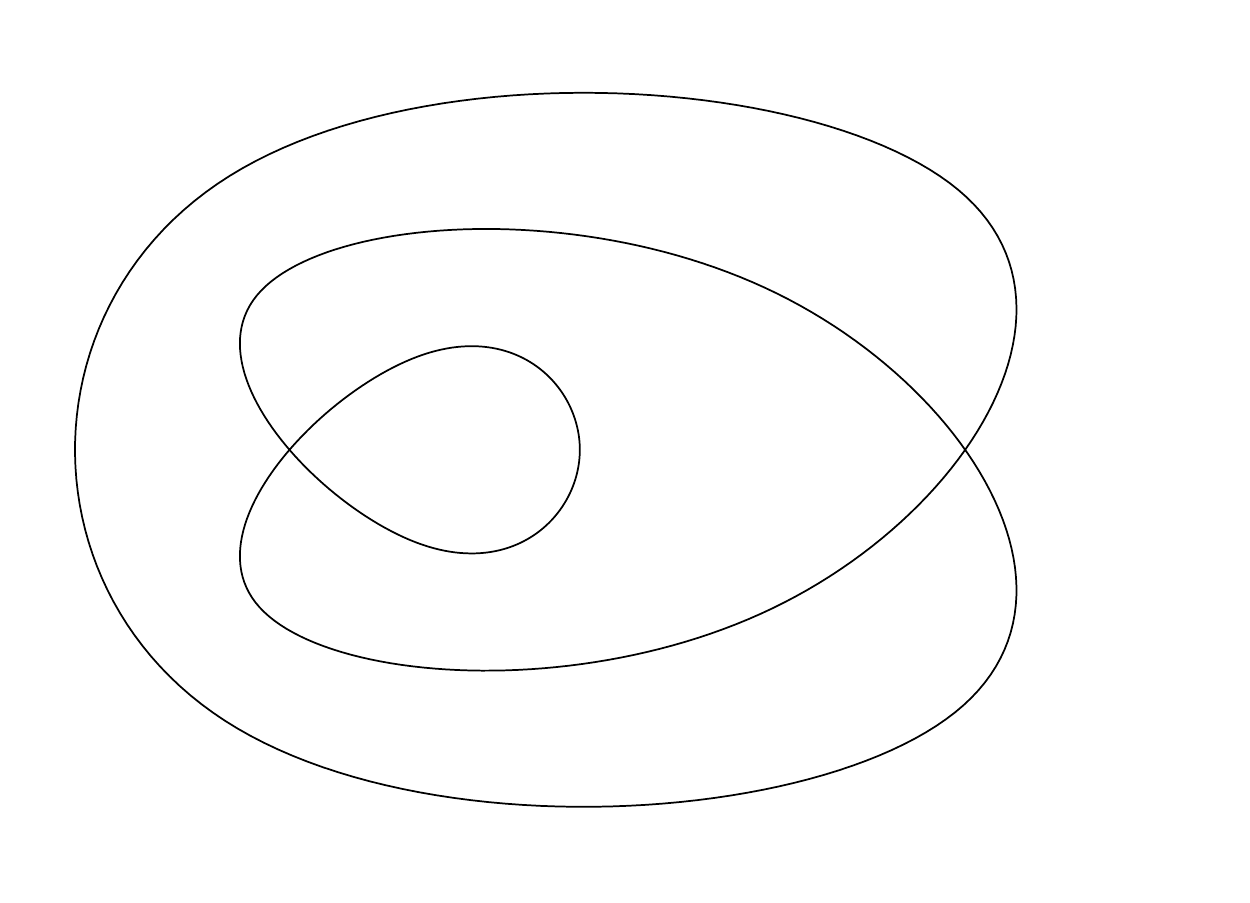}
\hspace{5mm}
\includegraphics[trim = 25mm 0mm 25mm 0mm, clip,width=13mm]{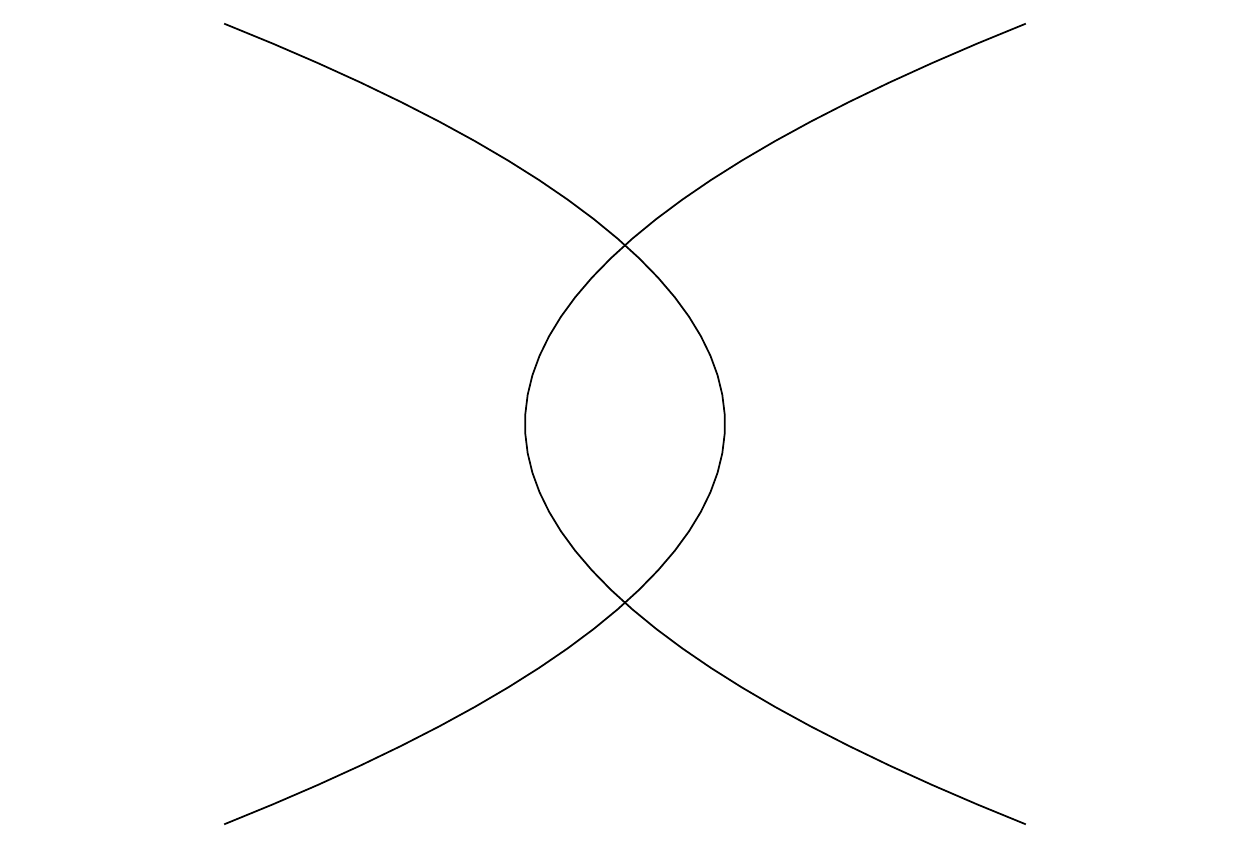}
}
\vspace{10mm}
\centerline{
\includegraphics[width=18mm]{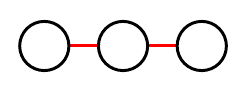}
\hspace{3mm}
\includegraphics[width=18mm]{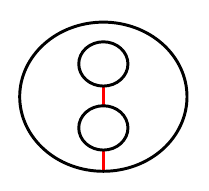}
\hspace{3mm}
\includegraphics[width=18mm]{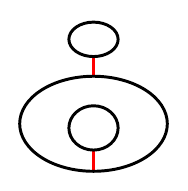}
\hspace{3mm}
\includegraphics[width=18mm]{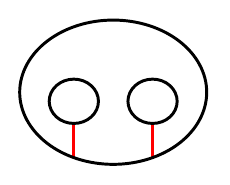}
\hspace{3mm}
\includegraphics[width=18mm]{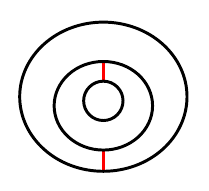}
\hspace{3mm}
\includegraphics[width=18mm]{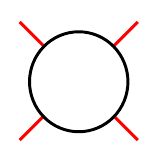}
}
\vspace{10mm}

\caption{Graphs, diagrams and corresponding immersions for
$n=2$, $d=0$
\label{n2classd0}}
\end{figure}

As the final example of complete classifications of all generic immersions
based on graphs $\Gamma(K)$ we consider the case $n=2$, $d=0$.
In this case $\Gamma(K)$ may have three vertices or a single vertex.
The classification is depicted on Figure \ref{n2classd0}
together with the smoothing diagrams themselves.
Note that in the case of a single vertex, there are a priori
two choices for the cyclic order at this vertex,
but only the one shown in Figure \ref{n2classd0}
satisfies to \ref{rem:deven}.
For the second choice, $\Sigma_A$ is a Klein bottle with holes
($V - E + B = 1 - 2 + 1 = 0$, cf.\ \ref{combprop}).

For higher $n$ the classification
of generic immersions gets rather complicated fast (cf. \cite{A}).
However, there is a class of relatively simple immersions
for arbitrary $n$ which we would like to distinguish.
\begin{defn}\label{def-arboreal}
A generic immersion $K\subset\rp^2$ of a circle is called
{\em arboreal} if the corresponding graph $\Gamma(K)$ is
a tree.
\end{defn}
\begin{prop}
If the degree $d$ is odd, then
the arboreal immersions are in 1-1 correspondence
with the ribbon rooted trees $T$
enhanced with the orientation of some
of its edges towards the root as well as the
data encoding the side change {\rm (}denoted by crosses
on the horizontal line in our pictures{\rm )}
for the edges adjacent to the root vertex of $T$.
We require all edges adjacent to the root to
be oriented.

If $d$ is even, then the arboreal immersions are
in 1-1 correspondence with the ribbon {\rm (}unrooted{\rm )} trees $T$
enhanced with orientations of some
of its edges in such a way that there exists at least one vertex
to which all the orientations point out.
\end{prop}
\begin{proof}
If $d$ is odd, we start by the pseudoline representing
the root and add the vanishing cycles adjacent to it
according to the side data.
To these cycles we attach the (even degree) immersions
corresponding to the connected components of $T$ minus
the union of its root and the open edges adjacent to the root.

If $d$ is even, the immersion $K$
may be deformed to $\R^2$ as $\Gamma(K)$ has no cycles.
Take the subtree $T'\subset T$
formed by the vertices on the headside of all oriented
edges. The subtree $T'$ also contains
all edges between such vertices (which must be non-oriented).
The subtree $T'$ is represented by a collection of non-nested
ovals in $\R^2$ and vanishing cycles between them.
For each oriented edge adjacent to $T'$ we take
a vanishing cycle adjacent to the corresponding oval from
its interior and proceed inductively.
\end{proof}

\ignore{
\begin{prop}
There is at most one non-orientable connected component of $\rp^2\setminus C$.

For an arboreal immersion $K\subset\rp^2$ of a circle
the complement $\rp^2\setminus K$ always has a non-orientable
connected component.
\end{prop}
\begin{proof}
Proposition follows by induction. Consider a new tree obtained
by removing a one-valent vertex from $\Gamma(K)$.
By Theorem \ref{m-graph} it corresponds to a generically immersed circle
$K'\subset\rp^2$ with a smaller number of nodes. The complement
$\rp^2\setminus K$ is homeomorphic to the disjoint union
of $\rp^2\setminus K'$ and a disk.
\end{proof}

\begin{defn}\label{infty-comp}
The non-orientable connected component of $\rp^2\setminus K$ (if it exists)
is called the {\em infinity component}.
\end{defn}

We may also distinguish a case (not necessarily arboreal)
when no orientations of edges appear (unless they are forced
by the root of the tree).
\begin{defn}
We say that an immersion of a circle is
{\em unnested} if all edges of $\Gamma(K)$ are unoriented,
unless they are adjacent to the root of $\Gamma(K)$
(if $d\neq 0$).
\end{defn}
In other words, $K$ is unnested if there are now nested
pairs of ovals in its smoothing diagrams $\Delta(K)$.

Unnested arboreal immersions provide an especially simple
class of generic immersions of circles into $\rp^1$.
}
\subsection{Classification of generic real rational curves of degree \texorpdfstring{$4$}{<= 4}}\label{sec:degree4}
Here we illustrate how enhanced graphs
and related smoothing diagrams can be used for classification
of
real nodal rational quartic curves in $\rp^2$.
In this degree the classification is already known
(see e.g. \cite{DMello}
for its recently found description in terms of chord diagrams).
Let us consider several ways in which this classification can be formulated.


Recall that for a nodal rational curve in $\rp^2$, we denote by
$h$, $e$, and $c$ the numbers of hyperbolic, elliptic, and
imaginary nodes of the curve (see Section \ref{sec-rational}).

\begin{thm}[see \cite{DMello}] \label{deg4thm}
There are $13$ isotopy types of nodal rational curves of degree $4$
in $\rp^2$. These isotopy types are
listed in Table \ref{quartic}.
{\rm (}The table contains each isotopy type on the right together
with its ``chord diagram'', explained in
Proposition \ref{13types} and Remark \ref{rem-chorddiagrams}, on the left.{\rm )}
\end{thm}

\begin{table}%
\centering
     $h=3$

    \input{files/gmpics/chorddiagram1.TpX}
    \input{files/gmpics/quartic1.TpX}
		\hspace{2cm}
    \input{files/gmpics/chorddiagram2.TpX}
    \input{files/gmpics/quartic2.TpX}

    \input{files/gmpics/chorddiagram3.TpX}
    \input{files/gmpics/quartic3.TpX}
		\hspace{2cm}
    \input{files/gmpics/chorddiagram4.TpX}
    \input{files/gmpics/quartic4.TpX}

    \input{files/gmpics/chorddiagram5.TpX}
    \input{files/gmpics/quartic5.TpX}

    $e>0, c=0$

    \input{files/gmpics/chorddiagram7.TpX}
    \input{files/gmpics/quartic7.TpX}
		\hspace{2cm}
    \input{files/gmpics/chorddiagram6.TpX}
    \input{files/gmpics/quartic6.TpX}

    \input{files/gmpics/chorddiagram8.TpX}
    \input{files/gmpics/quartic8.TpX}
		\hspace{2cm}
    \input{files/gmpics/chorddiagram9.TpX}
    \input{files/gmpics/quartic9.TpX}

    $c=2, h=1$

    \input{files/gmpics/chorddiagram10.TpX}
    \input{files/gmpics/quartic10.TpX}
		\hspace{2cm}
    \input{files/gmpics/chorddiagram11.TpX}
    \input{files/gmpics/quartic11.TpX}

    $c=2, e=1$

    \input{files/gmpics/chorddiagram12.TpX}
    \input{files/gmpics/quartic12.TpX}
		\hspace{2cm}
    \input{files/gmpics/chorddiagram13.TpX}
    \input{files/gmpics/quartic13.TpX}
	\caption{Nodal rational quartics} 
	\label{quartic}%
\end{table}

\begin{proof}
Assume first that $C$ does not have imaginary nodes
and that $K=\R C\setminus E$
(where $E$ is the set of elliptic nodes of $C$) is not an arboreal immersion.
Then, we have $n=3-e$ and $l=n-1=2-e$, so $e=0,1$
(as usual, $n$ is the number of edges of $\Gamma(K)$, while
$l$ is the number of its vertices; these two numbers have opposite parity
for an immersion of a circle).
%

If $e=1$, then the smoothing diagram $\Delta_K$ consists of a single oval
with two vanishing cycles. The corresponding immersion
is unique by our $n=2$ classification. If the elliptic node
sits outside the oval of $\Delta_K$ we get a contradiction
to the B\'ezout theorem by tracing a line through an elliptic
and one of the two hyperbolic nodes of $C$.

If $e=0$, then the graph $\Gamma(K)$ (which we
assumed to be non-arboreal) has two vertices
and three edges.
Suppose that there is no
vanishing cycle intersecting
the auxiliary pseudoline $J$.
Then, the corresponding two ovals of $\Delta_K$
must be nested.
Otherwise, the line $L$ connecting
any two of the three hyperbolic nodes
must intersect $\R C$ also somewhere else
by topological reasons.
We get a contradiction with the B\'ezout theorem as
each node already contributes two to the intersection
number of $L$ and $C$.

If there is
a vanishing cycle intersecting $J$, then
each such vanishing cycle
must correspond to a loop of $\Gamma(K)$.
Indeed, if the two ovals of $\Delta_K$ are nested,
then only the exterior oval may be adjacent to
a vanishing cycle intersecting $J$.
The unnested components
correspond to vertices of different color, so
a vanishing cycle intersecting $J$ cannot connect them
in a way coherent with the orientation.

If we have a single
loop at a vertex $v\in\Gamma(K)$,
then there are two edges $e_1,e_2$
connecting $v$ to a vertex $v'$.
These are the only edges adjacent to $v'$.
Therefore, $e_1$ and $e_2$ must be separated from each other by
the two endpoints of the loop edge, as otherwise $K$ would have
several components after normalization.
Once again this excludes the possibility that the two ovals
of $\Delta_K$ are unnested.
The unique nested configuration
is listed in Table \ref{quartic}.

If there are two loops at a vertex $v\in\Gamma(K)$, then
we have a single edge connecting $v$ to $v'$.
As in the $e=1$ case, the cyclic ordering of the loop edges at $v$
is unique by the $n=2$ classification.
When inserting the edge connecting $v$ to $v'$, there is again only one choice
due to symmetries.


We are left to consider the case when $C$ has
a pair of imaginary nodes.
There is exactly one remaining node. If it is hyperbolic, then
the corresponding vanishing cycle connects two ovals that can be either nested
or unnested. If the real node is elliptic, then $\R C\setminus E$ is an oval.
The elliptic node can be either inside or outside this oval.
\end{proof}

All thirteen topological types of $(\rp^2,\R C)$ described above
can be easily realized by quadratic (Cremona)
transformations of conics as specified in the following statement.

\begin{prop}\label{13types}
The 13 types of Theorem \ref{deg4thm} can be obtained
from conics by the quadratic
transformation
\begin{equation}\label{cremona}
[x_0:x_1:x_2]\mapsto [x_1x_2:x_0x_2:x_0x_1].
\end{equation}
For each type, we show a suitable arrangement of a conic and the coordinate
axes next to the curve in Table \ref{quartic}.
\end{prop}

The proof of this proposition is straightforward.

\begin{remark}[D'Mello, Viro, see \cite{DMello}] \label{rem-chorddiagrams}
Nodal rational quartics in $\rp^2$ correspond to {\it topological
chord diagrams}, that is, topological types of embeddings
of a disjoint union of zero-dimensional spheres
into the circle $S^1$ {\rm (}in the figures, each zero-dimensional sphere is represented
by the chord joining the images of the two points of the sphere{\rm )}.
The number of chords of a diagram is the
number of hyperbolic nodes of the curve.
\begin{itemize}
\item
The nine classes of nodal rational quartic curves without imaginary
nodes
{\rm (}cf. Table \ref{quartic}{\rm )}
correspond to nine possible
chord diagrams with no more than three chords.
\item
The four classes of nodal rational quartic curves with a pair
of imaginary nodes
{\rm (}cf. Table \ref{quartic}{\rm )}
correspond
to four possible chord diagrams with no more than one chord
enhanced with a single {\em imaginary chord} data.
The latter data
says whether the imaginary node corresponds
to the intersection of the same or different halves
of $\C \tilde C\setminus \R \tilde C$, where $\tilde C$ is the normalization
of $C$ {\rm (}corresponding to the values $\sigma(C) = 0$ and $\sigma(C) = 2$
in Proposition \ref{Hilbert16_nodal} below{\rm )}.
\end{itemize}

As it was noticed by Viro,
the correspondence with the chord diagrams is provided by the
real point set of the conic $Q$ obtained from $C$ by the quadratic transformation
centered in the nodes of $C$  {\rm (}the  quadratic transformation
of Proposition \ref{13types}{\rm )}
together with the parts of the real axes of $\rp^2$ in the interior
of the ellipse $\R Q$.
\end{remark}

The 13 types of Theorem \ref{deg4thm} can be decomposed
into groups according to the topological type of the smoothing $\R C_\circ$
(see Definition \ref{smoothingdiagram} for $\R C$, plus perturbing each elliptic
node into an oval, cf.\ Proposition \ref{nodal_perturbation} below).

\begin{lemma}\label{lemtypeI}
The smoothing $\R C_\circ \subset\rp^2$ of the real point set of a nodal rational quartic
$C$ in $\rp^2$ is isotopic to the real point set of a smooth quartic of type
$<4>, <1<1>>$ or $<2>$. If $c=0$, only the first two cases
appear.
\end{lemma}

\begin{proof}
A nodal rational quartic is of type I in the sense of Definition \ref{typeI_nodal}.
By Proposition \ref{nodal_perturbation}, the ``oriented''
small perturbation $C_\circ$ is also of type I.
In degree $4$,  the genus $g(C_\circ)$ of  $C_\circ$ is equal to $3$, and hence the number of connected components
of $\R C_\circ$ is $l=2,4$. Thus, the first statement follows from the classification of
smooth quartics in example \ref{classd4}.
Note that, by the complex orientation formula \ref{Hilbert16_nodal} for nodal curves,
we have $\sigma(C) = 0$ for $<4>, <1<1>>$
and $\sigma(C) = 2$ for $<2>$.
 %
%
\end{proof}


The following statement can be checked case by case from Table \ref{quartic}.

\begin{prop}
Let $C$ be a nodal rational quartic in $\rp^2$, and let $\R C_\circ$ be the smoothing
of the real point set of $C$.
Furthermore, let $Q,L_1, L_2,L_3$ be the arrangement of conic and three lines obtained
by the quadratic transformation centered in the nodes of $C$
{\rm (}cf.\ Remark \ref{rem-chorddiagrams} and Proposition \ref{13types}{\rm )}.
Then $\R C_\circ$ is of type $<1<1>>$ if and only if the interior
$\text{Int}(\R Q)$
of the oval $\R Q$ contains at least one of the intersection points
of the three lines.
%
%
%
\end{prop}

We conclude this section with yet another reformulation
(with a help
of the $J_-$-invariant)
of the classification of Theorem \ref{deg4thm} in the case of
three hyperbolic nodes.

\begin{theorem}\label{deg4reform}
A generic immersion $K$ of a circle $S^1 \to \rp^2$
with $n(K)=3$ is realizable by
a real nodal rational curve of degree $4$
if and only if $J_-(K) = -3$
and the isotopy type of $K$ is different from the one
depicted below:

\begin{minipage}{\textwidth}
\centering
    \input{files/gmpics/quartic_nonreal.TpX}
\end{minipage}
\end{theorem}
\begin{proof}
By Definitions \ref{Jminus} and \ref{Orinv},
$J_-(K) = -3$ is equivalent to $\Or(K) = 4$.
If $K$ is realizable by
a real nodal rational curve of degree $4$,
then by complex orientation formula \ref{Hilbert16}
we have $\Or(K) = 4$.
(Also we can verify this directly from
the classification of Theorem \ref{deg4thm}.)
Conversely, let $K$ be a generic immersion with three nodes and
such that
$\Or(K)=4$. The number of connected components of the smoothing
$K_\circ$ is $l(K)=4$ or $l(K)=2$. It easy to check that in order to
have $\Or(K)=4$, the smoothing $K_\circ$ must be of type $<4>$, $<1<1>>$ (a negative
injective pair) or $<1\sqcup1<2>>$ (with opposite orientations on
the two interior ovals). For $<4>$ there are two possible smoothing diagrams,
for $<1<1>>$ there are three, and all of them appear in Table \ref{quartic}.
For $<1\sqcup1<2>>$, the orientations only allow
for a single smoothing diagram.
The corresponding immersed circle is depicted
above.
\end{proof}

%
%
%
%
%
%


\section{Several restrictions on the topology of real algebraic curves}\label{restrictions}
\subsection{Complex orientations}\label{complex_orientations}

Recall that a \emph{real curve} is a pair $(\Sigma, \varphi)$, where
$\Sigma$ is a Riemann surface and
$\varphi: \Sigma \to \Sigma$ is an antiholomorphic involution.
The curve is irreducible if $\Sigma$ is connected.
The fixed point set of $\varphi$ is called the \emph{real part} of $\Sigma$ and is denoted by $\R\Sigma$.
An example of real curves is provided by nonsingular algebraic curves $\R C\subset\rp^2$:
the restriction of the involution of complex conjugation
$\conj: \cp^2 \to \cp^2$ to
the complex point set $\C C$ of such a curve $C$ is an antiholomorphic involution
on the Riemann surface $\C C$.

If $(\Sigma, \varphi)$ is an irreducible real curve, then either $\Sigma \setminus \R\Sigma$
consists of two connected components exchanged by $\varphi$, or $\Sigma \setminus \R\Sigma$ is connected.
In the first case, $\Sigma$ is said to be of \emph{type I} (or \emph{separating});
in the latter case, $\Sigma$ is said to be of \emph{type II}.
If $\Sigma$ is of type I, the two halves of $\Sigma \setminus \R\Sigma$ induce two opposite orientations
of $\R \Sigma$. These orientations are called \emph{complex orientations}.

The {\em complex scheme} of a nonsingular algebraic curve $C$ in $\rp^2$
is the topological type of the pair $(\rp^2,\R C)$ enhanced with the information of the type (I or II)
of the curve and, in the case of type I,
with one of two complex orientations
of $\R C$.
Namely, we say that two nonsingular algebraic curves $C$ and $C'$ in $\rp^2$
have the same complex scheme if they have the same type and there exists a homeomorphism
of pairs $(\rp^2,\R C)$ and $(\rp^2,\R C')$ that is consistent with complex orientations in the
case of type I. 

A powerful restriction on complex orientations of a nonsingular curve of type I in $\rp^2$
is provided by Rokhlin's complex orientation formula.
We present here this formula in the form proposed by O.\ Viro.
Choose a complex orientation of $\R C$.
Then the invariant $\Or(\R C)=\int\limits_{\rp^2}\ind^2_{\R C}d\chi$ from Definition \ref{Orinv} is well-defined.
Note that if we choose the opposite complex orientation of $\R C$
then $\Or(\R C)$ stays the same.
Thus, $\Or(\R C)$ is an invariant of the complex scheme of $\R C$.

\begin{theorem}[Rokhlin's complex orientation formula, cf. \cite{Viro-chi}]
\label{Hilbert16}
Let $C$ be a nonsingular curve of degree $d$ and type I in $\rp^2$. Then,
$$
\Or(\R C) = \frac{d^2}{4}.
$$
\end{theorem}

For curves of odd degree, the statement of Theorem \ref{Hilbert16} can be reformulated
in the following way. Let $C$ be a nonsingular curve of odd degree and type I in $\rp^2$.
Denote by $J$ the pseudoline of $\R C$, and equip $J$ with an orientation.
This determines one of the two complex orientations of $\R C$.
Let $O$ be an oval of $\R C$. Denote by $[O]$ and $[J]$ the classes
in $H_1(\rp^2 \setminus Int(O); \Z)$
(where $Int(O)$ is the interior of $O$) which are realized by $O$ and $J$, respectively.
One has $2[J] = \pm[O]$.
Recall that the oval $O$ is 
positive 
(respectively, negative) 
if $2[J] = -[O]$ (respectively, $2[J] = [O]$).
Notice that positivity or negativity of an oval does not depend
on the choice of a complex orientation of $\R C$.
A pair of ovals of $\R C$ is called \emph{injective} if one of these ovals is contained
in the interior of the other one.
An injective pair of ovals is called \emph{positive} 
if some orientation of the annulus bounded by these ovals
induces a complex orientation of the ovals; otherwise, the injective pair is called \emph{negative}. 
For a nonsingular curve $C$ of degree $2k + 1$ and type I in $\rp^2$,
Rokhlin's complex orientation formula is equivalent to the equality
$$
2(\Pi_+ - \Pi_-) + \Lambda_+ - \Lambda_- = l - 1 - k(k + 1),
$$
where $\Pi_+$ (respectively, $\Pi_-$) is the number of positive (respectively, negative)
injective pairs, $\Lambda_+$ (respectively, $\Lambda_-$) is the number of positive (respectively, negative)
ovals, and $l$ is the total number of connected components of $\R C$.

\begin{definition}\label{typeI_nodal}
An irreducible nodal algebraic curve in $\rp^2$ 
is said to be of type I, if its normalization is of type I. 
\end{definition} 

Theorem \ref{Hilbert16} can be generalized to the case of
nodal real curves,
including those
with imaginary nodes.
Consider a nodal curve $C$ in $\rp^2$ such that all the nodes of $C$
are imaginary and $C$ is of type I.
Denote by $\widehat C$ the normalization of $C$,
and denote by $\widehat C_\pm$
the two connected components of $\C \widehat C \setminus \R {\widehat C}$.
\ignore{GM: rephrased on Sep 12 for consistency with Section 5}
Denote by
\begin{equation}\label{def-sigma}
\sigma(C)=\#(C_-\cap C_+)
\end{equation}
the number of nodes
of $C$ resulted as the intersection of the images
$C_\pm\subset\C C$
of $\widehat C_\pm$ under the restriction
$\widehat C_\pm\to \C C$
of the normalization map $\C \widehat C\to\C C$.

The following statement is a slight generalization of Rokhlin's complex orientation formula,
cf. \cite{Rokhlin, Viro-Progress,Viro-chi,Viro-immersions}.
The proof literally coincides with Rokhlin's proof of the complex orientation formula (see \cite{Rokhlin}).

\begin{theorem}[Rokhlin's complex orientation formula for nodal curves]
\label{Hilbert16_nodal}
Let $C$ be a nodal cur\-ve in $\rp^2$ such that all the nodes of $C$
are imaginary and $C$ is of type I.
Then,
$$
	\Or(\R C) = \frac{d^2}{4} - \sigma(C).
$$
\end{theorem}

Consider the pencil of real lines passing through the intersection point of two distinct real lines $L_0$ and $L_1$
in $\rp^2$.
This pencil is divided by $L_0$ and $L_1$ into two segments of the form $\{L_t\}$, $t \in [0, 1]$, where $L_t$
is defined by the linear form
$$
(1 - t)\ell_0 + t\ell_1  = 0
$$
under certain choice of linear forms $\ell_0$ and $\ell_1$ defining $L_0$ and $L_1$, respectively.
A point of tangency of two oriented curves is said to be \emph{positive} if the orientations of the curves
define the same orientation of the common tangent line at the point, and \emph{negative} otherwise.

\begin{theorem}[Fiedler's alternation of orientations, cf. \cite{Fie82}]\label{alternation}
Let $C$ be a nonsingular curve of type I in $\rp^2$.
Let $L_0$ and $L_1$ be real lines tangent to $\R C$ at points $p_0$ and $p_1$, respectively,
which are not points of inflection of $C$. Let $\{L_t\}$, $t \in [0,1]$, be a segment of the line pencil, connecting $L_0$ with $L_1$.
Orient the lines $\R L_0$ and $\R L_1$
coherently in $\{L_t\}$.
If there exists a path $f: [0, 1] \to \C C$ connecting the points $p_0$ and $p_1$, such that for any $t \in (0, 1)$,
the point $f(t)$ belongs to $\C C \setminus \R C$ and is a point of transversal intersection of $\C L_t$ with $\C C$,
then the points $p_0$ and $p_1$ are either both positive or both negative points of tangency of $\R C$
with $\R L_0$ and $\R L_1$, respectively.
\end{theorem}

\subsection{Small perturbations}\label{perturbations}

Let $p$ be a real node of a nodal algebraic curve $C$ in $\rp^2$.
For a small disk $D(p)$ centered at $p$, the intersection $\R C \cap D(p)$ consists either
of two intersecting arcs (in the case of hyperbolic node), or of the point $p$
(in the case of elliptic node).
A topological type of smoothing of $p$ is given by a topological type of a pair $(D(p), S)$
for an appropriate subset $S \subset D(p)$.
If $p$ is hyperbolic, the subset $S$ is formed by two non-intersecting arcs whose extremal points coincide
with the four points of $\R C$ on the boundary of $D(p)$;
there are two such topological types of smoothing of $p$.
If $p$ is elliptic, there are also two topological types of smoothing of $p$:
for one of them, $S$ is empty, for the other one, $S$ is an oval entirely contained in $D(p)$.

The following statement is known in topology of real algebraic curves under the name of \emph{classical small perturbation}.

\begin{theorem}[Brusotti theorem, cf. \cite{Brusotti}]\label{Brusotti}
Let $C$ be a nodal curve {\rm (}not necessarily irreducible{\rm )} of degree $d$ in $\rp^2$. 
Let $U$ be a regular neighborhood of $\C C$ in $\cp^2$, represented as the union of a neighborhood $U_0$
of the set of singular points of $C$
and a tubular neighborhood $U_1$ of the submanifold $\C C \setminus U_0$
in  $\cp^2 \setminus U_0$. Assume that $U_0 \cap \rp^2 = \cup_p D(p)$,
where the union is taken over all real nodes of $C$ and $D(p)$ is a small disk centered at $p$.
For each real node $p$ of $C$, choose either to keep $p$, or to smooth it; in the second case, choose
one of the two possible topological types of smoothing of $p$ in $D(p)$.
For each pair of imaginary conjugate points $C$, choose either to keep them, or to smooth both of them.
Then, for any neighborhood ${\cal U}_C$ of $C$
in the space $\R {\cal C}_d$ of all curves of degree $d$ in $\rp^2$,
there exists a nodal curve $\widetilde C \in {\cal U}_C$ of degree $d$ in $\rp^2$ such that
\begin{enumerate}
\item $\C {\widetilde C} \subset U$;
\item for each connected component $u$ of $U_0$, the intersection $\C {\widetilde C} \cap u$
is embedded in $u$ according to the choice made for the corresponding nodal point of $C$; 
\item $\C {\widetilde C} \setminus U_0$ is a section of the tubular fibration $U_1 \to (\C C \setminus U_0)$.
\end{enumerate}
\end{theorem}

A curve $\widetilde C$ as in Theorem \ref{Brusotti} is called a \emph{small perturbation} of $C$.

The following two statements concern the relation between the type (I or II) of a nodal curve in $\rp^2$
and the type of small perturbations of this curve.
They are proved by local considerations
at neighborhoods of the nodes.

\begin{figure}[ht]
\centerline{
\includegraphics[height=35mm]{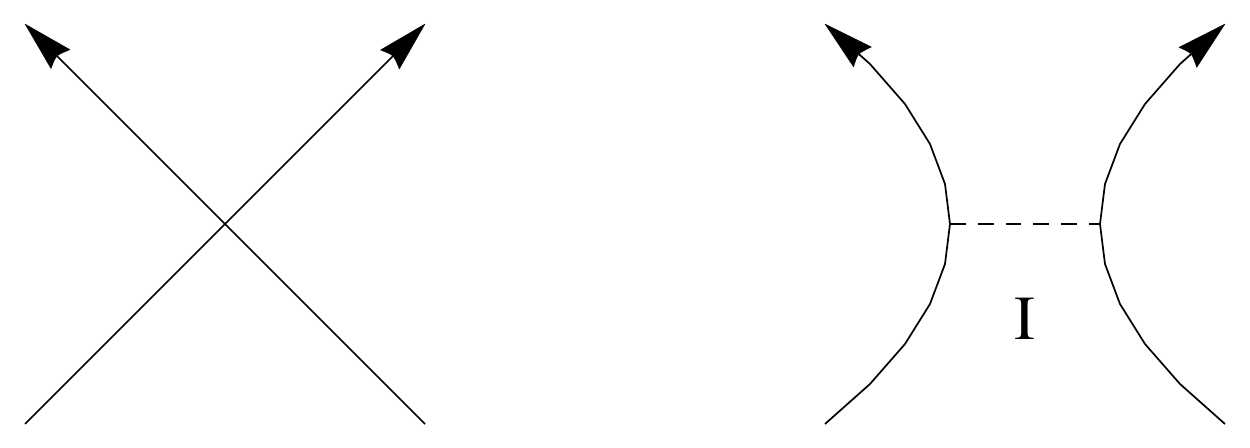}
}
\caption{Smoothing according to a complex orientation}
\label{or-smoothing1}
\end{figure}

\begin{proposition}[cf. \cite{Rokhlin, Mar80}]\label{nodal_perturbation}
Let $C$ be an irreducible nodal curve of type I in $\rp^2$.
Let ${\widetilde C}$ be a
curve obtained by a small perturbation of $C$ such that
each imaginary node of $C$ is kept,
each elliptic node of $C$ is turned to an oval of ${\widetilde C}$, and
each hyperbolic node of $C$ is perturbed according to the complex orientations
(see Figure \ref{or-smoothing1}).
Then, ${\widetilde C}$ is of type I. 
\end{proposition} 

\ignore{GM: Do we want to keep some of the real nodes
resulting from the intersection points?
\endgraf It: I think, no}
\begin{proposition}[cf. \cite{Rokhlin, Mar80}]\label{reducible_perturbation}
Let $C_1$, $\ldots$, $C_n$ be nonsingular curves of degrees $d_1$, $\ldots$, $d_n$ in $\rp^2$
such that no three of them pass through the same point, and $C_i$ intersects transversally $C_j$ in
$d_id_j$ 
points for any $1 \leq i < j \leq n$.
Let ${\widetilde C}$ be a
curve obtained by a small perturbation of the union $C_1 \cup \ldots \cup C_n$ in such a way that
all the imaginary intersection points of the curves $C_1$, $\ldots$, $C_n$ are kept,
and all the real intersection points of these curves are smoothed.
Assume that ${\widetilde C}$ is irreducible.
Then, ${\widetilde C}$ is of type I if and only if all the curves $C_1$, $\ldots$, $C_n$ are of type I
and there exists an orientation of $\R{\widetilde C}$
which agrees with some complex orientations of
$\R C_1$, $\ldots$, $\R C_n$
{\rm (}it means that the deformation turning
$C_1 \cup \ldots \cup C_n$ into ${\widetilde C}$
brings the chosen complex orientations of $C_i$ to the orientations
of the corresponding pieces of  $\R{\widetilde C}$
induced by a single orientation
of the whole $\R{\widetilde C}${\rm )}. 
In such case
this orientation
of $\R{\widetilde C}$ is
one of the complex orientations of ${\widetilde C}$.
\end{proposition}

\subsection{Rigid isotopies of nonsingular curves of degree \texorpdfstring{$5$}{5}}\label{rigid_isotopies}

Any curve of degree $5$ in $\rp^2$ is defined by a homogeneous real polynomial
in three variables of degree $5$. The multiplication of this polynomial by a non-zero real constant
gives rise to a polynomial defining the same curve.
Thus, the space $\R {\mathcal C}_5$ of all curves of degree $5$ in $\rp^2$ can be identified with
the real projective space of dimension $20$. The \emph{discriminant} $\Delta \subset \R {\mathcal C}_5$
is formed by the points of $\R {\mathcal C}_5$ which correspond to singular curves.
Two nonsingular curves of degree $5$ in $\rp^2$ are \emph{rigidly isotopic}
if the corresponding points belong to the same connected component of $\R {\mathcal C}_5 \setminus \Delta$.
It turns out that the rigid isotopy type of a nonsingular curve $C$ of degree $5$ in $\rp^2$
is determined by the topological arrangement of components of $\R C$ and the type (I or II) of $C$.

\begin{theorem}[Kharlamov, see \cite{Kh}]\label{rigid_classification}
There are nine rigid isotopy types of nonsingular curves of degree $5$ in $\rp^2$:
$<J \sqcup 6>_I$, $<J \sqcup 5>_{II}$, $<J \sqcup 4>_I$, $<J \sqcup 4>_{II}$, $<J \sqcup 3>_{II}$,
$<J \sqcup 2>_{II}$, $<J \sqcup 1>_{II}$, $<J>_{II}$, and
$<J\sqcup 1<1>>_I$,
where the subscript I or II indicates the type of curves.
\end{theorem}

We will need a more detailed information concerning the position of ovals
of nonsingular curves of degree $5$ and type I.
Let $C$ be a nonsingular curve of degree $5$ in $\rp^2$.
Let $x_1$ and $x_2$ be two distinct points in $\rp^2 \setminus J$,
where $J$ is the pseudoline of $\R C$.
The line passing through $x_1$ and $x_2$ intersects $J$ in odd number of points
(if the line is not transversal to $J$, we count intersection points with multiplicities).
Thus, exactly one of the two segments with endpoints $x_1$ and $x_2$ intersects $J$
in even number of points; we denote this segment by $[x_1, x_2]_C$.
A subset $S \subset (\rp^2 \setminus J)$ is \emph{convex with respect to $C$}
(or just \emph{convex} if $C$ is understood),
if for any two distinct points
$x_1$ and $x_2$ belonging to $S$ the segment $[x_1, x_2]_C$ is contained in $S$.
If a subset $S' \subset (\rp^2 \setminus J)$ is contained in a convex subset of $\rp^2 \setminus J$,
then we can consider the \emph{convex hull} of $S'$, that is, the smallest convex set containing $S'$.

\begin{proposition}\label{triangle}
Let $C$ be a nonsingular curve of degree $5$ in $\rp^2$
such that $\R C$ has at least three ovals.
Let $x_1$, $x_2$, and $x_3$ be points in the interiors of three distinct ovals of $\R C$. Then,
the union of the segments
$[x_1, x_2]_C$, $[x_1, x_3]_C$, and $[x_2, x_3]_C$
bounds a disc in $\rp^2\setminus J$.
This disc is the convex hull of $x_1$, $x_2$, and $x_3$.
\end{proposition}

\begin{proof}
The line passing through $x_i$ and $x_j$ (where $1 \leq i < j \leq 3$) intersects $\R C$ transversally in $5$ points:
two points of the oval whose interior contains $x_i$, two points of the oval whose interior contains $x_j$,
and one point of $J$.
Thus, the segment $[x_i, x_j]_C$ does not intersect $J$.
Since the union of our three segments does not intersect $J$,
this union bounds a disc in $\rp^2\setminus J$.
This disc coincides with one of the four triangles defined by the straight lines
passing through $x_1$ and $x_2$, $x_1$ and $x_3$, $x_2$ and $x_3$.
Clearly, it is convex and is contained in any convex set containing $x_1$, $x_2$, and $x_3$.
\end{proof}

The disc of Proposition \ref{triangle} is called the \emph{triangle} with vertices $x_1$, $x_2$, and $x_3$.
Let $O_1$, $\ldots$, $O_n$, $n \geq 3$, be a collection of ovals of a nonsingular curve $C$
of degree $5$ in $\rp^2$, and let $x_1$, $\ldots$, $x_n$ be points
in the interior of $O_1$, $\ldots$, $O_n$, respectively.
We say that the ovals $O_1$, $\ldots$, $O_n$ are \emph{in convex position}
if for any choice of indices $1 \leq i < j <  k \leq n$, the triangle with vertices $x_i$, $x_j$, and $x_k$
does not contain in its interior any point $x_1$, $\ldots$, $x_n$.
(B\'ezout theorem implies that the notion of convex position depends only on ovals $O_1$, $\ldots$, $O_n$
and not on the choice of points $x_1$, $\ldots$, $x_n$ inside these ovals.)

\begin{proposition}\label{quatre_ovales}
Let $C$ be a nonsingular curve of degree $5$ in $\rp^2$
such that $\R C$ has exactly four ovals.
Then, these ovals are in convex position if and only if $C$ is of type II.
\end{proposition}

\begin{figure}[ht]
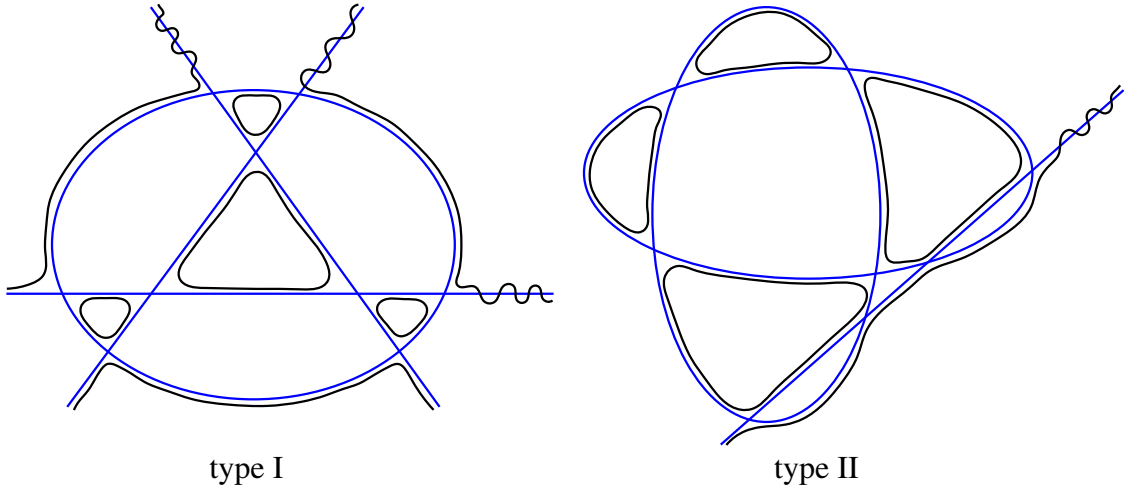

\centering
\input{files/pics/typeIdeformation.TpX}%
\input{files/pics/typeIIdeformation.TpX}
\caption{The construction of type I and II curves with 4 ovals from two ellipses and a line}
\label{two_types}
\end{figure}

\begin{proof}
The convexity of the position of four ovals of $C$
is invariant under rigid isotopies
by B\'ezout theorem.
Thus, Theorem \ref{rigid_classification} implies that, to prove the statement of the proposition,
it is enough to construct
\begin{itemize}
\item a nonsingular curve $C_1$ of degree $5$ and type I
in $\rp^2$ such that $\R C_1$ has exactly four ovals and these ovals are not in convex position, and
\item a nonsingular curve $C_2$ of degree $5$ and type II
in $\rp^2$ such that $\R C_2$ has exactly four ovals and these ovals are in convex position.
\end{itemize}
A construction of such curves is presented on Figure \ref{two_types}.
The fact that the first curve is of type I and the second curve is of type II follows
from Proposition \ref{reducible_perturbation}.
\end{proof}

\begin{proposition}\label{five_ovales}
Let $C$ be a nonsingular curve of degree $5$ in $\rp^2$
such that $\R C$ has at least five ovals.
Then, the ovals of $C$ are in convex position.
\end{proposition}

\begin{proof}
Let $x_1$, $\ldots$, $x_5$ be points in the interiors of five distinct ovals of $\R C$.
There exists a unique conic $A$ which passes through the five points $x_1$, $\ldots$, $x_5$.
B\'ezout theorem implies that this conic does not intersect $J$. The real part of $A$ is an oval,
and the disc bounded by it is convex with respect to $C$. Thus, this disc contains
the triangle with vertices $x_i$, $x_j$, and $x_k$ for any $1 \leq i < j < k \leq 5$. 
In particular, the conic $A$
does not have points inside the triangle.
\end{proof}

\ignore{GM: Was an M-curve defined?
\endgraf It: Yes}
Let $C$ be a nonsingular $M$-curve of degree $5$ in $\rp^2$,
and let $O_1$, $\ldots$, $O_6$ be the ovals of $\R C$.
Pick a point $x_i$ inside each oval $O_i$, $i = 1$, $\ldots$, $6$.
Points $x_i$ and $x_j$ are \emph{neighbors viewed from $x_k$}
(where $O_i$ and $O_j$ are two distinct ovals of $\R C$, and $O_k$ is another oval
of $\R C$) if
\begin{itemize} \label{def:neighbors}
\item
one of the segments of the line pencil ${\cal L}^k_{i, j}$ connecting the lines $x_kx_i$ and $x_kx_j$
does not contain any line which intersects an oval different from $O_i$, $O_j$, $O_k$;
denote this segment by ${\cal S}^k_{i, j}$;
\item there is a path $\sigma^k_{i, j} \subset \rp^2$ connecting $x_i$ and $x_j$ such that
any point of intersection of $\sigma^k_{i, j}$ with $\R C$ belongs either
to $O_i$ or to $O_j$, each line of ${\cal S}^k_{i, j}$ intersects $\sigma^k_{i, j}$ in one point,
and each line of ${\cal L}^k_{i, j } \setminus {\cal S}^k_{i, j}$ does not intersect $\sigma^k_{i, j}$.
\end{itemize}
B\'ezout theorem implies that, if $x_i$ and $x_j$ are neighbors viewed from $x_k$,
then for any choice of points $x'_i$, $x'_j$, and $x'_k$ inside the ovals $O_i$, $O_j$, and $O_k$, respectively,
the points $x'_i$ and $x'_j$ are neighbors viewed from $x'_k$.
In this case, we say that the ovals $O_i$ and $O_j$ are \emph{neighbors viewed from $O_k$}.

\begin{lemma}\label{conic_order}
Let $C$ be a nonsingular $M$-curve of degree $5$ in $\rp^2$,
and let $O_1$, $\ldots$, $O_6$ be the six ovals of $\R C$.
Assume that $O_i$ and $O_j$ are neighbors viewed from $O_k$.
Let $O_n$ and $O_m$ be two ovals different from $O_i$, $O_j$, and $O_k$.
Then, the real part $\R A$ of the conic $A$ which passes through the points $x_i$, $x_j$, $x_k$, $x_m$, and $x_n$
contains an arc which have endpoints $x_i$, $x_j$
and does not contain any of the points $x_k$, $x_m$, $x_n$.
\end{lemma} 

\begin{proof}
The conic $A$ does not intersect $J$, and
the points of $\R A$ are in a natural bijection with the lines of the pencil ${\cal L}^k$
centered at $x_k$. Thus, $\R A$
contains an arc $a$ with endpoints $x_i$, $x_j$ and such that
$a$ does not contain any of the points $x_m$, $x_n$.
Assume that $a$ contains the point $x_k$, and denote by $c_{i, j}$ the chord
connecting $x_i$ and $x_j$ in the interior of $\R A$.
There exists a path $\sigma^k_{i, j} \subset \rp^2$ connecting $x_i$ and $x_j$
and certifying that $O_i$ and $O_j$ are neighbors viewed from $O_k$.
The union of $\sigma$ and $c_{i, j}$ is a cycle that intersects once each line of ${\cal L}^k$.
Thus, this cycle is not homologous to $0 \in H_1(\rp^2; \Z)$. Hence, the cycle constructed
intersects $J$, which gives a contradiction.
\end{proof}

A \emph{reversible linear order} (respectively, \emph{reversible cyclic order}) on some collection of ovals
is a pair of opposite linear (respectively, cyclic) orders on this collection.

\begin{proposition}\label{convexity_maximal}
Let $C$ be a nonsingular $M$-curve of degree $5$ in $\rp^2$,
and let $O_1$, $\ldots$, $O_6$ be the six ovals of $\R C$.
\begin{enumerate}
\item For any oval $O_k$ of $\R C$, there exists a reversible linear order on the other five ovals of $\R C$
such that two ovals neighboring with respect to this reversible linear order are necessarily neighbors viewed from $O_k$.
\item If $O_i$ and $O_j$ are neighbors viewed from $O_k$, then
$O_i$ and $O_j$ are neighbors viewed from any oval $O_m$
such that $m$ is different from $i$ and $j$. 
\item If $O_i$ and $O_j$ are neighbors viewed from $O_k$, then one of the ovals $O_i$ and $O_j$
is positive, and the other one is negative.
\end{enumerate}
\end{proposition}

\begin{proof}
Pick a point $x_i$ inside each oval $O_i$, $i = 1$, $\ldots$, $6$.
The pseudoline $J$ of $\R C$ intersect each of the lines $x_kx_i$, $i = 1$, $\ldots$, $6$, $i \ne k$,
at exactly one point; denote this point by $y^k_i$.
To prove the statement (a), notice that the line pencil ${\cal L}^k$ centered at $x_k$
provides a reversible cyclic order on the five ovals different from $O_k$.
This pencil is formed by $5$ segments (with pairwise non-intersecting interiors), indexed by
pairs of ovals $(O_i, O_j)$ which are neighbors with respect to this order;
the segment $S^k_{i, j}$ indexed by $(O_i, O_j)$ connects the lines $x_kx_i$ and $x_kx_i$.
For such a segment $S^k_{i, j}$, the ovals $O_i$ and $O_j$
are neighbors viewed from $O_k$ if and only if
the orientations of the lines $x_kx_i$ and $x_kx_j$ provided by the triples
of points $(x_k, x_i, y^k_i)$ and $(x_k, x_j, y^k_j)$, respectively,
turn one to the one another through the segment $S^k_{i, j}$.
Our purpose is to prove that, among the five segments $S^k_{i, j}$,
there exists exactly one segment such that the corresponding ovals
are not neighbors viewed from $O_k$.
First, assume that there are two such segments $S^k_{i, j}$ and $S^k_{i', j'}$. 
Denote by $A$ the conic which passes through $x_k$, $x_i$, $x_j$, $x_{i'}$, $x_{j'}$
(in the case where the indices $i$, $j$, $i'$, $j'$ are not pairwise distinct,
we choose an oval $O_m$ different from $O_i$, $O_j$, $O_{i'}$, $O_{j'}$,
and suppose that $A$ passes through $x_m$).
The five marked points divide
the real part $\R A$ of $A$ into five arcs.
Either $x_i$, $x_j$ or $x_{i'}$, $x_{j'}$ are endpoints of such an arc
which contradicts the fact that $O_i$ and $O_j$, as well as $O_{i'}$ and $O_{j'}$,
are not neighbors viewed from $O_k$.
Furthermore, if any pair of ovals which are neighbors with respect
to the reversible cyclic order provided by ${\cal L}^k$
are neighbors viewed from $O_k$, then we get five paths whose union is a cycle $c$
intersecting once each line of ${\cal L}^k$; the cycle $c$ is not homologous to $0 \in H_1(\rp^2; \Z)$,
thus $c$ intersects $J$.

To prove the statement (b), consider the segment $S$ of the line pencil
connecting the lines $x_mx_i$ and $x_mx_j$ such that no line of $S$ intersects $O_k$,
and assume that some line of $S$ intersects an oval $O_n$, where $n$ is different from $i$, $j$, $k$, and $m$.
Trace a conic $B$ through the points $x_i$, $x_j$, $x_k$, $x_m$, and $x_n$.
Since $O_i$ and $O_j$ are neighbors viewed from $O_k$,
the real part $\R B$ of $B$ contains an arc $a$ such that its endpoints are $x_i$, $x_j$,
and the interior of $a$ does not contain any of the points $x_m$, $x_n$, $x_k$
(see Lemma \ref{conic_order}),
which gives a contradiction.

To prove the statement (c), consider the segment
${\cal S}^k_{i, j}$
of the line pencil ${\cal L}^k_{i, j}$ connecting the lines $x_kx_i$ and $x_kx_j$
such that any line of ${\cal S}^k_{i, j}$ does not intersect any oval different from $O_i$, $O_j$, $O_k$. 
Orient all the real parts of the lines of the segment in such a way that the orientations turn to one another under the isotopy
given by these real parts, and choose a complex orientation of $\R C$.
Let $\widetilde{\cal S}^k_{i, j}$ be a subsegment of ${\cal S}^k_{i, j}$
such that the endpoints of $\widetilde{\cal S}^k_{i, j}$ correspond to lines
tangent, respectively, to $O_i$ and $O_j$, and the interior points of $\widetilde{\cal S}^k_{i, j}$
correspond to lines which do not intersect any oval except $O_k$.
The points $\widetilde{\cal S}^k_{i, j}$ which correspond to lines tangent to $\R C$
divide $\widetilde{\cal S}^k_{i, j}$ into segments
$S_1$, $\ldots$, $S_r$ (each of these lines is tangent to $\R C$ at exactly one point).
The endpoints of each of the segments $S_1$, $\ldots$, $S_r$ correspond to two tangency points $p$, $q$ of $\R C$ with lines
of $\widetilde{\cal S}^k_{i, j}$; these two points
are either both positive or both negative.
Indeed, if the interior points of a segment $S_t$ under consideration
correspond to lines intersecting $\R C$ at $3$ points,
then the statement follows from Theorem \ref{alternation}.
Assume that the interior points of a segment $S_t$ under consideration
correspond to lines intersecting $\R C$ at $5$ points.
Then, for each line corresponding to an interior point of $S_t$,
either $3$ intersection points with $\R C$ belong to $J$ and the other two intersection
points belong to $O_k$, or one intersection point belongs to $J$ and the other four intersection
points belong to $O_k$. In the first case, B\'ezout theorem implies that $p$ and $q$ are connected by an arc $\alpha$ of $J$
such that $\alpha$ has exactly one common point with any line of $\widetilde{\cal S}^k_{i, j}$;
furthermore, $p$ and $q$ are either both positive or both negative.
In the second case, since $x_k$ is inside of $O_k$, the points $p$ and $q$ are connected by an arc $\beta$ of $O_k$
such that $\beta$ has exactly one common point with any line of $\widetilde{\cal S}^k_{i, j}$;
furthermore, $p$ and $q$ are either both positive or both negative.
\end{proof}

\begin{corollary}\label{cyclic_order}
Let $C$ be a nonsingular $M$-curve of degree $5$ in $\rp^2$.
Then, there exists a reversible cyclic order
of the six ovals of $\R C$ such that
any two ovals which are neighbors with respect to this order are
neighbors viewed from any other oval of $\R C$.
This reversible cyclic order is invariant up to rigid isotopy of the curve;
positive and negative ovals alternate with respect to this order.
\qed
\end{corollary}

\begin{figure}[ht]
\centering
\input{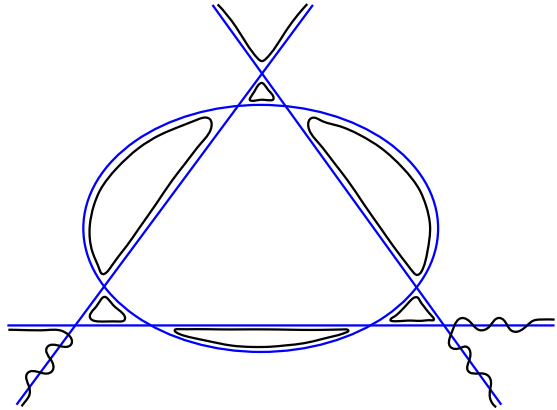}%
\caption{The construction of a $M$-curve from three lines and an ellipse}
\label{maximal_construction}
\end{figure}

\begin{remark}
Corollary \ref{cyclic_order} can as well be deduced from the rigid isotopy classification
of nonsingular curves of degree $5$ in $\rp^2$ (Theorem \ref{rigid_classification}),
the fact that the existence of a reversible cyclic order required in the corollary
is invariant under rigid isotopies, and a particular construction
of a maximal curve of degree $5$ in $\rp^2$, for example, a small
perturbation of three lines and an ellipse shown on Figure \ref{maximal_construction}.
\end{remark}

\section{Tropical constructions}\label{sec:tropical}
This section is devoted to the combinatorial patchworking
and its tropical interpretation.
The patchworking technique was invented by O. Viro at the end of the 1970's.
This technique provides a powerful tool to construct real plane algebraic curves
(and, more generally, real algebraic hypersurfaces in toric varieties).
We discuss here only certain particular cases of the general patchworking theorem.

\subsection{Nodal tropical curves}\label{nodal}
We give here the definitions required for the combinatorial patchworking construction presented below.
An introduction to tropical geometry and
a detailed information on tropical curves can be found, for example, in \cite{BIMS}.
A {\it
tropical curve} in $\R^2$ is a finite weighted rectilinear graph $\Gamma$ in $\R^2$
(some of the edges of $\Gamma$ are not bounded) such that
\begin{itemize}
\item each edge $e$ of $\Gamma$ has a rational slope
and prescribed a positive integer weight
$w(e)$,
\item at each vertex $v$ of $\Gamma$
the following {\it balancing condition} is satisfied:
\begin{equation}
\label{bcondition}
\sum_{e \supset v}w(e)u_v(e) = 0,
\end{equation}
where the sum is taken over all edges $e$ adjacent to $v$,
the vector $u_v(e)$
is the primitive integer vector
({\it i.e.}, vector with integer relatively prime coordinates)
in the direction of $e$ and pointing outward of $v$.
\end{itemize}

\ignore{GM: The definition of the Newton polygon
modified}
Consider the collection $\mathcal C$
of integer vectors $w(e)u_v(e)$
where $e$ runs over all non-bounded edges of $\Gamma$
and $v$ is the vertex adjacent to $e$.
By \eqref{bcondition} the sum of all vectors in the collection
is zero.
Thus, there exists a convex polygon $\Delta$ 
with integer vertices in $\R^2$ dual to
this collection.
This means that each vector $w(e)u_v(e)\in{\mathcal C}$
is an outward normal to a side $E\subset\Delta$ and
$$\#(E\cap\Z^2)-1=\sum w(e),$$
where the sum is taken over all $w(e)u_v(e)\in{\mathcal C}$
that are outward normal vectors to $E$.
The polygon $\Delta$ is called {\it Newton polygon} of $\Gamma$.
It is defined up to translation.
The quantity $\#(E\cap\Z^2)-1$ is called {\em the integer
length of the interval $E$} (recall that the endpoints
of $E$ are from $\Z^2$).

If $\Delta$ can be chosen to coincide with the triangle with vertices $(0, 0)$, $(d, 0)$, $(0,d)$
for some positive integer $d$, then we say that $\Gamma$ is
{\it projective of degree $d$}.
\ignore{GM: changed to "projective" of degree $d$}
The latter means that each vector $u_v(e)$, where $e$ is a non-bounded edge of $\Gamma$,
is either $(-1, 0)$, or $(0, -1)$, or $(1, 1)$, and the number of such vectors (counted with weights $w(e)$)
in each direction is equal to $d$.

\ignore{It: definition added}
A tropical curve $\Gamma$ in $\R^2$ is said to be {\it irreducible}, if it cannot be presented as a union
of two tropical curves 
different from $\Gamma$.

The space ${\mathcal T}(\Delta)$ of tropical curves
with a given Newton polygon $\Delta \subset \R^2$
is equipped with a natural topology
induced by the Hausdorff distance
$$
d(\Gamma_1, \Gamma_2) = \max\{\sup_{p \in \Gamma_1}\inf_{q \in \Gamma_2}\text{\rm dist}(p, q), \
\sup_{q \in \Gamma_2}\inf_{p \in \Gamma_1}\text{\rm dist(p, q)}\},
$$
where $\text{\rm dist}(p, q)$
is the Euclidean distance between points $p$ and $q$.
The condition that $\Gamma_1$ and $\Gamma_2$ have
the same Newton polygon ensures finiteness of
this distance.

A
tropical curve $\Gamma$ in $\R^2$ is said to be {\it nonsingular} if
\begin{itemize}
\item each edge of $\Gamma$ is of weight $1$,
\item each vertex $v$ of $\Gamma$ is $3$-valent and
the primitive integer vectors in the directions of three edges adjacent to $v$
generate (over $\Z$) the lattice $\Z^2 \subset \R^2$ of vectors with integer coordinates.
\end{itemize}

A tropical curve $\Gamma$ in $\R^2$ is said to be {\it nodal}
(or {\it simple}, cf. \cite{Mi05})
\ignore{GM: inserted an alias "simple" as this was
the term used in \cite{Mi05}}
if
\begin{itemize}
\item each vertex of $\Gamma$ is either $3$-valent or $4$-valent,
\item for each $4$-valent vertex of $\Gamma$, the union of four edges adjacent to this vertex
is contained in a union of two straight lines.
\end{itemize}

For our constructions, we use some particular degenerations of nonsingular tropical curves
to nodal ones. These degenerations contract certain edges of nonsingular tropical curves in $\R^2$,
as well as some "triangles". A {\it triangle} in a nonsingular tropical curve $\Gamma$ in $\R^2$
is a collection of three edges of $\Gamma$ which form a cycle such that no vertex of $\Gamma$
is inside this cycle.

Let 
$\Delta \subset \R^2$ be a convex polygon with integer vertices,
and let $\gamma:\to {\mathcal T}(\Delta)$ be a path
such that
\begin{itemize}
\item $\gamma(t)$ is a nonsingular tropical curve for any $t \in (0, 1]$;
\item $\gamma(0)$ is a nodal tropical curve;
\item the underlying graph of the tropical curve $\gamma(0)$ can be obtained
from the underlying graph of $\gamma(t)$ for any $t \in (0, 1]$
by contraction of a collection ${\mathcal C}_{\text{\rm edges}}$ of pairwise disjoint edges ({\it i.e.}, no two edges
of ${\mathcal C}_{\text{\rm edges}}$ have
common endpoint) and a collection ${\mathcal C}_{\text{\rm triangles}}$ of pairwise disjoint triangles
({\it i.e.}, no two edges
of different contracted triangles have
common endpoint); furthermore, no edge of ${\mathcal C}_{\text{\rm edges}}$ has common endpoint
with any edge of the contracted triangles.
\end{itemize}
In this case, we say that the nodal tropical curve $\gamma(0)$
is an {\it immediate degeneration}
\ignore{GM: changed to an "immediate degeneration"
for consistency with the classical case where
degeneration can be arbitrary and not restricted to
disjoint edges/triangles.}
of $\gamma(1)$ (and $\gamma(1)$ is an {\it immediate perturbation}
of $\gamma(0)$).

\subsection{Combinatorial patchworking of real nodal curves}\label{nodal_patchworking}

Let $\Gamma$ be a nonsingular tropical curve in $\R^2$. A {\it real structure} on $\Gamma$
is given by a collection $T$ of bounded edges of $\Gamma$ which satisfy the following condition:
\begin{itemize}
\item for any cycle of $\Gamma$, denote by $e_1$, $\ldots$, $e_\ell$
the edges of the cycle that belong to $T$;
then, one has
\begin{equation}\label{trop-planarity}
\sum_{i=1}^\ell u_i = 0\quad \text{mod }2,
\end{equation}
where $u_1,\ldots, u_\ell$  are primitive integer vectors in the directions of $e_1$, $\ldots$, $e_\ell$, respectively.
\end{itemize}
Such a collection $T$ is called {\it twist-admissible}, and the edges of $T$ are called {\it twisted}.

To each nonsingular tropical curve $\Gamma$ in $\R^2$ and each twist-admissible collection $T$ of edges of $\Gamma$
we associate a smooth curve $C(\Gamma, T) \subset (\R^\times)^2$
using the following procedure.

\begin{itemize}
\item At each vertex of $\Gamma$, we draw three arcs as depicted in Figure
\ref{trop-picture1}.

\begin{figure}[ht]
\centering
\input{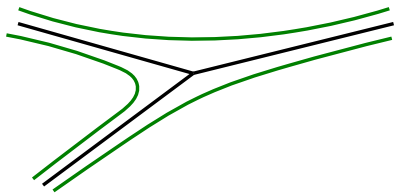}%
\caption{Three arcs at a vertex}
\label{trop-picture1}
\end{figure}

\item For each bounded edge $e$ of $\Gamma$
we join the two corresponding arcs of one endpoint of $e$ to the
two corresponding arcs of the other endpoint of $e$ in the following way:
if $e \notin T$, then join these arcs as depicted in Figure \ref{trop-picture2};
if $e \in T$, then join these arcs as depicted in Figure \ref{trop-picture3}.
Denote by ${\widetilde C}(\Gamma, T)$ the curve obtained.

\begin{figure}[ht]
\centering
\input{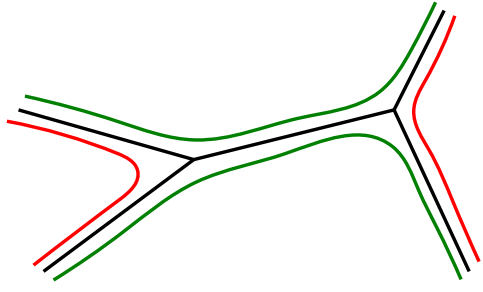}%
\caption{Arcs for an edge which does not belong to $T$}
\label{trop-picture2}
\end{figure}

\begin{figure}[ht]
\centering
\input{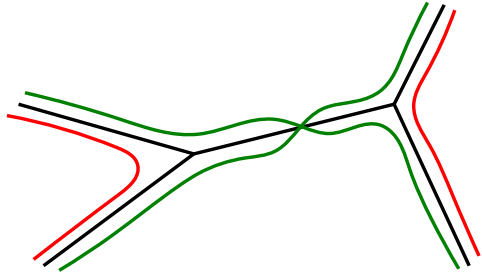}%
\caption{Arcs for an edge which belongs to $T$}
\label{trop-picture3}
\end{figure}

\item Choose arbitrarily a branch of ${\widetilde C}(\Gamma, T)$
and a pair of signs (each sign being $+$ or $-$) for this branch.

\item Associate pairs of signs to all branches of ${\widetilde C}(\Gamma, T)$:
for each edge $e$ with primitive integer direction $(u_1, u_2)$,
the pairs
of signs of the two branches of ${\widetilde C}(\Gamma, T)$ corresponding to $e$
differ by the factor $((-1)^{u_1},(-1)^{u_2})$.
The compatibility condition
\eqref{trop-planarity} ensures that the rule is consistent.

\item Map each branch of ${\widetilde C}(\Gamma, T)$ to $(\R^\times)^2$
by $(x,y)\mapsto (\epsilon_1 e^x, \epsilon_2 e^y)$, where
$(\epsilon_1, \epsilon_2)$ is the pair of signs
associated to the branch.
Denote by $C(\Gamma, T)$
the union of the images of all branches
of ${\widetilde C}(\Gamma, T)$.
\end{itemize}

The isotopy type of the curve $C(\Gamma, T) \subset (\R^\times)^2$
is determined by $\Gamma$ and $T$ up to axial symmetries.

If $\Gamma$ has $\Delta$ as Newton polygon, denote
by ${\overline C}(\Gamma, T)$
the closure of $C(\Gamma, T) \subset (\R^\times)^2 \subset \R Tor(\Delta)$
in the real part $\R Tor(\Delta)$ of the toric surface $Tor(\Delta)$ associated with $\Delta$.
In particular, if $\Gamma$ is of projective
degree $d$, then ${\overline C}(\Gamma, T) \subset \rp^2$.

\ignore{GM: paragraph rephrased, some standard
definitions removed (commented)
\endgraf It: I put back $A$ instead of $\R A$ for consistency with the rest of the paper,
but we can discuss it on Tuesday} 
An algebraic curve $A$ in $(\R^\times)^2$
is a real Laurent polynomial
in two variables well-defined up to multiplication by a monomial.
Such a polynomial has a zero locus in $\R A \subset (\R^\times)^2$.
The Newton polygon $\Delta$ of the polynomial is also
called the Newton polygon of 
$A$. 
Denote with $\R \bar A$ the
closure of $\R A$ in $Tor(\Delta)$.

The combinatorial patchworking
(a particular case of the Viro patchworking theorem \cite{Viro})
can be reformulated in terms of twist-admissible
collections as follows.

\begin{theorem}[cf. \cite{Viro}]
\label{thm:viro}
Let $\Gamma$ be a nonsingular tropical curve in $\R^2$,
and let $\Delta$ be a Newton polygon of $\Gamma$.
Then, for any twist-admissible collection $T$ of $\Gamma$,
there exists
a nonsingular real algebraic curve $A$ in $(\R^\times)^2$ of Newton polygon $\Delta$
such that the pairs $((\R^\times)^2, \R A)$ and $((\R^\times)^2, C(\Gamma, T))$
are homeomorphic. 
Furthermore, the pairs $(\R Tor(\Delta),\R\bar A)$
and \linebreak $(\R Tor(\Delta), {\overline C}(\Gamma, T))$ 
are also homeomorphic. 
\end{theorem}

E.g. if $\Gamma$ is 
projective degree $d$, then
there exists a nonsingular curve
$\R {\overline A}$ of degree $d$ in $\rp^2$
such that the topological pairs $(\rp^2, \R{\overline A})$
and $(\rp^2, {\overline C}(\Gamma, T))$
are homeomorphic.
A reformulation similar to Theorem \ref{thm:viro} was used by B. Haas in \cite{Haas}
for characterization of M-curves obtained by
combinatorial patchworking.

Notice that the empty collection of edges is always twist-admissible.
The resulting nonsingular real algebraic curves are called {\it simple Harnack curves}.
They were introduced in \cite{Mik11}.

\begin{figure}[ht]
\centering
\input{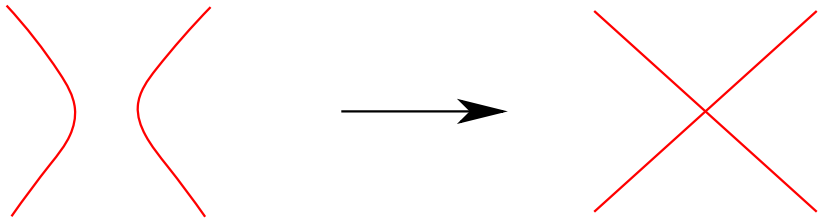}%
\caption{Replacing two arcs by a "cross"}
\label{trop-double}
\end{figure}

\begin{figure}[ht]
\centering
\input{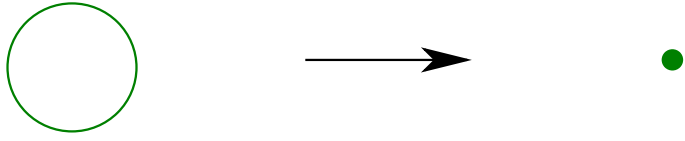}%
\caption{Contracting an oval}
\label{trop-solitary}
\end{figure}

\ignore{GM: paragraph rephrased}
Let $\Gamma$ be a nonsingular tropical curve in $\R^2$,
and let $\Gamma'$ be an immediate degeneration of $\Gamma$
which is a nodal tropical curve.
Then, the underlying graph of $\Gamma'$ is obtained from 
the underlying graph of $\Gamma$
by contracting 
a collection ${\mathcal C}_{\text{\rm edges}}=
\{e_1, \ldots, e_r\} $
of pairwise-disjoint edges 
as well as
a collection ${\mathcal C}_{\text{\rm triangles}}
=\{tr_1,\ldots, tr_s\}$ of triangles.
Suppose that $\Gamma$
is enhanced with such a real structure
that all $e_j$, $j=1,\dots,r$,
are twisted
while all
the triangles $tr_j$, $j=1,\dots,s$,
are composed of non-twisted edges.
%
For each $i = 1$, $\ldots$, $r$, the curve $C(\Gamma, T)$
contains two arcs $a_i$ and $b_i$
which are associated with two endpoints of $e_i$ but
correspond to edges different from $e_i$
(in Figure \ref{trop-picture3} such arcs are shown in red).
The condition \eqref{trop-planarity}
implies that these two arcs are contained
in the same quadrant of $(\R^\times)^2$.
Furthermore, for each $j = 1$, $\ldots$, $s$,
the curve $C(\Gamma, T)$ contains an oval $o_j$
corresponding to the triangle $tr_j$.
Denote by
$C'(\Gamma, T, {\mathcal C}_{\text{\rm edges}},
{\mathcal C}_{\text{\rm triangles}}) \subset (\R^\times)^2$
the subset
obtained from the curve $C(\Gamma, T)$
by replacing each pair of arcs $a_i$, $b_i$
with a "cross" as in Figure \ref{trop-double}
and contracting each oval $o_j$ (where $j = 1$, $\ldots$, $s$)
to a point as in Figure \ref{trop-solitary}.
The isotopy type of
$C'(\Gamma, T, {\mathcal C}_{\text{\rm edges}},
{\mathcal C}_{\text{\rm triangles}}) \subset (\R^\times)^2$
is determined up to axial symmetries by
the tropical curve $\Gamma$, the real structure $T$,
and the collections
${\mathcal C}_{\text{\rm edges}}$ and
${\mathcal C}_{\text{\rm triangles}}$.
Once again, if $\Gamma$ has $\Delta$ as Newton polygon,
denote by ${\overline C}'(\Gamma, T, {\mathcal C}_{\text{\rm edges}},
{\mathcal C}_{\text{\rm triangles}})$ the closure of
$C'(\Gamma, T, {\mathcal C}_{\text{\rm edges}},
{\mathcal C}_{\text{\rm triangles}})$
in $\R Tor(\Delta)$.

\ignore{GM: referenced to the theorem / its proof simplified}
The following theorem is a corollary
of the Viro patchworking theorem \cite{Viro}
for an appropriate $(r + s)$-dimensional family of curves.
It also can be viewed as a special case of
\cite{Sh1}. 

\begin{theorem}[cf. \cite{Viro}]
\label{thm:viro_nodal}
Let $\Gamma$ be a nonsingular tropical curve of Newton polygon $\Delta$ in $\R^2$,
and let $\Gamma'$ be a nodal tropical curve
obtained as an immediate degeneration of $\Gamma$.
Denote by ${\mathcal C}_{\text{\rm edges}}$ {\rm (}respectively, ${\mathcal C}_{\text{\rm triangles}}${\rm )} the
collection of edges
{\rm (}respectively, of triangles{\rm )} of $\Gamma$ that are contracted
in the degeneration of $\Gamma$ to $\Gamma'$.

Then, for any twist-admissible collection $T$ of $\Gamma$ such that
${\mathcal C}_{\text{\rm edges}} \subset T$ and no edge of the triangles in ${\mathcal C}_{\text{\rm triangles}}$ is in $T$,
there exists
a nodal real algebraic curve $\R A\subset (\R^\times)^2$
with Newton polygon $\Delta$ and $r + s$ nodes
such that the pairs 
$((\R^\times)^2, C'(\Gamma, T, {\mathcal C}_{\text{\rm edges}}, {\mathcal C}_{\text{\rm triangles}}))$
and $((\R^\times)^2, \R A)$ 
are homeomorphic.
Furthermore, the pairs
$(\R Tor(\Delta),
\overline{C}'(\Gamma, T, {\mathcal C}_{\text{\rm edges}},
{\mathcal C}_{\text{\rm triangles}}))$
and $(\R Tor(\Delta), \R \bar A)$ are also homeomorphic; 
if in addition $\Gamma$ is irreducible,
then ${\overline A}$ can be chosen irreducible.
\end{theorem} 

\section{Restrictions} \label{sec:restrictions}

In this section, we start the proof of theorem \ref{maintheorem}. 
In fact, we obtain a more detailed classification of nodal rational curves
of degree $5$ in $\rp^2$: we describe all possible complex schemes of these curves.
The {\em complex scheme} of a nodal rational curve $C$ in $\rp^2$
is the topological type of the pair $(\rp^2,\R C)$ 
equipped with one of two complex orientations of $\R C$.
Recall that given a half of the normalization $\psi : \widehat{C} \to C$, we get
a complex orientation for each elliptic node $p \in \R C$
that is the local orientation of $\rp^2$ at $p$ such that the 
intersection at $p$ of $\rp^2$ and the image under $\psi$ of the chosen half
is positive.


\subsection{Possible complex schemes of small perturbations}\label{possible_schemes} 

Let $C$ be a 
nodal rational curve of degree 5 in $\rp^2$, 
and let $C_\circ$ be a small perturbation of $C$
(see Brusotti's theorem \ref{Brusotti}),
such that all hyperbolic nodes are smoothed according to the complex orientations,
each elliptic node is smoothed into an oval, and all complex conjugated nodes are kept.  
%
%
Since $C$ is of type I,
the perturbation $C_\circ$ is also of type I (see Proposition \ref{nodal_perturbation}).
Furthermore, the complex orientations of $\R C_\circ$
are induced by the complex orientations of $\R C$.
Denote by $l$ the number of connected components of $\R C_\circ$.
Recall that $\sigma$ stands for the number 
of imaginary nodes
of $C$ resulted as the intersection of the images
of ${\widehat C}_+$ and ${\widehat C}_-$
under the normalization map $\psi: \widehat C \to C$,
where ${\widehat C}_\pm$ are the two connected components of $\C \widehat C \setminus \R {\widehat C}$. 

\begin{proposition} \label{prop:complexschemes}
The complex scheme of $C_\circ$ is one of the seven schemes
listed in the table \ref{complexschemes}.
	Moreover,
	\begin{itemize}
		\item
		  if $l=7$, then $C_\circ$ is a nonsingular $M$-curve, and hence satisfies the convexity properties 
		  \ref{five_ovales}, \ref{conic_order}, \ref{convexity_maximal}, \ref{cyclic_order},
		\item
		  if $l=5, \sigma=0$, the single positive oval of $\R C_\circ$ is contained in the triangle formed by the
			three negative ovals {\rm (}in particular, the four ovals are in non-convex position{\rm )}, 
		\item
		  if $l=5, \sigma=2$, then the four ovals of $\R C_\circ$ are in convex position. 
	\end{itemize}
  \begin{table}[ht]
	\centering
  \begin{tabular}{@{}llcccc@{}}
    \toprule
    $\sigma$ & $l$ & 7 & 5 & 3 & 1 \\ \cmidrule(l){2-6}
    0        &     & \input{files/pics/scheme33.TpX} & \input{files/pics/scheme13.TpX} & \input{files/pics/schemenest.TpX} &   \\
    2        &     &   & \input{files/pics/scheme22.TpX} & \input{files/pics/scheme02.TpX} &   \\
    4        &     &   &  & \input{files/pics/scheme11.TpX} &   \\
    6        &     &   &   &   & \input{files/pics/scheme00.TpX} \\ \bottomrule
  \end{tabular}
	\caption{Seven possible complex schemes of $\R C_\circ$}
	\label{complexschemes}
  \end{table}
\end{proposition}

\begin{proof}
  The numbers $l,h,e,c,\sigma$ satisfy some straightforward relations.
	First, we have $h + e + c = 6$, the total number of nodes of $C$.
  Second, from lemma \ref{parity-Kcirc} and the fact that each elliptic node gives rise to an oval of $\R C_\circ$
	we can deduce that $l$ is odd and $1+e \leq l \leq 7-c$ (as $6-c$ is the number of real nodes of $C$).
  Moreover, we obviously have $0 \leq \sigma \leq c$. These inequalities explain the 
upper-triangular
	shape of the table. 
	Let us now smooth the complex conjugated nodes of $C_\circ$ in order to obtain a nonsingular curve $C'$.
	Note that $C'$ is of type I if and only if $\sigma=0$. Therefore, according to the rigid isotopy classification \ref{rigid_classification},
the complex scheme of $\R C'$, and hence also the complex scheme of $\R C_\circ$,
is completely determined by $l$ and $\sigma$, namely
	$<J \sqcup 6>_I$, $<J \sqcup 4>_I$ or $<J\sqcup 1<1>>_I$ for $\sigma = 0$
	and $<J \sqcup 4>_{II}$, $<J \sqcup 2>_{II}$ or $<J>_{II}$ for $\sigma \neq 0$.
	As explained before, $C_\circ$ is of type I by Proposition \ref{nodal_perturbation} 
	and therefore satisfies the complex orientation formula \ref{Hilbert16_nodal}.
	Applied to our case, when $l$ and $\sigma$ are fixed, this formula in fact uniquely determines the complex orientations,
	as depicted in table \ref{complexschemes}.
	In particular, if $l=1$ the formula implies $\sigma=6$.
	Finally, the convexity statements follow from the corresponding statements for $C'$ (using proposition \ref{quatre_ovales}).
	In the non-convex case, we also use Fiedler's alternation rule \ref{alternation}.
	Consider a triangle spanned by three ovals and containing the fourth one.
	Applying \ref{alternation} to a pencil of lines with base point in the fourth oval, for example,
	we see that the three outer ovals must have the same orientation.
	This finishes the proof.
\end{proof}
\ignore{GM: Does the last type (with $\sigma=6$) exist?}
\ignore{JR: Yes, perturbation of union of line and empty quartic}

In what follows, we go through these seven cases and study which smoothing diagrams
we can get for each complex scheme.

\subsection{General restrictions on smoothing diagrams}

Let $\Delta_C=(C_\circ; \; \bigcup_p I_p)$ be a smoothing diagram of
real nodal rational curve $C$ of degree $5$ in $\rp^2$.
As above, $C_\circ$ is
a type I small perturbation of $C$
and $\bigcup_p I_p$ is a union of $h$ vanishing cycles.
The following pieces of terminology will be useful.

\begin{definition}
  \label{def:isolatedirred}
  Let $\Delta = (L;I)$ be a smoothing diagram
	(cf.\ Definition \ref{def-coherent-membrane}).
	A connected component of $L$ which does not
	intersect any of the vanishing cycles in $I$
	is called an \emph{isolated oval}.
	Let $\Delta^*$ be the smoothing diagram
	obtained from $\Delta$ after removing all
	isolated ovals. We call $\Delta$ \emph{irreducible}
	if $\Delta^*$ is the smoothing diagram
	of a \emph{single} immersed circle
	(cf.\ Proposition \ref{prop-DeltaK}).
\end{definition}

The following proposition recollects straightforward
properties of smoothing diagrams of real rational curves
(cf. discussions in
subsections \ref{backtosmoothovals} and \ref{immersiongraph}).

\ignore{GM: restated}
\begin{proposition} \label{general_restrictions}
The smoothing diagram $\Delta_C$ satisfies the following properties.
\begin{enumerate}
	\item
	  The underlying graph $\Gamma(\Delta_c)$ has $e$
isolated vertices and one further connected component
		with $h+e+1-l = 7-c-l$ cycles.
	\item
	  All vanishing cycles of $\Delta_C$ are coherent (cf.\ Definition \ref{def-coherent-membrane}).
		This means that a vanishing cycle $I_p$ connects either
		\begin{itemize}
			\item
			  the pseudoline and a negative oval,
			\item
			  or a positive and a negative oval, if they are unnested,
			\item
			  or a negative injective pair of ovals.
		\end{itemize}
  \item
	  The smoothing diagram $\Delta_C$ is irreducible (cf.\ Definition \ref{def:isolatedirred}).
	\item
	  The interior of any isolated oval is empty.
\end{enumerate}
\end{proposition}

In addition to these purely topological properties, further general restrictions
in the case of algebraic curves are provided, by the B\'ezout theorem.
To formulate these restrictions, we introduce the concept of a quasiline.
Let $\Delta=(C_\circ; \; \bigcup_p I_p)$ be a smoothing diagram,
and let $L$ be a non-contractible smoothly
embedded closed 1-submanifold
of $\R\P^2$.
We say that $L$ is a \emph{quasiline} with respect to $\Delta$
if the following conditions hold:
\ignore{GM: Do we want transversal
intersections of $L$?
Added the corresponding requirement.}
\ignore{JR: I don't think we really need it, but maybe it clarifies the idea of pseudolines to require it}
\begin{itemize}
\item The intersection of $L$ and $C_\circ$ is transverse,

\item
	  for each (closed) interval $I_p$,
we have either $|L \cap I_p| \leq 1$ or $L \cap I_p = I_p$,
in the former case $L$ and $I_p$ intersect transversally,
and $L$ is disjoint from the endpoints of $I_p$,
	\item
	  let $r$ be the number of intervals $I_p$ with $L \cap I_p \neq  \emptyset$, and let $s$
		be the cardinality of the intersection $L \cap (C_\circ \setminus \bigcup_p I_p)$;
		then
		\begin{equation}
		\label{eq:quasiline}
		  L.\Delta := 2r + s \leq 5.
		\end{equation}
\end{itemize}

\begin{lemma} \label{bezoutrestriction}
Let $\Delta_C=(C_\circ; \; \bigcup_p I_p)$ be a smoothing diagram of real nodal rational curve $C$
of degree $5$ in $\rp^2$.
	Then, for any two points $p_1, p_2 \in \R\P^2$, there exists a quasiline with respect to $\Delta_C$
	passing through the points.
	If several pairs of points are given, all these quasilines can be chosen
	in such a way that any
	two of them
	intersect either in a single point transversally or
	along a single interval $I_p$.
\end{lemma}

\begin{proof}
	Let $C^t, t \in [0,1]$ be  family of curves such that $C^0 = C$, $C^1 = C_\circ$, and
	such that for $t>0$ this family forms an rigid isotopy.
	Let us fix small discs $D$ around the hyperbolic nodes of $C$.
	Then there exists $s \in (0,1]$
  such that the intersection of any of the discs
	$D$ with $\R C^t$ for $t \in (0,s]$ consists of two embedded intervals.
  Let $I_p^t$ be a corresponding family of vanishing cycles connecting these intervals
	such that $I_p^0 = \{p\}$ is a node of $C$. Then, for each $t \in (0, s]$,
	the tuple $(\R C^t;\bigcup_p I^t_p)$ is equivalent to $\Delta_C$ up to isotopy.
	For each point $p_i$ appearing in the statement, we choose a path $p_i^t$ such
	$\R C^t \cup \bigcup_p I_p \cup p^t_i$ is isotopic to $\Delta_C \cup p_i$ for all
	$t \in (0, s]$.
	In particular, if $p_i \in I_p$, then $p_i^0$
	coincides with the corresponding node $p$ of $C$.
	Let us now prove the existence of quasilines through $p_1,p_2$.
	Start with the (honest) lines $L^t := \overline{p_1^t p_2^t}$.
	By  B\'ezout's theorem, we have $| \R L^0 \cap \R C | \leq 5$.
	We would like to show that $L^t . (\R C^t; \; \bigcup_p I_p) \leq | \R L^0 \cap \R C |$.
	As we consider small perturbations of $C$ and $L^0$, the statement is clear
	if $L^0$ does not contain any hyperbolic node of $C$.
	Hence, assume that $p \in L^0$ is a hyperbolic node of $C$, and let $D$ be the disc around $p$.
	Suppose first that $L_0$ intersects both branches of $C$
at $p$ transversally (so that
	the local intersection multiplicity is $m_p = 2$).
	Depending on whether the two points which form the intersection of $L^t$
	with the boundary $\partial D$ of $D$ are connected in $D \setminus C^t$ or not, we may replace
	the segment $L^t \cap D$ by a path which intersects $I_p$ in a single point or a path which contains $I_p$.
	In both cases, the local contribution to \eqref{eq:quasiline} is $2$,
	so we constructed a quasiline with the required properties.
	For $m_p > 2$, note that any two points in $\partial D \setminus C_t$
	can be joined via a path traversing or containing $I_p$
	with local contribution to \eqref{eq:quasiline} less than $3$.
	Finally, the statement about several lines obviously follows from the fact that two
	(honest) lines intersect in a single point. Indeed, even if this intersection point is a node of $C$, the above procedure can be performed such that
	the two paths intersect either in a single point or in the interval $I_p$.
\end{proof}

\begin{remark} \label{quasilines-through-nodes}
  Similarly, we can show the following statement:
	Given an interval $I_p$, there exists
	\begin{itemize}
		\item
		 a quasiline intersecting
the interval $I_p$ transversally in a single interior point,
		\item
		 a quasiline containing  the interval.
	\end{itemize}
	This corresponds to choosing $L_0$ to be a line intersecting the corresponding node
	transversally and lying in the corresponding ``quadrants of the node''.
Clearly, we may assume that the
pairwise intersections of such quasilines
are also either single points or the interval $I_p$ itself.
\ignore{GM: rephrased}
\end{remark}

\subsection{\texorpdfstring{$M$}{M}-curve case (\texorpdfstring{$l=7$}{l=7})} \label{subsec:restr:Mcurve}

In the M-curve case $l=7$, we have $\sigma=c=0$, and all necessary further restrictions can be summarized in
Proposition \ref{mcurvediagrams} below.
Let $\Delta$ be the following diagram (in the sense of definition
\ref{def-coherent-membrane}):

\begin{minipage}{\textwidth}
	\centering
  \input{files/pics/harnack.TpX}
\end{minipage}

A \emph{subdiagram} of $\Delta$ is a diagram obtained from $\Delta$ by removing
some of the (dashed) intervals. A diagram $\Delta'$
is called a \emph{tree} (or {\em arboreal} cf. Definition
\ref{def-arboreal}) if the underlying graph
$\Gamma(\Delta)$ is a (connected) tree
plus possibly some isolated points.

\begin{proposition}\label{mcurvediagrams}
  If $l = 7$ (i.e.\ if $C_\circ$ is a M-curve), then the smoothing diagram $\Delta_C$ is
	(up to isotopy) a subtree of $\Delta$.
	Moreover, table \ref{list3} lists all 32 possible choices of such subtrees.
\end{proposition}

\begin{proof}
  We will prove the following properties of $\Delta_C$ which, together with corollary \ref{cyclic_order},
	imply the statement.
	\begin{enumerate}
		\item
		  $\Delta_C$ is a tree (in particular, it does not have double edges or loops).
		\item
		  Two ovals connected by a vanishing cycle must be neighbors in the sense of subsection \ref{rigid_isotopies}.
		\item
		  An oval connected to the pseudoline of $C_\circ$ by a vanishing cycle must be negative.
		\item
		  Assume that there are three vanishing cycles attached to the pseudoline of $C_\circ$.
			Then when moving along $J$, the vanishing cycles lie on alternating sides.
	\end{enumerate}
	Part (a) follows from Proposition \ref{general_restrictions} (a).
	\\
	For (b), let us consider two ovals $O_1$ and $O_2$ connected
by a vanishing cycle.
	Let $O_3$ be an arbitrary third oval of $C_\circ$,
and let $x_1, x_2, x_3$ be points in their respective interiors.
	We want to show that $x_1$ and $x_2$ are neighbors
viewed from $x_3$ according to the definition on page
	\pageref{def:neighbors}.
To do so, we define the path $\sigma_{1,2}^3$ to be the line segment
	$[x_1,x_2]_{C_\circ}$ and we define ${\cal S}^k_{i, j}$ to be the segment of
	the pencil of lines through $x_3$ which intersect $[x_1,x_2]_{C_\circ}$.
	With these choices, the second condition for being neighbors is obviously satisfied, and it remains to show that
	the lines in ${\cal S}^k_{i, j}$ do not intersect any other oval than $O_1, O_2, O_3$.
  For this we might have to replace $C_\circ$ by a deformation which is sufficiently close to $C$.
  More precisely, let $C^t, t \in [0,1]$ be  family of curves such that $C^0 = C$, $C^1 = C_\circ$, and
	such that for $t>0$ this family forms an rigid isotopy.
	We can proceed as above for each $t$, with the additional assumption that the chosen points
	$x^t_1, x^t_2$ converge to the node $p$ when $t$ goes to $0$. We now want to show that
	the lines in the segment ${\cal S}^{k,t}_{i, j}$ do not intersect any other oval of $C^t$
	for sufficiently small $t$. If this was not true, when $t$ goes to $0$ this would imply
	the existence of a  line passing through $p$ and at least two other ovals/nodes
	of $C$, which contradicts the B\'ezout theorem.
	\\
	Part (c) immediately follows from Proposition \ref{general_restrictions} (b).
	\\
	For (d), let us choose a point in the interior of each positive oval of $C_\circ$, and consider the three
	straight lines passing each through two of the chosen points.
	These three lines divide $\rp^2$ in four regions.
One of these regions is the convex triangle
	with vertices at the chosen points. Propositions \ref{five_ovales} and \ref{convexity_maximal} imply
	that each of the three other regions contains a negative oval of $C_\circ$.
	Furthermore, each of these three regions is cut by
	the pseudoline of $C_\circ$ into two sub-regions. We call such a sub-region
	\emph{triangular} (respectively, \emph{quadrangular}) if it is adjacent to two (respectively, three)
	straight lines passing through the chosen points.
	Proposition \ref{convexity_maximal} implies that
	each negative oval of $C_\circ$ is contained in
a quadrangular sub-region.
As in the proof of part (b)
the vanishing cycles must be disjoint from
our three straight lines.
	Thus the three quadrangular sub-regions
are adjacent to the pseudoline from alternating sides,
	which finishes the proof.

  \begin{minipage}{\textwidth}
	  \centering
    \input{files/pics/harnacklines.TpX}
  \end{minipage}
\end{proof}

\subsection{The Hyperbolic Curve (\texorpdfstring{$l=3, \sigma = 0$}{l=3, sigma=0})} \label{restr:nestedcurve}

Let us now consider smoothing diagrams which are based on a hyperbolic curve (i.e.\ $l=3, \sigma = 0$).
By Proposition \ref{general_restrictions} (b), such a smoothing diagrams has vanishing cycles of two types,
connecting the pseudoline and the outer oval or connecting the nested ovals.
For now, let us remove the latter ones and
classify the remaining diagrams up to isotopy.
For this it is useful to use the language of immersion graphs
from Definition \ref{enhanced-Gamma}
(even though for visual convenience we stick
to more intuitive smoothing diagrams in the pictures).
In our situation (after removing the inner oval
and adjacent vanishing cycles), we are left with
immersion graphs $\Gamma$
on two vertices: the root vertex $w$
and a second vertex $v$. All edges are directed from $v$ to $w$.
The number of edges is $1 \leq t \leq 5$. The ribbon structure of $\Gamma$ gives us two cyclic orderings of the edges.
Finally, the projective enhancement of $\Gamma$ (the crosses on the segments of $w$) can be encoded in a vector
$T \in \Z_2^t$, which contains an entry $1$ for each crossed segment. This vector $T$ is well-defined up to cyclic reordering
and if $|T|$ denotes the sum of the entries of $T$, we have
\begin{equation}
  |T| \equiv 1 \mod 2.
\end{equation}
From the compatibility property of the two cyclic orderings in proposition \ref{clusters} it follows
that $\Gamma$ is in fact completely determined by $T$. Moreover, by proposition \ref{criteria-odddegree}
the isotopy type of the reduced smoothing diagram is determined by $\Gamma$.
Hence the list of possible reduced smoothing diagrams corresponds to the list of possible vectors
$T$ and looks as follows.

\begin{minipage}{\textwidth}
	\centering
\input{files/pics/nested1.TpX}
\input{files/pics/nested01.TpX}
\input{files/pics/nested001.TpX}
\input{files/pics/nested111.TpX}
\end{minipage}
\begin{minipage}{\textwidth}
  \centering
\input{files/pics/nested0001.TpX}
\input{files/pics/nested0111.TpX}
\input{files/pics/nested00001.TpX}
\input{files/pics/nested00111.TpX}
\end{minipage}
\begin{minipage}{\textwidth}
  \centering
\input{files/pics/nested01011.TpX}
\input{files/pics/nested11111.TpX}
\end{minipage}

Given such a reduced diagram, the outer oval is subdivided into $t$ segments.
Again by proposition \ref{clusters}, the isotopy type
of the full diagram is determined by the data to which segments the $h-t$ ``inner'' vanishing cycles
are attached.

Many of the choices for attaching the vanishing cycles
are prohibited by the irreducibility condition of
Proposition \ref{general_restrictions} (c). The typical situation is depicted
in the following picture.

\begin{minipage}{\textwidth}
	\centering
  \input{files/pics/red.TpX}
\end{minipage}

It shows a component $K$ together with a sequence $p_0, p_1, \ldots, p_n=p_0$ of pairs of vanishing cycles $
p_i = \{I_i, I'_i\}$ attached to it. Each pair $p_i$ connects $K$ to the same component $K_i$.
Moreover, the vanishing cycles following $I'_i$ with respect to the cyclic ordering on $K$ is $I_{i+1}$,
and the vanishing cycles following $I_i$ with respect to the cyclic ordering on $K_i$ is $I'_i$
(as indicated by the bold segments here).
Such a diagram violates Proposition \ref{general_restrictions} (c),
as the bold segments together with the vanishing cycles
form an embedded circle which does not cover the whole curve.
We refer to this situation by $(Red)$.
Note that $K$ can be any component of $C_\circ$,
even though we depicted it as an oval in the picture.

We now list all diagrams which satisfy to the properties listed
in Proposition \ref{general_restrictions}
(in particular, to the $(Red)$-rule).
As explained above, given a reduced diagram from the above the list, we need to specify to which segments
the ``inner'' vanishing cycles are attached. Of course, it suffices to consider these choices up to symmetries
of the reduced diagram. For example, for $T=(0,0,1)$,
the two ``lower'' segments of the outer oval are symmetric.
Note that by Proposition \ref{general_restrictions} (d) we have $e=0,1$ (only the inner oval can correspond to an elliptic node).
Let us start with $e=0$.

We start with the case $e=0$ and go through all possible choices for $T$.
For $T=(1)$, there is only one segment to choose, so we get the following three types.

\begin{minipage}{\textwidth}
\centering
\begin{tabular}{lll}
$c=0$ & $c=2$ & $c=4$ \\
\input{files/pics/nested1a.TpX}      & \input{files/pics/cc1o.TpX}      & \input{files/pics/cc2alt.TpX}
\end{tabular}
\end{minipage}

For $T=(0,1)$, there are two segments to choose from.
Note that assigning an even number of vanishing cycles to each of the segments is forbidden $(Red)$.
All other options are depicted below.

\begin{minipage}{\textwidth}
\centering
\begin{tabular}{lll}
$c=0$ &  & $c=2$ \\
\input{files/pics/nested01a.TpX}      & \input{files/pics/nested01b.TpX}      & \input{files/pics/cc1p.TpX}
\end{tabular}
\end{minipage}

For $T=(0,0,1)$, we get three segments,
two of which are symmetric.
All choices satisfy to Proposition \ref{general_restrictions},
except for attaching
one vanishing cycle to each segment (such a smoothing diagram violates the irreducibility condition (c)).

\begin{minipage}{\textwidth}
\centering
\begin{tabular}{llll}
$c=0$ &  &  &  \\
\input{files/pics/nested001a.TpX} & \input{files/pics/nested001b.TpX} & \input{files/pics/nested001c.TpX} & \input{files/pics/nested001d.TpX} \\
      & $c=2$ &  &  \\
\input{files/pics/nested001e.TpX} & \input{files/pics/cc1q.TpX} & \input{files/pics/cc1r.TpX} &
\end{tabular}
\end{minipage}

For $T=(1,1,1)$, all three segments are symmetric.
Note also that in this case, the reduced diagram itself
is not irreducible in the sense
of Proposition \ref{general_restrictions} (c),
but consists of three components containing one of the segments each.
In particular, the completed diagram can
satisfy to the condition (c) only
if the additional vanishing cycles are attached to all
the three segments. Hence, the only possible type is as follows:

\begin{minipage}{\textwidth}
\centering
\input{files/pics/nested111a.TpX}
\end{minipage}

From now on we assume that $t \geq 4$ (hence $c=0$).
For $T=(0,0,0,1)$, we get a $(Red)$-situation whenever we attach the two inner vanishing cycles to the same or opposite segments.
Two possibilities (up to symmetries) remain:

\begin{minipage}{\textwidth}
\centering
\input{files/pics/nested0001a.TpX}
\input{files/pics/nested0001b.TpX}
\end{minipage}

For $T=(0,1,1,1)$, we have to attach exactly one vanishing cycle to the ``lower'' segment of the outer oval. In any other case, we get a
$(Red)$-situation (if we attach no vanishing cycle to the lower segment, $K$ in $(Red)$
can be chosen to be the pseudoline).
Again, only two possibilities remain:

\begin{minipage}{\textwidth}
\centering
\input{files/pics/nested0111a.TpX}
\input{files/pics/nested0111b.TpX}
\end{minipage}

For $t = 5$, the inner oval is attached to the reduced diagram via a single vanishing cycle.
Hence a diagram is not irreducible if and only if
the reduced diagram does.
The choice $T=(0,0,1,1,1)$ violates the property and therefore can be ignored.
For all other choices of $T$, we get the following list of possibilities (up to symmetries).

\begin{minipage}{\textwidth}
\centering
\input{files/pics/nested00001a.TpX}
\input{files/pics/nested00001b.TpX}
\input{files/pics/nested00001c.TpX}
\end{minipage}
\begin{minipage}{\textwidth}
  \centering
\input{files/pics/nested01011a.TpX}
\input{files/pics/nested01011b.TpX}
\input{files/pics/nested01011c.TpX}
\end{minipage}
\begin{minipage}{\textwidth}
  \centering
\input{files/pics/nested11111a.TpX}
\end{minipage}

So far, we only applied the restrictions given by
Proposition \ref{general_restrictions}.
However, by applying B\'ezout's theorem
we can prohibit some more diagrams
in this list.
The following three
subdiagrams cannot be contained in
the smoothing diagram of
a real rational quintic curve
by Lemma \ref{bezoutrestriction}.
The diagram also specify the points
that should be connected
by quasilines.

\begin{minipage}{\textwidth}
\centering
\input{files/pics/bezouta.TpX}
\input{files/pics/bezoutb.TpX}
\input{files/pics/bezoutc.TpX}
\end{minipage}

In the first two cases it is impossible to draw a quasiline
through the specified pair of points.
In the third case, individual quasilines can be drawn,
but
they violate the second part
of Lemma \ref{bezoutrestriction}.
Hence, any diagram which is an extension of the three diagrams above is prohibited.
This kills the second diagram for $T=(0,1)$ diagrams number 3,4,5 for for $(0,0,1)$,
the second diagram for $T=(0,0,0,0,1)$, and finally diagrams number 3 and 6 for $t=5$.
The remaining diagrams are contained in tables \ref{list1}, \ref{list4} and \ref{list6}.

Let us now turn to $e=1$. In this case, there are
no inner vanishing cycles (i.e.\ $t=h$) and
diagram's irreducibility
only depends on the choice of $T$.
Moreover, $t=h$ must be odd, thus we can ignore $t=2,4$. According to our previous
considerations, we are left with $T=(0,0,0,0,1)$, $T=(0,1,0,1,1)$, $T=(1,1,1,1,1)$,
$T=(0,0,1)$ and $T=(1)$. The five corresponding diagrams are contained in
tables \ref{list1}, \ref{list4} and \ref{list6}.

\subsection{The non-convex 4-Oval Curve (\texorpdfstring{$l=5, \sigma = 0$}{l=5, sigma=0})}

Let us, first, collect some restrictions for this
case in the following proposition.

\begin{proposition} \label{4ovalsdiagrams}
  Let $\Delta_C$ be a smoothing diagram of $C$ with $l=5, \sigma=0$. Then,
	\begin{enumerate}
		\item
    	the positive oval can be connected to each negative oval by at most one vanishing cycle,
		\item
		  if two vanishing cycles connect the same (negative) oval to the pseudoline,
			then the disc they bound does not contain other ovals.
	\end{enumerate}
\end{proposition}

\begin{proof}
Both statements follow from Lemma \ref{bezoutrestriction}.
Indeed, the pairs of points in the following pictures cannot be connected by a quasiline.
Here,
the blue oval on the right in the second picture
is an arbitrarily chosen third oval (which might be contained
in the disc as well).

	\begin{minipage}{\textwidth} \label{4oval_nonconvex}
    \centering
		\begin{tabular}{ll}
			(a) & (b) \\
			\input{files/pics/4ovalnodouble.TpX}		&  \input{files/pics/4ovalbezA2.TpX}
		\end{tabular}
  \end{minipage}
\end{proof}

Consider the graph $\Gamma'$ obtained form $\Gamma(\Delta_C)$
by subsequently removing all zero- and one-valent vertices until
each vertex has valence at least $2$.

\begin{proposition} \label{4ovals_6cases}
  The graph $\Gamma'$ is either empty (when $c=2$) or has genus 2 (when $c=0$).
	In the latter case, there are 5 possible graphs which can occur as $\Gamma'$.
	Moreover, for each $\Gamma'$
its immersion to $\rp^2$ is topologically unique.
	The list of corresponding smoothing diagrams is as follows:

\begin{minipage}{\textwidth}
  \centering
	\resizebox{\textwidth}{!}{%
  \begin{tabular}{ccccc}
    A & B & C & D & E \\
    \input{files/pics/4oval300.TpX} & \input{files/pics/4oval21.TpX} & \input{files/pics/4oval21B.TpX} & \input{files/pics/4oval21C.TpX} & \input{files/pics/4oval111.TpX}
  \end{tabular}
	}
\end{minipage}
\end{proposition}

\begin{proof}
  The first statement follows from
Proposition \ref{general_restrictions} (a).
	The second statement follows from
Propositions \ref{general_restrictions} (b) and
	\ref{4ovalsdiagrams} (a).
	For the third statement, let us fix one of the five graphs for $\Gamma'$ and study the
	possible enhancements.
	First, the orientations are fixed:
all edges adjacent to the root vertex are oriented towards it.
	All other edges are not oriented.
	As before, let us describe the projective enhancement of $\Gamma$ by a vector
	$T \in \Z_2^n$. In each case, we have two choices, namely $T=(0,0,1)$ or $T=(1,1,1)$ for
	A,B,E and $T=(0,0,0,1)$ or $T=(0,1,1,1)$ for C and D.
	Note that removing zero- or one-valent vertices does not affect the irreducibility
	property \ref{general_restrictions} (c) of the corresponding graph/diagram,
	hence the enhancement of $\Gamma'$ is required to be irreducible.
	This excludes the possibilities $T=(1,1,1)$ for A,B, $T=(0,0,0,1)$ for C,D,
	and $T=(0,0,1)$ for E.
	For the opposite choices of $T$, we claim that $T$ determines the enhancement of
	$\Gamma$ completely. In A,B and E, there is only a single cluster except for the root vertex,
	and hence the ribbon structure is determined by Proposition \ref{clusters}.
	In the cases C and D, Proposition \ref{clusters} shows that the two clusters
	are unlinked with respect to the ``two-sided'' cyclic	order on the pseudoline.
	In each case, there are two possibilities, one of which is reducible by $(Red)$.
	By Proposition \ref{criteria-odddegree}, the immersion graph determines
	the (reduced) smoothing diagram completely.
\end{proof}

Let us now go through the six cases (A,B,C,D,E and $c=2$) one by one.

\paragraph{Diagram A.}

Let $\Delta'$ be the diagram obtained from $\Delta_C$ be removing all ovals (respectively, vanishing cycles)
which are not connected (respectively, not attached) to the pseudoline.
The following list shows all such diagrams (up to isotopy):

\begin{minipage}{\textwidth}
  \centering
\input{files/pics/4oval300.TpX}
\input{files/pics/4oval310a.TpX}
\input{files/pics/4oval310b.TpX}
\input{files/pics/4oval311a.TpX}
\input{files/pics/4oval311b.TpX}
\end{minipage}
\begin{minipage}{\textwidth}
  \centering
\input{files/pics/4oval311c.TpX}
\input{files/pics/4oval311d.TpX}
\input{files/pics/4oval311e.TpX}
\input{files/pics/4oval311f.TpX}
\input{files/pics/4oval311g.TpX}
\end{minipage}

Four of these diagrams can be prohibited. Three isotopy types of diagrams
can be excluded by the following convexity argument
(denoted by $(Conv)$).
Choose a point in the interior of the positive oval and a point
in the relative interior of the ``middle''
vanishing cycle
(i.e.\ the one separating the two zeroes in $T=(0,0,1)$).
By Lemma \ref{bezoutrestriction} there exists a quasiline through this pair of points which
(up to isotopy) looks as the vertical line in the following picture:

\begin{minipage}{\textwidth}
  \centering
  \input{files/pics/4ovalconv.TpX}
\end{minipage}

The quasiline and the pseudoline split $\rp^2$ into two regions, and it follows from
the non-convexity property
(see Proposition \ref{quatre_ovales})
that the quasiline can be chosen such that
the remaining two negative ovals do not lie in the same region.
This restriction excludes the diagrams 5, 7, and 10.

We show now that the diagram 4 can also be prohibited.
To do so, choose a point $p$ on the ``middle'' vanishing cycle $I$ as before,
and choose points $q_1, q_2$ in the interior of each ``small'' oval.
Then, we
get a contradiction with existence
of three quasilines: $\overline{p q_1}$, $\overline{p q_1}$
and a quasiline intersecting $I$ in a single point.
The situation is illustrated in the following
picture (where the quasilines $\overline{p q_1}$ and $\overline{p q_1}$ are drawn
as dashed blue lines).

\begin{minipage}{\textwidth}
  \centering
  \input{files/pics/4ovalbezB.TpX}
\end{minipage}

We end up with a list of 6 (incomplete) diagrams. It remains to count the number of ways in which
these diagrams can be completed to full diagrams.
Whenever there is an ambiguity in completion,
 we draw the various possibilities in the same diagram
by dashed lines and write the corresponding ``multiplicity'' next to it.

\begin{minipage}{\textwidth}
  \centering
	\input{files/pics/4oval300+.TpX}
	\input{files/pics/4oval310a+.TpX}
	\input{files/pics/4oval310a++.TpX}
	\input{files/pics/4oval310b+.TpX}
\end{minipage}
\begin{minipage}{\textwidth}
  \centering
	\input{files/pics/4oval310b++.TpX}
	\input{files/pics/4oval311c+.TpX}
	\input{files/pics/4oval311e+.TpX}
	\input{files/pics/4oval311f+.TpX}
\end{minipage}
\begin{minipage}{\textwidth}
  \centering
	\input{files/pics/4oval300ell.TpX}
	\input{files/pics/4oval310aell.TpX}
	\input{files/pics/4oval310bell.TpX}
	\input{files/pics/4oval311cell.TpX}
\end{minipage}
\begin{minipage}{\textwidth}
  \centering
	\input{files/pics/4oval311eell.TpX}
	\input{files/pics/4oval311fell.TpX}
\end{minipage}
\begin{minipage}{\textwidth}
  \centering
	\input{files/pics/4oval300ell2.TpX}
	\input{files/pics/4oval310ell2.TpX}
	\input{files/pics/4oval300ell3.TpX}
\end{minipage}

All these diagrams are contained in table \ref{list2}.

\paragraph{Diagram B.}

For diagram B, it remains to add one negative oval to the pictures. A priori,
it might be isolated (3 connected components of the complement give rise to 3 choices),
connected to the positive oval (2 segments of the positive oval give rise to 2 choices)
or connected to a segment of the pseudoline (3 segments $\times$ 2 sides = 6 choices).
Here are the forbidden choices:

\begin{minipage}{\textwidth}
  \centering
  \input{files/pics/4oval21Adisc.TpX}
  \input{files/pics/4oval21Abez.TpX}
  \input{files/pics/4oval21Aconv.TpX}
\end{minipage}

The first picture is forbidden by
Proposition \ref{4ovalsdiagrams} (b), in the second case no quasiline can be drawn through
the indicated pair of points, and the third case is prohibited by the convexity argument $(Conv)$ with respect to quasilines
through the indicated pair of points. It remains the following 5 cases, which are contained in table \ref{list2}:\\[2ex]

\begin{minipage}{\textwidth}
  \centering
  \input{files/pics/4oval21A.TpX}
	\input{files/pics/4oval21Aell.TpX}
\end{minipage}

\paragraph{Diagram C and E.}

Both diagrams are in fact complete, so they corresponds to exactly two smoothing diagram which are contained in table \ref{list2}.

\paragraph{Diagram D.}

Proposition \ref{4ovalsdiagrams} (b)
provides a restriction on the position of the remaining ovals.
Furthermore, the following way of attaching the third negative oval to the pseudoline is forbidden
by the convexity rule $(Conv)$:

\begin{minipage}{\textwidth}
  \centering
  \input{files/pics/4oval21Cd.TpX}
\end{minipage}

The following list shows all other possibilities to complete diagram D (note that the two ovals in
diagram D are symmetric). All these diagrams are contained in table \ref{list2}.

\begin{minipage}{\textwidth}
  \centering
	\input{files/pics/4oval21Ca.TpX}
  \input{files/pics/4oval21Cb+.TpX}
  \input{files/pics/4oval21Cc+.TpX}
	\input{files/pics/4oval22Cell.TpX}
  \input{files/pics/4oval22ell2.TpX}
\end{minipage}

\paragraph{Case $c=2$.}

In this case, recall that by Proposition \ref{4ovals_6cases}
(or Proposition \ref{general_restrictions} (a))
the smoothing diagram $\Delta_C$ is a tree.
The possibilities for such trees can be easily listed.

\begin{minipage}{\textwidth}
  \centering
  \input{files/pics/cc1f.TpX}
  \input{files/pics/cc1g.TpX}
  \input{files/pics/cc1h.TpX}
  \input{files/pics/cc1uu.TpX}
  \input{files/pics/cc1uuu.TpX}
\end{minipage}

\begin{minipage}{\textwidth}
  \centering
  \input{files/pics/cc1c.TpX}
  \input{files/pics/cc1u.TpX}
  \input{files/pics/cc1x.TpX}
  \input{files/pics/cc1i.TpX}
  \input{files/pics/cc1s.TpX}
\end{minipage}

\begin{minipage}{\textwidth}
  \centering
  \input{files/pics/cc1w.TpX}
  \input{files/pics/cc1y.TpX}
  \input{files/pics/cc1z.TpX}
\end{minipage}

All these diagrams except for the diagrams 3 and 5 are contained in table \ref{list5}.
In order to prohibit diagrams 3 and 5, we need a B\'ezout type argument similar to
Lemma \ref{bezoutrestriction}, but involving a conic.
Namely, let $H$ be the conic passing through
5 nodes of $C$, the two complex conjugated nodes and the three hyperbolic nodes
corresponding to the vanishing cycles attached to the pseudoline.
By B\'ezout's theorem, $C$ and $H$ do not have further intersection points
and the intersection multiplicity at each node is exactly 2.
Thus, there exists a \emph{quasiconic} (an embedded contractible loop in $\rp^2$
intersecting $\Delta_C$ only in the three vanishing cycles attached to the
pseudoline). The following picture shows such a quasiconic (unique up to
isotopy) for diagrams 3 and 5:

\begin{minipage}{\textwidth}
  \centering
  \input{files/pics/bezoutwithconic.TpX} \hspace{0.2\textwidth}
  \input{files/pics/bezoutwithconic2.TpX}
\end{minipage}

In addition, let us now consider the real line $L$ which passes through the two complex conjugated nodes.
It intersects $H$ in the two nodes and hence nowhere else. Moreover, $L$ must intersect
$\R C$ in exactly one point with intersection multiplicity $1$.
Now again, after small perturbations, $L$ gives rise to a quasiline
in $\rp^2$ which does not intersect the quasiconic and intersects $\Delta_C$
only in one point on the pseudoline.
This is clearly impossible, so diagram 3 and 5 are forbidden.

\subsection{The convex 4-Oval Curve (\texorpdfstring{$l=5, \sigma = 2$}{l=5, sigma=2})}

Let us collect all necessary restrictions in the following statement.

\begin{proposition}
  If $l = 5, \sigma=2$ (i.e.\ if $C_\circ$ is the curve with four ovals in convex position),
	then the smoothing diagram $\Delta_C$ is (up to isotopy) a subtree of the following diagram.

\begin{minipage}{\textwidth}
  \centering
  \input{files/pics/convex4oval.TpX}
\end{minipage}

	Moreover, all 11 such subtrees are contained in table \ref{list5}
	(as diagrams 5--11 and 13--16).
\end{proposition}

\begin{proof}
  It follows from Proposition \ref{general_restrictions} (a) that $\Delta_C$ is tree.
	Proposition \ref{general_restrictions} (b) implies that
	$\Gamma(\Delta_C)$ is a subgraph of the graph underlying the depicted
	diagram.
Since $\Delta_C$ is a tree with at most two vanishing cycles
	attached to the pseudoline, the diagram is completely determined
	by the underlying graph.
\end{proof}

\subsection{The remaining cases}

In the case $l=3, \sigma=2, c=2$, the two ovals are both negative and therefore
cannot be connected by a vanishing cycle.
If one oval is connected to the pseudoline by three vanishing cycles,
the $T=(1,1,1)$ pattern is prohibited by
the irreducibility.
The same holds true for the $T=(0,0,0,1)$ pattern when both ovals are connected to the
pseudoline by two vanishing cycles each.
Hence $T=(0,0,1)$ or $T=(0,1,1,1)$, and by Propositions \ref{clusters} and \ref{criteria-odddegree}
this determines the (reduced) smoothing diagram completely (as in the proof of Proposition \ref{4ovals_6cases}).
We get the following cases, contained in table \ref{list4}:

\begin{minipage}{\textwidth}
  \centering
  \input{files/pics/cc1lnew.TpX}
	\input{files/pics/cc1n.TpX}
  \input{files/pics/cc1v.TpX}
\end{minipage}

In all other cases, we have $c \geq 4$ and the combinatorics become very easy.
We get the following 5 cases contained in table \ref{list6}.

\begin{minipage}{\textwidth}
  \centering
  \input{files/pics/cc2a.TpX}
  \input{files/pics/cc2b.TpX}
  \input{files/pics/cc2c.TpX}
  \input{files/pics/cc2d.TpX}
  \input{files/pics/cc3.TpX}
\end{minipage}

\section{Constructions}\label{sec:constructions}

In this section, we show that every smoothing diagram shown in the tables on pages \pageref{list1} -- \pageref{list6}
is the smoothing diagram of some real nodal rational
curve of degree $5$ in $\rp^2$.
This is the final ``construction'' part of the proof of Theorem \ref{maintheorem}.
We
start with the construction of curves with only hyperbolic nodes and add elliptic (respectively,
complex conjugated) nodes later.

\subsection{\texorpdfstring{$M$}{M}-curve (\texorpdfstring{$l=7, h=6$}{l=7, h=6})}

In the $M$-curve case, constructions are straightforward.
We consider the following arrangement of a conic and 3 lines
(with the indicated choice of a complex orientation).

\begin{minipage}{\textwidth}
  \centering
  \input{files/pics/conic3lines.TpX} \hspace{10ex}
  \input{files/pics/mcurvesmoothing.TpX}
\end{minipage}

Let us smooth all 9 nodes according
to these orientations.
The resulting smoothing diagram is indicated on the right hand side of the figure.
Note that this pattern coincides with the ``universal'' smoothing diagram from Proposition
\ref{mcurvediagrams}. Moreover, note that by Theorem \ref{Brusotti} instead of smoothing all 9 nodes we
may keep some of them, which leads to curves whose smoothing diagram is any subdiagram of the universal one.
The irreducibility of these curves is equivalent to the irreducibility of the corresponding smoothing diagram
in the sense of Proposition \ref{general_restrictions} (c).
Hence we can construct the 9 isotopy types with $l=7, h=6$ in table \ref{list3}.

\subsection{Non-convex 4-Oval Curve (\texorpdfstring{$l=5, h=6$}{l=5, h=6})}

In this case we can use the small deformations of five lines.
The pictures are equally easily drawn in the classical and tropical world, and so we will give
both descriptions.
Let us start with tropical pictures and consider the smooth tropical quintic
$B \subseteq \T\P^2$ on the left of the following picture.

\begin{minipage}{\textwidth}
  \centering
  \input{files/pics/bathroomtiling.TpX} \hspace{10ex}
  \input{files/pics/5tropicallines.TpX}
\end{minipage}

We equip $B$ with a real structure by twisting all bounded diagonal edges,
i.e.\ $T= \{$diagonal edges$\}$. Note that $T$ is twist-admissible
(see condition \ref{trop-planarity}).
In the picture, the twists are indicated by small crosses.
%
%
%
Note that all twisted edges are pairwise disjoint (i.e.\ now two of them have endpoints in common).
Hence we may degenerate $B$ to a tropical nodal curve by shrinking \emph{any} collection of pairwise disjoint twisted edges.
In particular, if we shrink all twisted edges, we obtain a collection of five tropical lines intersecting transversally,
as depicted on the right hand side of the above picture.

The next pictures depict the ``real versions'' constructed in section \ref{nodal_patchworking}.
On the left hand, we see $\widetilde{C}(B,T) \subset \R^2$, on the right hand side we find
$\overline{C}(B,T) \subset \rp^2$.
%
Additionally, for each twisted edge we draw the vanishing cycle we obtain by shrinking this edge and creating a hyperbolic node.
Note that the smoothing diagram on the right is more or less directly visible from the tropical curve above.

\begin{minipage}{\textwidth}
  \centering
  \input{files/pics/constr4ovalreal1.TpX} \hspace{10ex}
  \input{files/pics/constr4ovalreal2.TpX}
\end{minipage}

Using an isotopy we may redraw the picture on the right as follows:

\begin{minipage}{\textwidth}
  \centering
  \input{files/pics/constr4ovaluniversal.TpX}
\end{minipage}

Hence, Theorem \ref{thm:viro_nodal} ensures that any irreducible subdiagram of the smoothing diagram above
occurs as the smoothing diagram of an irreducible rational nodal curve $\overline{A}$.
So it only remains to check
that all isotopy types with $l=5, h=6$ in table \ref{list2}
are indeed subdiagrams of the above one.

%

Let us also briefly mention the classical construction.
Consider $5$ real lines in the real projective plane such
that no three among them intersect. The topological type of such a real line arrangement is unique
and looks as follows:

\begin{minipage}{\textwidth}
  \centering
  \input{files/pics/constr5lines1.TpX} \hspace{10ex}
  \input{files/pics/constr5lines2.TpX}
\end{minipage}

There are $4$ possibilities to orient the lines (up to permutation of the lines).
Indeed, moving along the sides of the pentagon, at each vertex
we can choose to keep or to switch orientation and the total number of switches must be even. This can be described
by the four cyclic vectors $T$ of length $5$ which  appeared in subsection \ref{restr:nestedcurve}
(but with an even number of non-zero coordinates).
For three of the choices, the type I small perturbation gives the hyperbolic curve.
Only the choice of orientations displayed above
leads to the 4-oval curve. Moreover, the
smoothing diagram associated to this choice is exactly the ``universal'' one displayed above.

\subsection{Hyperbolic curve}

In the case of the hyperbolic curve, we could start from the line arrangement
as above and use the three other possible orientations to obtain all the possible
topological types.
Instead of this, we
use the tropical approach, which
gives a more elegant and unified construction method in this case.
In fact, as for the $4$-oval curve, a single smooth tropical curve will suffice
to construct all the possible topological types. The curve is
a {\it honeycomb} curve $H$
of degree $5$ with \emph{all} bounded edges being twisted, i.e.\ 
with the set $T$ given by the set of all bounded edges,
cf. \cite{Speyer}.

\begin{minipage}{\textwidth}
  \centering
  \input{files/pics/honeycomb1.TpX}
\end{minipage}

The following picture shows $\widetilde{C}(H,T) \subset \R^2$ on the left hand side and
$\overline{C}(H,T) \subset \rp^2$ on the left hand side.
%
Note that among the $5$ ``strings'' appearing on the right hand side, the
outer two are glued to form the inner oval, the second and fourth are glued to form
the outer oval and the middle string represents the pseudoline.

\begin{minipage}{\textwidth}
  \centering
  \input{files/pics/honeycomb2.TpX}	 \hspace{10ex}
  \input{files/pics/honeycomb3.TpX}
\end{minipage}

As before, on the right hand side we also draw all the vanishing cycles corresponding to shrinking a twisted edge of $H$
(the  directions of the vanishing cycles are identical to those of the corresponding edges).
However, in contrast to the $4$-oval curve, some pairs of twisted edges share an endpoint.
Hence we cannot realize any subdiagram of the smoothing diagram above automatically.
Nevertheless, we may still shrink any collection of pairwise disjoint twisted edges, and
Theorem \ref{thm:viro_nodal} ensures that the corresponding subdiagram is the smoothing diagram of
a real nodal curve.
This is in fact enough to construct any smoothing diagram with $l=3, h=6$ in table \ref{list1}.
Following the ordering of that table, we present corresponding collections of $6$ pairwise disjoint edges in $H$
(note that the choice of collections is, of course, not unique):

\begin{minipage}{\textwidth}
  \centering
\input{files/pics/dimer1.TpX}
\input{files/pics/dimer2.TpX}
\input{files/pics/dimer3.TpX}
\input{files/pics/dimer4.TpX}
\end{minipage}

\begin{minipage}{\textwidth}
  \centering
\input{files/pics/dimer5.TpX}
\input{files/pics/dimer6.TpX}
\input{files/pics/dimer7.TpX}
\input{files/pics/dimer8.TpX}
\end{minipage}

\begin{minipage}{\textwidth}
  \centering
\input{files/pics/dimer9.TpX}
\input{files/pics/dimer10.TpX}
\input{files/pics/dimer11.TpX}
\input{files/pics/dimer12.TpX}
\end{minipage}

\begin{minipage}{\textwidth}
  \centering
\input{files/pics/dimer13.TpX}
\end{minipage}

\subsection{Elliptic Nodes}

Let us now include elliptic nodes (still assuming that $c=0$).\\[-5ex]

\subsubsection{\texorpdfstring{$M$}{M}-curve}

For $e=1$, we can for example use the following constructions. To shrink a positive oval,
we start with the union of a rational cubic with elliptic node and two lines as depicted below.
Note that this arrangement can be obtained by perturbing the tangent lines of two inflection
points of the cubic.
When choosing the ``upward'' orientations for the lines and the ``downward'' orientation for the
cubic, after smoothing according to the orientations we obtain the smoothing diagram on the right hand side.

\begin{minipage}{\textwidth}
  \centering
  \input{files/pics/cubic2linesB.TpX}	\hspace{0.2\textwidth}
  \input{files/pics/cubic2linesB2.TpX}
\end{minipage}

In order to shrink a negative oval, we use the following union of a rational quartic and a line
as depicted below.
This arrangement can be obtained by perturbing a line passing through the two hyperbolic nodes
of the quartic.
When smoothing this reducible curve according to the displayed orientations we obtain the smoothing
diagram on the right hand side.

\begin{minipage}{\textwidth}
  \centering
  \input{files/pics/quarticlineB.TpX}	\hspace{0.2\textwidth}
  \input{files/pics/quarticlineB2.TpX}
\end{minipage}

Using these two constructions we get all smoothing diagrams from table \ref{list3} with $l=7, e=1$.

For $e=2$, we use again patchworking in the form of Theorem \ref{thm:viro_nodal}.
In the following picture, the left hand side depicts a tropical quartic $Q$
with a single twisted edge. We choose ${\mathcal C}_{\text{\rm edges}}$
to be the single twisted edge and ${\mathcal C}_{\text{\rm triangles}}$ to consist
of the two triangles disjoint from the twisted edge.
The right hand side shows the smoothing diagram of
$C'(Q, T, {\mathcal C}_{\text{\rm edges}}, {\mathcal C}_{\text{\rm triangles}})$.
Theorem \ref{thm:viro_nodal} guarantees that it is the smoothing diagram
of a rational nodal quartic $\overline{A}$ in $\rp^2$.
%

\begin{minipage}{\textwidth} \label{tropconstrquartic1}
  \centering
  \input{files/pics/tropicalquartic.TpX}
  \hspace{10ex}
  \input{files/pics/tropicalquarticreal.TpX}
\end{minipage}

We take the union of this quartic with one of the coordinate lines and smooth
the hyperbolic nodes.
We get the
following two pictures of a smooth curve with vanishing cycles. Choosing the axis $x=0$ or $y=0$
leads to the left hand side picture, the line at infinity gives rise to the right hand side picture.

\begin{minipage}{\textwidth}
  \centering
  \input{files/pics/mcurveD.TpX}
  \input{files/pics/mcurveE.TpX}
\end{minipage}

Using these two curves, we can
realize all cases in table \ref{list3} with $l=7, e=2$.

Finally, let us consider the case $e \geq 3$.
We start with a rational quartic and
then perform a quadratic (Cremona's) transformation.
If the three base points of the quadratic transformation are smooth points on the quartic, we obtain
a curve of degree $2 \cdot 4 - 3 = 5$, that is, a rational quintic. Moreover, if the a line through
two of the base points has no other real intersection with the quartic, the two complex
conjugated intersection points are mapped to an elliptic node in the quadratic transform.
Applying this strategy we can construct all the curves under consideration
as images of the following quartics under the a quadratic transformation which contracts the
depicted lines.

\begin{minipage}{\textwidth}
  \centering
  A \hspace{-3ex} \input{files/pics/quadtrans1.TpX}	\hspace{3ex}
  B \hspace{-3ex} \input{files/pics/quadtrans2.TpX}	\hspace{3ex}
  C \hspace{-3ex} \input{files/pics/quadtrans3.TpX}	\hspace{3ex}
  D \hspace{-3ex} \input{files/pics/quadtrans4.TpX}
\end{minipage}

\begin{minipage}{\textwidth}
  \centering
  E \hspace{-3ex} \input{files/pics/quadtrans5.TpX}	\hspace{3ex}
  F \hspace{-3ex} \input{files/pics/quadtrans6.TpX}	\hspace{3ex}
  G \hspace{-3ex} \input{files/pics/quadtrans7.TpX}	\hspace{3ex}
  H \hspace{-3ex} \input{files/pics/quadtrans8.TpX}
\end{minipage}

\begin{figure}[ht]
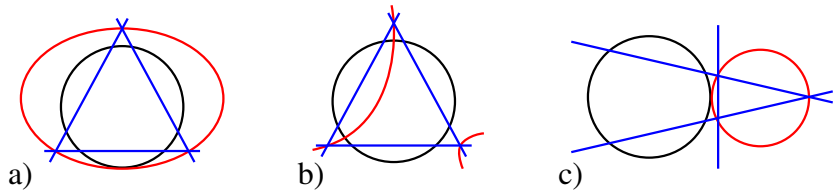

  \centering
  a) \hspace{-3ex} \input{files/pics/conic1.TpX} \hspace{3ex}
  b) \hspace{-3ex} \input{files/pics/conic2.TpX} \hspace{3ex}
  c) \hspace{-3ex} \input{files/pics/conic3.TpX}
	\caption{Chord diagram arrangements for cases A,D and F}
	\label{chorddiagramarrangements1}
\end{figure}

\label{pageexplanations}
It remains to explain why such arrangements of a quartic and three lines exist.
This can easily verified case by case.
For the construction of the quartics, we refer to subsection \ref{sec:degree4}.
The arrangements B and C can be explicitly constructed starting from two ellipses
and three lines intersecting correspondingly.
The arrangements E,G, and H can be constructed by perturbing lines
passing through an elliptic node and the corresponding segment of the curve.
Finally, for arrangements A,D and F, it is helpful to recall from
subsection \ref{sec:degree4} that the quartic curves can be constructed
as quadratic transforms of the corresponding chord diagrams
(see Table \ref{quartic}). Via this quadratic transformation, lines
not passing through the nodes of the quartic correspond to
irreducible conics passing through the three points
of indeterminacy of the quadratic transformation.
Hence the existence of the arrangements can be translated to
the existence of a corresponding arrangement of conics, which can
be proven explicitly. (We refer to such arrangements as
\emph{chord diagram arrangements} in the following.)
For example, for arrangement A it suffices to construct
a line whose only intersection with the quartic is a point of tangency
at a loop. The corresponding chord diagram arrangement can be
constructed explicitly by starting from two tangent ellipses
as depicted in Figure \ref{chorddiagramarrangements1} a).
Arrangement D can be obtained from the chord diagram arrangement
in Figure \ref{chorddiagramarrangements1} b), symmetrized
with respect to rotation around the chord diagram circle by $120^\circ$.
Finally, for arrangement F it is enough to construct a line
whose only intersection with the quartic is a point of tangency
at the middle oval. The two ``horizontal'' lines can then be obtained
by perturbing lines passing through this point of tangency and the
elliptic node. The tangent line can be constructed via the
chord diagram arrangement Figure \ref{chorddiagramarrangements1} c).

\subsubsection{Hyperbolic Curve and non-convex 4-Oval Curve}

Consider the union of a singular cubic with an elliptic node with two lines,
each intersecting the cubic in three points as it is shown on the figure below:

\begin{minipage}{\textwidth}
  \centering
	\input{files/pics/cubic2lines.TpX}
  \hspace{10ex}
  \input{files/pics/cubic2lines1.TpX}
  \input{files/pics/cubic2lines2.TpX}
  \input{files/pics/cubic2lines3.TpX}
\end{minipage}

There are three possibilities (up to symmetries) to orient the components of this reducible quintic
and we get the smoothing diagrams on the right hand side.
Using the first diagram we can realize all smoothing diagrams with $l=3, h=5, e=1$ (cf.\ table \ref{list1}).
The remaining two diagrams give rise to all types with with $l=5, h=5, e=1$ (cf.\ table \ref{list2}).

Let us consider the remaining 4-oval curve cases. For $e=2$, start with the union
of a quartic with 2 elliptic nodes and a line. For $e=3$, we take a quartic with
three elliptic nodes instead. The intersection patterns we need look as follows.

\begin{minipage}{\textwidth}
  \centering
  \input{files/pics/quarticline1.TpX}
  \input{files/pics/quarticline2.TpX}
  \input{files/pics/quarticline3.TpX}
\end{minipage}

Deforming these curves, we can realize all the remaining cases of 4-oval curves with
$e=2$ and $e=3$. The reducible curves can be constructed as follows.
For the construction of the quartics, we again refer to subsection \ref{sec:degree4}.
The first reducible curve was already constructed by tropical methods
on page \pageref{tropconstrquartic1}. The second one is straightforward.
The third one is given by slightly modifying the previous tropical construction
by dropping the twisted edge and including the third triangle instead of it as depicted below.

\begin{minipage}{\textwidth} \label{tropconstrquartic2}
  \centering
  \input{files/pics/tropicalquartic2.TpX}
  \hspace{10ex}
  \input{files/pics/tropicalquarticreal2.TpX}
\end{minipage}

The constructed quartic has three elliptic nodes and its union with with any of
the coordinate lines gives the third reducible curve from above.

\subsection{Complex-conjugated nodes}

Finally, let us consider the case $c>0$.
For $c=2$, all isotopy types in tables \ref{list4} and \ref{list5} except for the last one ($h=0, e=4$)
can be obtained from the union of a rational quartic and a line, intersecting in only two
real points. For such a reducible curve, we have two possibilities to pick orientations,
corresponding to $\sigma  = 0$ or $\sigma = 2$. Furthermore, we can choose which of the two real
intersection points we perturb (while keeping the other one) in order to obtain an irreducible rational quintic.
The following lists show the reducible curves and
the isotopy types realized.
We start with the case $e=0$.

\begin{center}
\begin{longtable}{lll}
  & $\sigma = 0$ & $\sigma = 2$ \\
  A \input{files/pics/redquintic1.TpX} &
	\input{files/pics/cc1c.TpX} &
	\input{files/pics/cc1a.TpX} \\
  B \input{files/pics/redquintic2.TpX} &
	\input{files/pics/cc1f.TpX} \input{files/pics/cc1g.TpX} &
	\input{files/pics/cc1b.TpX} \\
  C \input{files/pics/redquintic3.TpX} &
	\input{files/pics/cc1c.TpX} &
	\input{files/pics/cc1e.TpX} \\
  D \input{files/pics/redquintic4.TpX} &
	\input{files/pics/cc1i.TpX} &
	\input{files/pics/cc1d.TpX} \\
  E \input{files/pics/redquintic5.TpX} &
	\input{files/pics/cc1o.TpX} &
	\input{files/pics/cc1lnew1.TpX} \\
  F \input{files/pics/redquintic6.TpX} &
	\input{files/pics/cc1p.TpX} &
	\input{files/pics/cc1n.TpX} \\
  G \input{files/pics/redquintic7.TpX} &
	\input{files/pics/cc1q.TpX} \input{files/pics/cc1r.TpX} &
	\input{files/pics/cc1lnew.TpX}
\end{longtable}
\end{center}

\begin{figure}[ht]
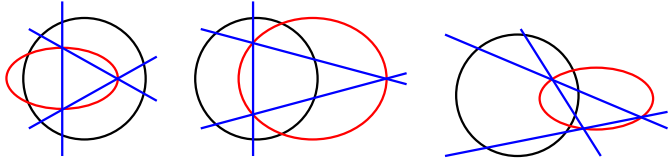

  \centering
  \input{files/pics/conic4.TpX}
  \input{files/pics/conic5.TpX}
  \input{files/pics/conic6.TpX}
	\caption{Chord diagram arrangements for cases E,F and G}
	\label{chorddiagramarrangements2}
\end{figure}

The reducible curves on the left hand side can be constructed as follows.
Again, the construction of the quartics can be found in subsection \ref{sec:degree4}.
The reducible curves in line A,B,C and D were constructed before
at the end of subsection \ref{pageexplanations} (on page
\pageref{pageexplanations}). Moreover, the same methods
apply to the reducible curves E,F and G.
They can be obtained explicitly via quadratic transformation
from the chord diagram arrangements depicted in
Figure \ref{chorddiagramarrangements2}.

The following list deals with the cases $e=1,2,3$.

\begin{center}
\begin{longtable}{lll}
  \label{constructionsccnodes}
  A \input{files/pics/redquintic8.TpX} &
	\input{files/pics/cc1s.TpX} &
	\input{files/pics/cc1u.TpX} \\
  B \input{files/pics/redquintic9.TpX} &
	\input{files/pics/cc1t.TpX} &
	\input{files/pics/cc1uu.TpX} \\
  C \input{files/pics/redquintic10.TpX} &
  \input{files/pics/cc1v.TpX} &
	\input{files/pics/cc1vv.TpX} \\
  D \input{files/pics/redquintic11.TpX} &
	\input{files/pics/cc1w.TpX} &
	\input{files/pics/cc1x.TpX} \\
  E \input{files/pics/redquintic12.TpX} &
	\input{files/pics/cc1y.TpX} &
\end{longtable}
\end{center}

Again, we should explain the constructions of the reducible curves.
The reducible curves A, B, D and E were constructed before
at the end of subsection \ref{pageexplanations} (on page
\pageref{pageexplanations}).
The curve C can be obtained by perturbing a line passing
through the elliptic and a hyperbolic node of the quartic.

The missing case $c=2,e=4$ can be
realized
via quadratic transformation from
a quartic with one elliptic and a pair of complex conjugated nodes.

\begin{minipage}{\textwidth}
  \centering
	\input{files/pics/quadtransnew.TpX}
\end{minipage}

Finally, for $c > 2$, all isotopy types from table \ref{list6} can be easily
realized using
conics and lines (including complex conjugated pairs
of those).
The empty set in the last picture
stands for a pair of complex conjugated
conics without real points. 

\begin{minipage}{\textwidth}
  \centering
  \input{files/pics/c2constra.TpX}
  \input{files/pics/c2constrcc.TpX}
  \input{files/pics/c2constrd.TpX}
  \input{files/pics/c2constre.TpX}
\end{minipage}

\begin{remark} \label{differentsigma}
  As mentioned at the beginning of section \ref{sec:restrictions},
	we actually prove a slightly stronger result than an isotopy classification.
	We classify \emph{complex schemes} of nodal rational curves.
	As explained at loc.\ cit., after choosing one of the two possible
	orientations for the immersed circle, the additional information consists
	of a local orientation of $\rp^2$ at each elliptic node. 
	After smoothing, the induced complex orientation of $C_\circ$ is such that
	each oval obtained from an elliptic node is oriented positively with respect
	to the local orientation. 
	Recalling the list of possible complex schemes for $C_\circ$ from 
	Proposition \ref{prop:complexschemes}, we see that there 
	are two pairs of isotopic arrangements which can be equipped 
	with complex orientations in two ways. They can be distinguished by
	their values of $\sigma$, namely we have the 4-oval curve with 
	$\sigma = 0$ or $2$ and the (unnested) 2-oval curve with $\sigma = 2$ or $4$. 
  Our claim (which implies the classification of complex schemes of nodal rational
  curves) is that whenever an isotopy type can be equipped with two different 
	complex orientations, than both complex schemes are realizable by nodal rational curves.
	This concerns the isotopy types 9, 11, and 13 -- 16 in table \ref{list5}
	as well as types 4 and 6 in table \ref{list6}. 
	Let us explain how to prove this claim. In the table on page \pageref{constructionsccnodes},
	the claim concerns constructions A, D, and E. 
	In cases A and D, the line in the reducible curve can obtained from perturbing a line which
	passes through an elliptic node of the quartic. Depending on which perturbation we choose,
	the elliptic node corresponds to a positive or negative oval in $C_\circ$. 
	Both perturbations give the same isotopy type (as long as we do not change the orientation
	of the line) and hence realize both possible complex schemes in each case. 
	In case E, flipping the orientation of the line obviously switches $\sigma = 0$ and $\sigma =2$.
	In the case $c=2, e=4$, the quadratic transformation of the depicted quartic gives rise
	to $\sigma =2$, whereas using the quartic with elliptic node inside the oval we get $\sigma =0$.
	This follows from the fact that the quadratic transformation does not change the value 
	of $\sigma$ (it can be also checked by hands that, in the first case the elliptic nodes 
	end up in a convex position, while in the second case they are in a non-convex position). 
	For type 4 in table \ref{list6}, we used the union of a real line, a real conic and two complex conjugated
	lines (see second curve in the last picture above). Let us first perturb the real line 
	and the real conic to a nodal cubic such that the elliptic node lies outside of its oval.
	Then, each of the complex conjugated lines intersects one half of the cubic in two points and the other half 
	in one point (otherwise, smoothing a pair of these nodes would give a curve with $\sigma = 0$, 
	which is impossible for this isotopy type). Hence, depending on which pair of complex conjugated 
	nodes we smooth, we obtain a curve with $\sigma=2$ or $\sigma=4$.
	Finally, we realize type 6 in table \ref{list6} from the union a real line and two pairs
	of complex conjugated lines. An even simpler argument as before shows that depending on which two pairs
	of complex conjugated nodes we smooth, we obtain $\sigma=2$ or $\sigma=4$. 
	%
\end{remark}


\printbibliography

\section*{Contact}

\begin{itemize}
	\item 
Ilia Itenberg, Universit\'e Pierre et Marie Curie, 
Institut de Math\'ematiques de Jussieu - Paris Rive Gauche, 
4 place Jussieu, 75252 Paris Cedex 5, France 
and \\
D\'epartement de Math\'emati\-ques et Applications,
Ecole Normale Sup\'erieure,\\
45 rue d'Ulm, 75230 Paris Cedex 5, France;
\href{mailto: ilia.itenberg@imj-prg.fr}{ ilia.itenberg AT imj-prg.fr}. 

  \item
  Grigory Mikhalkin, Section de Math\'ematiques, Universit\'e de Gen\`eve, Battelle Villa,
	1227 Carouge, Suisse; \href{mailto:grigory.mikhalkin@unige.ch}{grigory.mikhalkin AT unige.ch}.

  \item
  Johannes Rau, Fachrichtung Mathematik, Universit\"at der Saarlandes, 
  Postfach 151150, 66041 Saarbr\"ucken, Germany; \href{mailto:johannes.rau@math.uni-sb.de}{johannes.rau AT math.uni-sb.de}.
\end{itemize}

\end {document}